\documentclass[titlepage,11pt]{article}
\usepackage{latexsym}
\usepackage{amsfonts}
\usepackage{amsmath}
\usepackage{tikz}
\usepackage{float}
\usepackage{comment}
\usepackage{xcolor}
\usepackage{enumitem}
\usetikzlibrary{decorations.pathreplacing,decorations.markings}

\newcommand\blackslug{\hbox{\hskip 1pt \vrule width 4pt height 8pt depth 1.5pt
        \hskip 1pt}}
\newcommand\bbox{\hfill \quad \blackslug \bigbreak}
\def\DD{\hbox{-}}
\def\CC{\hbox{-}\cdots\hbox{-}}
\def\LL{,\ldots,}
\def\cupcup{\cup\cdots\cup}
\def\up{\overrightarrow}
\def\arb{arborescence}
\def\sub{subarborescence}
\def\arbs{arborescences}

\title{Graphs with all holes the same length}
\author{Linda Cook\footnotemark[1],
Jake Horsfield\footnotemark[2],
Myriam Preissmann\footnotemark[3],
Cl\'eoph\'ee Robin\footnotemark[3]$^{,}$\footnotemark[5]$^{,}$\footnotemark[7],
\and
Paul Seymour\footnotemark[4],
Ni Luh Dewi Sintiari\footnotemark[5]$^{,}$\footnotemark[8],
Nicolas Trotignon\footnotemark[5], 
Kristina Vu\v{s}kovi\'{c}\footnotemark[2]$^{,}$\footnotemark[6]
}
\footnotetext[1]{
Princeton University, Princeton, NJ 08544, USA. Partially supported by AFOSR grant
A9550-19-1-0187, and partially supported by the Institute for
Basic Science (IBS-R029-C1). Current affiliation: Discrete Mathematics Group, Institute for Basic Science,
Daejeon, Republic of Korea.
}
\footnotetext[2]{
School of Computing, University of Leeds, UK.
}
\footnotetext[3]{
Univ. Grenoble Alpes, CNRS, Grenoble INP, G-SCOP, 38000 Grenoble, France
}
\footnotetext[4]{
Princeton University, Princeton, NJ 08544, USA. Partially supported by AFOSR grants
A9550-19-1-0187 and FA9550-22-1-0234, and by NSF grants  DMS-1800053 and DMS-2154169.
}
\footnotetext[5]{
Univ. Lyon, EnsL, UCBL, CNRS, LIP, F-69342, Lyon Cedex 07, France. Partially
supported by the LABEX MILYON (ANR-10-LABX-0070) of Universit\'e de Lyon, within
the program Investissements d'Avenir (ANR-11-IDEX-0007) operated by the French National Research Agency (ANR) and by Agence Nationale de la Recherche (France) under
research grant ANR DIGRAPHS ANR-19-CE48-0013-01.
}
\footnotetext[6]{
Partially supported by EPSRC grant EP/N019660/1 and EP/V002813/1.
}
\footnotetext[7]{
Current affiliation: Normandie Univ., UNICAEN, ENSICAEN, CNRS, GREYC,
14000 Caen, France
}
\footnotetext[8]{
Current affiliation: Universitas Pendidikan Ganesha, Indonesia.
}

\date{May 8, 2021; revised \today}

\newtheorem{thm}{}[section]

\newcommand{\Proof}{\noindent{\bf Proof.}\ \ }

\begin{document}
\maketitle
\begin{abstract}
A graph is {\em $\ell$-holed} if all its induced cycles of length at least four have length exactly $\ell$.
We give a complete description of the $\ell$-holed graphs for each $\ell\ge 7$.
\end{abstract}

\section{Introduction}
A {\em hole} in a graph is an induced cycle of length at least four (graphs in this paper are finite and have no loops or 
parallel edges). 
Berge graphs and even-hole-free graphs have decomposition theorems that are deep and useful (see~\cite{SPGT,evenhole}), but not fully 
satisfactory,
and we do not know explicit constructions that will generate all the graphs of either type. What if we restrict the lengths of holes
much more severely; can we then give explicit constructions? 
In this direction, the simplest class is the class where holes of all lengths are excluded. This is the class of chordal graphs, 
and it has a structural description that is fully understood~\cite{dirac}. The next simplest is to exclude holes of all lengths
except one, and this is what we study here. As we will see, the description is complex, but when the permitted hole length is at 
least seven, our description fully describes the structure of 
all the graphs in the class.  Incidentally, when all holes have length five, we do not know a complete description, but 
such graphs were studied~\cite{woodruffe} in the context of algebraic combinatorics and commutative algebra.

If $\ell\ge 4$ is an integer, we say a graph is {\em $\ell$-holed} if all its holes have length exactly $\ell$.
How can we make the most general $\ell$-holed graph?
There are a few cases that come to mind immediately: chordal graphs; cycles; pyramids in which 
every path 
from apex to base has the same length; thetas in which every path between the two ends of the theta has the same length; and prisms
in which every path between the two triangles of the prism has the same length. (We will define all these terms later.) 
And we can enlarge these in trivial ways, for instance by overlapping them on clique cutsets, but finding further examples is not so easy.
It turns out that when $\ell\ge 7$, the most general example has a family resemblance to these easy ones.

Before we give the construction, let us digress a little: how can we test if a graph is $\ell$-holed in polynomial time? This is easy
for each given $\ell$, but obtaining an algorithm with running time a polynomial in $|G|$ and independent of $\ell$ is not so clear.
($|G|$ denotes the number of vertices of a graph $G$.)
But here is one simple way. There is an algorithm due to Berger, Seymour and Spirkl~\cite{berger} that, given two vertices $s,t$
of a graph $G$, tests whether there is an induced path between $s,t$ with length more than the distance between $s,t$. 
(Its running time is $\mathcal{O}(|G|^{18})$.) For each three-vertex induced path $a\DD b\DD c$ of $G$, delete $b$ and all its neighbours except 
$a,c$
from $G$, forming $G'$ say, and first test whether there is an $a-c$ path in $G'$ (if not, move on to the next three-vertex path);
check that the distance in $G'$ between $a,c$ is $\ell-2$ (if not, $G$ is not $\ell$-holed and we stop); 
and use the algorithm 
of~\cite{berger} to check that there is no induced path in $G'$ between $a,c$ with length more than $\ell-2$ (if there is such a path,
the graph 
is not $\ell$-holed and we stop). If after processing all three-vertex paths, we still have not determined that $G$ is not 
$\ell$-holed, then it is $\ell$-holed and we stop. This has running time $\mathcal{O}(|G|^{21})$. It can be done faster with 
more complication, making use of the structure theorem proved in this paper, and 
such an algorithm appears in Jake Horsfield's thesis~\cite{jake}.

Let us return to constructing $\ell$-holed graphs.
One way to make a larger $\ell$-holed graph from two smaller ones is via clique cutsets.
We say that $X\subseteq V(G)$ is a {\em clique cutset} of $G$ if $G[X]$ is a complete graph and $G\setminus X$ is disconnected. 
($G[X]$ denotes the subgraph of $G$ induced on $X$, and $G\setminus X$
is the subgraph of $G$ obtained by deleting $X$.) 
If $G_1,G_2$ are two $\ell$-holed graphs, and $X_i$ is a clique of $G_i$
for $i = 1,2$, both of the same cardinality and with $|X_i|<|G_i|$, and we identify each vertex of $X_1$ with a vertex of $X_2$ bijectively, the graph $G$
we produce is $\ell$-holed, and admits a clique cutset; and every graph with a clique cutset can be built by this operation.
Consequently, to understand $\ell$-holed graphs in general, it suffices to understand those
with no clique cutset.

There is another ``trivial'' way to make larger $\ell$-holed graphs from smaller ones; add a new vertex adjacent to all
old vertices. We call a vertex of $G$ adjacent to all the other vertices of $G$ a {\em universal vertex}. 
Thus, we would like
to describe all $\ell$-holed graph with no clique cutset and no universal vertex. We were not able to do this for $\ell=4,5,6$,
but we have a complete description for all $\ell\ge 7$.

There is a third way to enlarge $\ell$-holed graphs to larger ones; choose a vertex and replace it by a set of pairwise
adjacent vertices,
each with the same neighbours as the original vertex (and one another). If two adjacent vertices have exactly the same neighbours,
we call them {\em adjacent twins}. We could assume there are no adjacent twins, without loss of generality, but it makes little 
difference, so usually we will not do so.

Some terminology: two disjoint subsets $X,Y$ of a graph are {\em complete} (to each other) if every vertex of $X$ is adjacent to every vertex in $Y$, and {\em anticomplete}
(to each other) if there are no edges between $X,Y$. Incidentally,  we will sometimes use expressions such as ``$G$-adjacent'' and ``$G$-neighbour'' 
when we want to make clear which graph we are using, as we will often have different graphs with the same vertex set.

An {\em ordering} of a set $X$ means a sequence enumerating the members of $X$. Let $v_1\LL v_n$ be an ordering of $X\subseteq V(G)$.
We say a vertex $u\in V(G)\setminus X$ is {\em adjacent to an initial segment} of the ordering if for all $i,j\in \{1\LL n\}$ with $i<j$, if $u,v_j$ are adjacent then $u,v_i$ are adjacent.
An {\em ordered clique} means a clique together with 
some ordering of it. We will often use the same notation for an ordered clique and the (unordered) clique itself, leaving the 
ordering to be understood when it is needed.
A {\em half-graph} is a bipartite graph with no induced two-edge matching.
Let $X,Y$ be disjoint subsets of $V(G)$; we denote by $G[X,Y]$ the bipartite subgraph of $G$ with vertex set $X\cup Y$
and edge set the set of edges of $G$ between $X,Y$. 
Take orderings $x_1\LL x_m$ and $y_1\LL y_n$ of $X$ and $Y$ respectively. We say $G[X,Y]$ {\em obeys} these orderings if 
for all $i,i',j,j'$ with $1\le i\le i'\le m$ and $1\le j\le j'\le n$, if $x_{i'}y_{j'}$ is an edge
then $x_{i}y_{j}$ is an edge; or, equivalently, each vertex in $Y$ is adjacent to an initial segment of $x_1\LL x_m$,
and each vertex in $X$ is adjacent to an initial segment of $y_1\LL y_n$. 
Thus, $G[X,Y]$ is a half-graph if and only if there are orderings of $X$ and $Y$
that $G[X,Y]$ obeys.

Half-graphs will be of great importance in this paper. For instance, let $G$ be an $\ell$-holed graph where $\ell\ge 5$, 
and let $X,Y$
be disjoint cliques of $G$. Since there is no 4-hole 
with two vertices in $X$ and two in $Y$ (a {\em $k$-hole} means a hole of length $k$), it follows that
$G[X,Y]$ is a half-graph. If $X,Y,Z$ are disjoint cliques of $G$, we say that $G[X,Y], G[X,Z]$ are {\em compatible}
if $G[X,Y\cup Z]$ is a half-graph. If $Y,Z$ are anticomplete, then $G[X,Y], G[X,Z]$ are compatible
if and only if there is no induced four-vertex path in $G[X\cup Y\cup Z]$ with first vertex in $Y$, second and third in $X$, and 
fourth in $Z$; this latter property is important for keeping all holes the same length.
This is also equivalent to asking that there are orderings of $X,Y$ and $Z$ that $G[X,Y]$ and $G[X,Z]$ both obey.

\begin{figure}[H]
\centering

\begin{tikzpicture}[scale=.8,auto=left]
\tikzstyle{every node}=[inner sep=1.5pt, fill=black,circle,draw]
\def\q{1.5}
\def\r{2}
\def\s{2.5}
\def\t{3}
\def\angle{360/7}
\node (a1) at ({\q*cos(90)}, {\q*sin(90)}) {};
\node (a2) at ({\r*cos(90)}, {\r*sin(90)}) {};
\node (a3) at ({\s*cos(90)}, {\s*sin(90)}) {};
\node (b1) at ({\q*cos(90+\angle)}, {\q*sin(90+\angle)}) {};
\node (b2) at ({\r*cos(90+\angle)}, {\r*sin(90+\angle)}) {};
\node (b3) at ({\s*cos(90+\angle)}, {\s*sin(90+\angle)}) {};
\node (b4) at ({\t*cos(90+\angle)}, {\t*sin(90+\angle)}) {};
\node (c1) at ({\q*cos(90+2*\angle)}, {\q*sin(90+2*\angle)}) {};
\node (c2) at ({\r*cos(90+2*\angle)}, {\r*sin(90+2*\angle)}) {};
\node (c3) at ({\s*cos(90+2*\angle)}, {\s*sin(90+2*\angle)}) {};
\node (d1) at ({\q*cos(90+3*\angle)}, {\q*sin(90+3*\angle)}) {};
\node (d2) at ({\r*cos(90+3*\angle)}, {\r*sin(90+3*\angle)}) {};
\node (e1) at ({\q*cos(90+4*\angle)}, {\q*sin(90+4*\angle)}) {};
\node (e2) at ({\r*cos(90+4*\angle)}, {\r*sin(90+4*\angle)}) {};
\node (f1) at ({\q*cos(90+5*\angle)}, {\q*sin(90+5*\angle)}) {};
\node (f2) at ({\r*cos(90+5*\angle)}, {\r*sin(90+5*\angle)}) {};
\node (f3) at ({\s*cos(90+5*\angle)}, {\s*sin(90+5*\angle)}) {};
\node (g1) at ({\q*cos(90+6*\angle)}, {\q*sin(90+6*\angle)}) {};
\node (g2) at ({\r*cos(90+6*\angle)}, {\r*sin(90+6*\angle)}) {};
\node (g3) at ({\s*cos(90+6*\angle)}, {\s*sin(90+6*\angle)}) {};
\node (g4) at ({\t*cos(90+6*\angle)}, {\t*sin(90+6*\angle)}) {};

\foreach \from/\to in {a1/b1,a1/b2,a1/b3,a1/b4,a2/b1,a2/b2,a2/b3,a3/b1,a3/b2}
\draw [-] (\from) -- (\to);

\foreach \from/\to in {b1/c1,b1/c2,b1/c3,b2/c1,b2/c2,b3/c1,b3/c2,b4/c1}
\draw [-] (\from) -- (\to);

\foreach \from/\to in {c1/d1,c1/d2,c2/d1,c3/d1}
\draw [-] (\from) -- (\to);

\foreach \from/\to in {d1/e1,d1/e2,d2/e1}
\draw [-] (\from) -- (\to);

\foreach \from/\to in {e1/f1,e1/f2,e1/f3,e2/f1,e2/f2}
\draw [-] (\from) -- (\to);

\foreach \from/\to in {f1/g1,f1/g2,f1/g3,f1/g4,f2/g1,f2/g2,f2/g3,f3/g1}
\draw [-] (\from) -- (\to);

\foreach \from/\to in {g1/a1,g1/a2,g1/a3,g2/a1,g2/a2,g3/a1,g4/a1}
\draw [-] (\from) -- (\to);

\foreach \from/\to in {a1/a2,a2/a3,b1/b2,b2/b3,b3/b4,c1/c2,c2/c3,d1/d2,e1/e2,f1/f2,f2/f3,g1/g2,g2/g3,g3/g4}
\draw [-] (\from) -- (\to);

\foreach \from/\to in {a1/a3,b1/b3,b1/b4,b2/b4,c1/c3,f1/f3,g1/g3,g1/g4,g2/g4}
\draw [bend left =30] (\from) to (\to);

\end{tikzpicture}

\caption{A blow-up of a 7-cycle.} \label{fig:cycleblowup}
\end{figure}
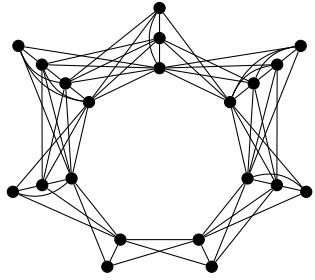

Let $G$ be a graph with vertex set partitioned into sets $W_1\LL W_{\ell}$, with the following properties:
\begin{itemize}
\item $W_1\LL W_{\ell}$ are non-null cliques;
\item for $1\le i\le \ell$, $G[W_{i-1}, W_{i}]$ is a half-graph (reading subscripts modulo $\ell$);
\item for all distinct $i,j\in \{1\LL \ell\}$, if there is an edge between $W_i, W_j$ then $j=i\pm 1$ (modulo $\ell$); and
\item for $1\le i\le \ell$, the graphs $G[W_i, W_{i+1}], G[W_i, W_{i-1}]$ are compatible.
\end{itemize}
(See figure \ref{fig:cycleblowup}.) We call such a graph a {\em blow-up of an $\ell$-cycle}.

To describe our main result we will need \arbs{}.
An {\em arborescence} is a tree with its edges directed in such a way that no two edges have a common head; or equivalently, such that
for some vertex $r(T)$ (called the {\em apex}), every edge is directed away from $r(T)$. A {\em leaf} is a vertex different from the 
apex, with outdegree zero, and
$L(T)$ denotes the set of leaves of the \arb{} $T$.

Our theorem says that if $G$ is $\ell$-holed, and $\ell\ge 7$, and $G$ has no clique cutset or 
universal vertex, then either $G$ is a blow-up of an $\ell$-cycle, or $G$ is an instance of a construction we will describe.
The construction is rather complicated, however, and we will give the description in stages. The cases of $\ell$ odd and $\ell$ even 
are different, and the case when $\ell$ is odd is simpler, so 
let us begin with that. The underlying structure is what we call an {\em $\ell$-framework}, and is best described with a figure.

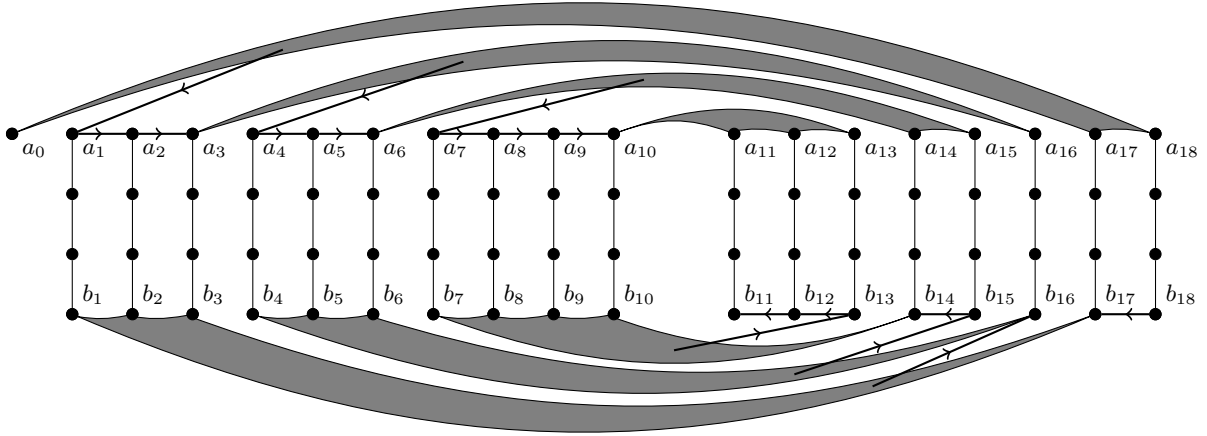
\begin{figure}[h]
\centering
\begin{tikzpicture}[scale=.8,auto=left]
\tikzstyle{every node}=[inner sep=1.5pt, fill=black,circle,draw]

\node (a0) at (0,3) {};
\node (a1) at (1,3) {};
\node (a2) at (2,3) {};
\node (a3) at (3,3) {};
\node (a4) at (4,3) {};
\node (a5) at (5,3) {};
\node (a6) at (6,3) {};
\node (a7) at (7,3) {};
\node (a8) at (8,3) {};
\node (a9) at (9,3) {};
\node (a10) at (10,3) {};
\node (a11) at (12,3) {};
\node (a12) at (13,3) {};
\node (a13) at (14,3) {};
\node (a14) at (15,3) {};
\node (a15) at (16,3) {};
\node (a16) at (17,3) {};
\node (a17) at (18,3) {};
\node (a18) at (19,3) {};
\node (b1) at (1,0) {};
\node (b2) at (2,0) {};
\node (b3) at (3,0) {};
\node (b4) at (4,0) {};
\node (b5) at (5,0) {};
\node (b6) at (6,0) {};
\node (b7) at (7,0) {};
\node (b8) at (8,0) {};
\node (b9) at (9,0) {};
\node (b10) at (10,0) {};
\node (b11) at (12,0) {};
\node (b12) at (13,0) {};
\node (b13) at (14,0) {};
\node (b14) at (15,0) {};
\node (b15) at (16,0) {};
\node (b16) at (17,0) {};
\node (b17) at (18,0) {};
\node (b18) at (19,0) {};
\node (c1) at (1,2) {};
\node (c2) at (2,2) {};
\node (c3) at (3,2) {};
\node (c4) at (4,2) {};
\node (c5) at (5,2) {};
\node (c6) at (6,2) {};
\node (c7) at (7,2) {};
\node (c8) at (8,2) {};
\node (c9) at (9,2) {};
\node (c10) at (10,2) {};
\node (c11) at (12,2) {};
\node (c12) at (13,2) {};
\node (c13) at (14,2) {};
\node (c14) at (15,2) {};
\node (c15) at (16,2) {};
\node (c16) at (17,2) {};
\node (c17) at (18,2) {};
\node (c18) at (19,2) {};
\node (d1) at (1,1) {};
\node (d2) at (2,1) {};
\node (d3) at (3,1) {};
\node (d4) at (4,1) {};
\node (d5) at (5,1) {};
\node (d6) at (6,1) {};
\node (d7) at (7,1) {};
\node (d8) at (8,1) {};
\node (d9) at (9,1) {};
\node (d10) at (10,1) {};
\node (d11) at (12,1) {};
\node (d12) at (13,1) {};
\node (d13) at (14,1) {};
\node (d14) at (15,1) {};
\node (d15) at (16,1) {};
\node (d16) at (17,1) {};
\node (d17) at (18,1) {};
\node (d18) at (19,1) {};

\foreach \from/\to in {a1/c1,a2/c2,a3/c3,a4/c4,a5/c5,a6/c6,a7/c7,a8/c8,a9/c9,a10/c10,a11/c11,a12/c12,a13/c13,a14/c14,a15/c15,a16/c16,
a17/c17,a18/c18}
\draw(\from) to (\to);
\foreach \from/\to in {c1/d1,c2/d2,c3/d3,c4/d4,c5/d5,c6/d6,c7/d7,c8/d8,c9/d9,c10/d10,c11/d11,c12/d12,c13/d13,c14/d14,c15/d15,c16/d16,
c17/d17,c18/d18}
\draw(\from) to (\to);
\foreach \from/\to in {d1/b1,d2/b2,d3/b3,d4/b4,d5/b5,d6/b6,d7/b7,d8/b8,d9/b9,d10/b10,d11/b11,d12/b12,d13/b13,d14/b14,d15/b15,d16/b16,
d17/b17,d18/b18}
\draw(\from) to (\to);

\draw[fill=gray] (a0) to [bend left=20] (a17) to [bend left=10] (a18) to [bend right=23] (a0);
\draw[fill=gray] (a3) to [bend left=22] (a16) to [bend right=17] (a3);
\draw[fill=gray] (a6) to [bend left=17] (a14) to [bend left=10] (a15) to [bend right=20] (a6);
\draw[fill=gray] (a10) to [bend left=20] (a11) to [bend left=10] (a12) to [bend left=10] (a13) to [bend right=20] (a10);

\draw[fill=gray] (b17) to [bend left=23] (b1) to [bend right=10] (b2) to [bend right=10] (b3) to [bend right=20] (b17);
\draw[fill=gray] (b16) to [bend left=20] (b4) to [bend right=10] (b5) to [bend right=10] (b6) to [bend right=17] (b16);
\draw[fill=gray] (b14) to [bend left=20] (b7) to [bend right=10] (b8) to [bend right=10] (b9)  to [bend right=10] (b10) to [bend right=20] (b14);

\node (a0) at (0,3) {};
\node (a1) at (1,3) {};
\node (a2) at (2,3) {};
\node (a3) at (3,3) {};
\node (a4) at (4,3) {};
\node (a5) at (5,3) {};
\node (a6) at (6,3) {};
\node (a7) at (7,3) {};
\node (a8) at (8,3) {};
\node (a9) at (9,3) {};
\node (a10) at (10,3) {};
\node (a11) at (12,3) {};
\node (a12) at (13,3) {};
\node (a13) at (14,3) {};
\node (a14) at (15,3) {};
\node (a15) at (16,3) {};
\node (a16) at (17,3) {};
\node (a17) at (18,3) {};
\node (a18) at (19,3) {};
\node (b1) at (1,0) {};
\node (b2) at (2,0) {};
\node (b3) at (3,0) {};
\node (b4) at (4,0) {};
\node (b5) at (5,0) {};
\node (b6) at (6,0) {};
\node (b7) at (7,0) {};
\node (b8) at (8,0) {};
\node (b9) at (9,0) {};
\node (b10) at (10,0) {};
\node (b11) at (12,0) {};
\node (b12) at (13,0) {};
\node (b13) at (14,0) {};
\node (b14) at (15,0) {};
\node (b15) at (16,0) {};
\node (b16) at (17,0) {};
\node (b17) at (18,0) {};
\node (b18) at (19,0) {};

\tikzstyle{every node}=[]
\draw (a0) node [below right]           {\footnotesize$a_0$};
\draw (a1) node [below right]           {\footnotesize$a_1$};
\draw (a2) node [below right]           {\footnotesize$a_2$};
\draw (a3) node [below right]           {\footnotesize$a_3$};
\draw (a4) node [below right ]           {\footnotesize$a_4$};
\draw (a5) node [below right ]           {\footnotesize$a_5$};
\draw (a6) node [below right]           {\footnotesize$a_6$};
\draw (a7) node [below right]           {\footnotesize$a_7$};
\draw (a8) node [below right]           {\footnotesize$a_8$};
\draw (a9) node [below right]           {\footnotesize$a_9$};
\draw (a10) node [below right]           {\footnotesize$a_{10}$};
\draw (a11) node [below right]           {\footnotesize$a_{11}$};
\draw (a12) node [below right]           {\footnotesize$a_{12}$};
\draw (a13) node [below right]           {\footnotesize$a_{13}$};
\draw (a14) node [below right]           {\footnotesize$a_{14}$};
\draw (a15) node [below right]           {\footnotesize$a_{15}$};
\draw (a16) node [below right]           {\footnotesize$a_{16}$};
\draw (a17) node [below right]           {\footnotesize$a_{17}$};
\draw (a18) node [below right]           {\footnotesize$a_{18}$};

\draw (b1) node [above right]           {\footnotesize$b_1$};
\draw (b2) node [above right]           {\footnotesize$b_2$};
\draw (b3) node [above right]           {\footnotesize$b_3$};
\draw (b4) node [above right]           {\footnotesize$b_4$};
\draw (b5) node [above right]           {\footnotesize$b_5$};
\draw (b6) node [above right]           {\footnotesize$b_6$};
\draw (b7) node [above right]           {\footnotesize$b_7$};
\draw (b8) node [above right]           {\footnotesize$b_8$};
\draw (b9) node [above right]           {\footnotesize$b_9$};
\draw (b10) node [above right]           {\footnotesize$b_{10}$};
\draw (b11) node [above right]           {\footnotesize$b_{11}$};
\draw (b12) node [above right]           {\footnotesize$b_{12}$};
\draw (b13) node [above right]           {\footnotesize$b_{13}$};
\draw (b14) node [above right]           {\footnotesize$b_{14}$};
\draw (b15) node [above right]           {\footnotesize$b_{15}$};
\draw (b16) node [above right]           {\footnotesize$b_{16}$};
\draw (b17) node [above right]           {\footnotesize$b_{17}$};
\draw (b18) node [above right]           {\footnotesize$b_{18}$};

\begin{scope}[thick, decoration={
    markings,
    mark=at position 0.5 with {\arrow{>}}}
    ]
    \draw[postaction={decorate}] (4.5,4.4) to (a1);
    \draw[postaction={decorate}] (7.5,4.2) to (a4);
    \draw[postaction={decorate}] (10.5,3.9) to (a7);
    \draw[postaction={decorate}] (14.3,-1.2) to (b16);
    \draw[postaction={decorate}] (13,-1.0) to (b15);
    \draw[postaction={decorate}] (11,-.6) to (b13);

\foreach \from/\to in {a1/a2,a2/a3,a4/a5,a5/a6,a7/a8,a8/a9,a9/a10, b18/b17, b15/b14,b13/b12,b12/b11}
\draw[postaction={decorate}] (\from) to (\to);
\end{scope}
\end{tikzpicture}
\caption{An 18-bar 9-framework.} \label{fig:oddframework}
\end{figure}

Let us describe the important features of figure \ref{fig:oddframework}. 
There are 19 vertices $a_0\LL a_{18}$ and 18 vertices 
$b_1\LL b_{18}$ (these could be any two numbers $k+1$ and $k$). 
For $1\le i\le k$ there is a vertical path $P_i$ of length $(\ell-3)/2$
between $a_i, b_i$. (In the case of the figure, $\ell=9$.) The numbers $0\LL k$ break into two intervals $\{0\LL m\}$ and $\{m+1\LL k\}$
(in the 
figure $m=10$).

Let us call the grey shaded areas ``tents''. 
The tents are disjoint subsets of the plane, and 
each of the (four, in the figure) upper tents contains one vertex in $\{a_0\LL a_m\}$ called its ``apex'', and 
contains 
a nonempty interval of $\{a_{m+1}\LL a_k\}$ called its ``base''. 
Each of $a_{m+1}\LL a_k$ belongs to the base of an upper tent.
The lower tents do the same with left and right switched.
There can be any positive number of tents, but there must be a tent with apex $a_0$. (There is an odd number of tents in the figure, but
there could be an even number.) Possibly $m=0$, and if so there are no lower tents.
The way the upper and lower tents interleave is important; for each upper tent (except the innermost when there is an odd number of tents), 
the leftmost vertex of its base is some $a_i$, and 
$b_i$ is the apex of one of the lower tents; and for each lower tent (except the innermost when there is an even number of tents), 
the rightmost vertex of its base corresponds to the apex for 
one of the upper tents. (This gives a sort of spiral running through all the apexes, in the figure with vertices
$$a_0\DD a_{17}\DD b_{17}\DD b_3\DD a_3\DD a_{16}\DD b_{16}\DD b_6\DD a_6\DD a_{14}\DD b_{14} \DD b_{10}\DD a_{10},$$
which might be helpful for visualization.)

Each tent is meant to be an \arb{} with the given apex and with set of leaves
the base of the tent, and with its other vertices not drawn. (We call such an \arb{} a {\em tent-\arb{}}.)
For each $i\in \{1\LL m\}$, if $a_{i-1}$ is the apex of an upper tent-\arb{} $T_{i-1}$ say, there is a directed edge
from some nonleaf vertex of $T_{i-1}$ (possibly from $a_{i-1}$) to
$a_{i}$; and if $a_{i-1}$ is not the apex of a tent, there is a directed edge from $a_{i-1}$ to $a_i$. So all these upper tent-\arbs{} and all the vertices $a_0\LL a_m$, are connected up in a sequence to form one big \arb{} $T$ with apex $a_0$, and with 
set of leaves either $\{a_{m+1}\LL a_k\}$ or $\{a_m\LL a_k\}$.
There is a directed path of $T$ that contains $a_0, a_1\LL a_m$ in order, possibly containing 
other vertices of $T$ between them.
Similarly for each $i\in \{m+1\LL k-1\}$, if $b_{i+1}$ is an apex of a lower tent-\arb{} $S_{i+1}$,
there is a directed edge
from some nonleaf vertex of $S_{i+1}$ to
$b_{i}$, and otherwise there is a directed edge $b_{i+1}b_i$.
So similarly the lower tent-\arbs{}, and the vertices $b_{m+1}\LL b_k$, are joined up to make one \arb{} $S$ with apex $b_k$ and with set of leaves either $\{b_1\LL b_m\}$ or $\{b_1\LL b_{m+1}\}$

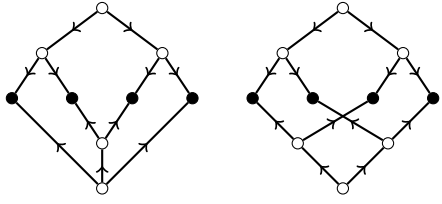
\begin{figure}[H]
\centering

\begin{tikzpicture}[scale=.8,auto=left]
\tikzstyle{every node}=[inner sep=1.5pt, fill=black,circle,draw]
\node (4) at (-1.5,0) {};
\node (5) at (-.5,0) {};
\node (6) at (.5,0) {};
\node (7) at (1.5,0) {};
\tikzstyle{every node}=[inner sep=1.5pt, fill=white,circle,draw]
\node (1) at (0,1.5) {};
\node (2) at (-1,.75) {};
\node (3) at (1,.75) {};
\node (8) at (0,-.75) {};
\node (9) at (0,-1.5) {};

\begin{scope}[thick, decoration={
    markings,
    mark=at position 0.5 with {\arrow{>}}}
    ]
\foreach \from/\to in {1/2, 1/3,2/4,2/5,3/6,3/7,9/8,9/4,9/7,8/5,8/6}
\draw [postaction={decorate}] (\from) to (\to);
\end{scope}
\tikzstyle{every node}=[inner sep=1.5pt, fill=black,circle,draw]
\node (4) at (2.5,0) {};
\node (5) at (3.5,0) {};
\node (6) at (4.5,0) {};
\node (7) at (5.5,0) {};
\tikzstyle{every node}=[inner sep=1.5pt, fill=white,circle,draw]
\node (1) at (4,1.5) {};
\node (2) at (3,.75) {};
\node (3) at (5,.75) {};
\node (8) at (3.25,-.75) {};
\node (9) at (4.75,-.75) {};
\node (10) at (4,-1.5) {};

\begin{scope}[thick, decoration={
    markings,
    mark=at position 0.5 with {\arrow{>}}}
    ]
\foreach \from/\to in {1/2, 1/3,2/4,2/5,3/6,3/7,8/4, 8/6,9/5,9/7,10/8,10/9}
\draw [postaction={decorate}] (\from) to (\to);
\end{scope}
\end{tikzpicture}
\caption{The first pair are coarboreal, the second pair are not.} \label{fig:coarb}
\end{figure}

Thus the figure describes a graph in which some of the edges are directed: each directed edge belongs to one of two \arbs{} $T,S$
and each undirected edge belongs to one of the paths $P_i$. We call such a graph an {\em $\ell$-framework}. (We will explain later
how $\ell$-frameworks describe the structure of $\ell$-holed graphs.)
We observe:
\begin{thm}\label{oddadj}
	Let $\ell\ge 5$ be odd and let $F$ be an $\ell$-framework, with notation as above. For $1\le i<j\le k$, either there
	is a directed path of $T$ between $a_i, a_j$, or there is a directed path of $S$ between $b_i, b_j$ and not both.
\end{thm}
(To clarify: ``directed path of $T$ between $a_i, a_j$'' means a directed path either from $a_i$ to $a_j$, or from $a_j$ to $a_i$.)
Next we will describe a similar object for when $\ell$ is even, but we need another concept.
Let $T$ be an arborescence.
For $v\in V(T)$, let $D_v$ be the set  of all vertices $w\in L(T)$ for which there is a directed path of $T$ from $v$ to $w$.
Let $S$ be a tree with $V(S)=L(T)$. We say that $T$ {\em lives in $S$} if for each $v\in V(T)$, the set $D_v$
is the vertex set of a subtree of $S$. Let $T,T'$ be arborescences with $L(T)=L(T')$.
We say they are {\em coarboreal} if there is a tree $S$ with $V(S)=L(T)=L(T')$ such that $T,T'$ both live in $S$.
For instance, the first pair of arborescences in figure \ref{fig:coarb} (with leaf set the four black vertices) are coarboreal, but the second pair are not.
Finally, let $T,T'$ be \arbs{} and let $\phi$ be a bijection from $L(T)$ onto $L(T')$. We say that $T,T'$ are
{\em coarboreal under $\phi$} if identifying each vertex of $L(T)$ with its image under $\phi$ gives a coarboreal pair.

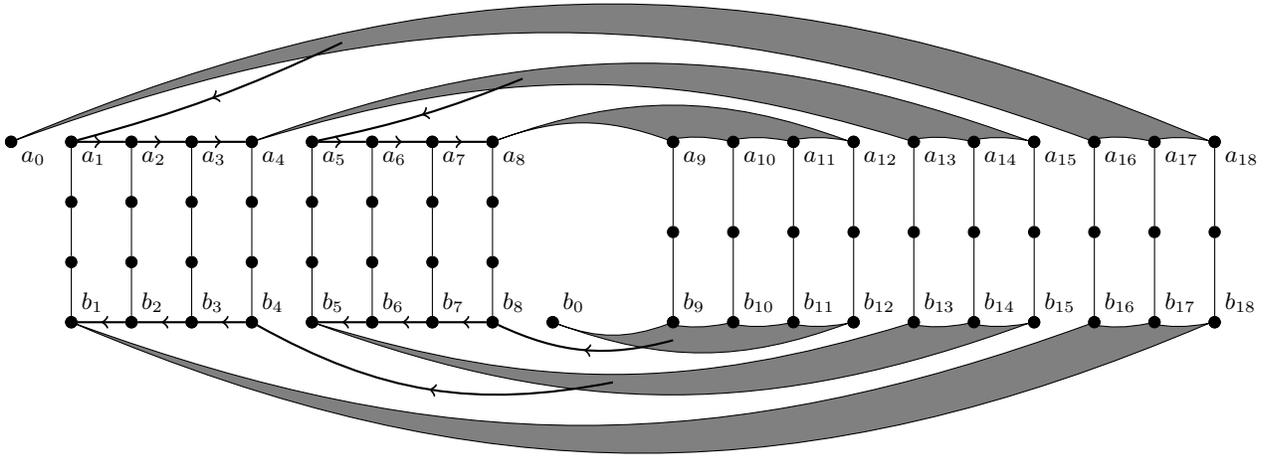
\begin{figure}[h]
\centering
\begin{tikzpicture}[scale=.8,auto=left]
\tikzstyle{every node}=[inner sep=1.5pt, fill=black,circle,draw]

\node (a0) at (-1,3) {};
\node (a1) at (0,3) {};
\node (a2) at (1,3) {};
\node (a3) at (2,3) {};
\node (a4) at (3,3) {};
\node (a5) at (4,3) {};
\node (a6) at (5,3) {};
\node (a7) at (6,3) {};
\node (a8) at (7,3) {};

\node (a9) at (10,3) {};
\node (a10) at (11,3) {};
\node (a11) at (12,3) {};
\node (a12) at (13,3) {};
\node (a13) at (14,3) {};
\node (a14) at (15,3) {};
\node (a15) at (16,3) {};
\node (a16) at (17,3) {};
\node (a17) at (18,3) {};
\node (a18) at (19,3) {};

\node (b1) at (0,0) {};
\node (b2) at (1,0) {};
\node (b3) at (2,0) {};
\node (b4) at (3,0) {};
\node (b5) at (4,0) {};
\node (b6) at (5,0) {};
\node (b7) at (6,0) {};
\node (b8) at (7,0) {};

\node (b9) at (10,0) {};
\node (b10) at (11,0) {};
\node (b11) at (12,0) {};
\node (b12) at (13,0) {};
\node (b13) at (14,0) {};
\node (b14) at (15,0) {};
\node (b15) at (16,0) {};
\node (b16) at (17,0) {};
\node (b17) at (18,0) {};
\node (b18) at (19,0) {};

\node (c1) at (0,2) {};
\node (c2) at (1,2) {};
\node (c3) at (2,2) {};
\node (c4) at (3,2) {};
\node (c5) at (4,2) {};
\node (c6) at (5,2) {};
\node (c7) at (6,2) {};
\node (c8) at (7,2) {};

\node (c9) at (10,1.5) {};
\node (c10) at (11,1.5) {};
\node (c11) at (12,1.5) {};
\node (c12) at (13,1.5) {};
\node (c13) at (14,1.5) {};
\node (c14) at (15,1.5) {};
\node (c15) at (16,1.5) {};
\node (c16) at (17,1.5) {};
\node (c17) at (18,1.5) {};
\node (c18) at (19,1.5) {};

\node (d1) at (0,1) {};
\node (d2) at (1,1) {};
\node (d3) at (2,1) {};
\node (d4) at (3,1) {};
\node (d5) at (4,1) {};
\node (d6) at (5,1) {};
\node (d7) at (6,1) {};
\node (d8) at (7,1) {};

\node (b0) at (8,0) {};

\foreach \from/\to in {a1/c1,a2/c2,a3/c3,a4/c4,a5/c5,a6/c6,a7/c7,a8/c8,a9/c9,a10/c10,a11/c11,a12/c12,a13/c13,a14/c14,a15/c15,a16/c16,
a17/c17,a18/c18}
\draw(\from) to (\to);
\foreach \from/\to in {c1/d1,c2/d2,c3/d3,c4/d4,c5/d5,c6/d6,c7/d7,c8/d8}
\draw(\from) to (\to);
\foreach \from/\to in {d1/b1,d2/b2,d3/b3,d4/b4,d5/b5,d6/b6,d7/b7,d8/b8,c9/b9,c10/b10,c11/b11,c12/b12,c13/b13,c14/b14,c15/b15,c16/b16,
c17/b17,c18/b18}
\draw(\from) to (\to);

\draw[fill=gray] (a0) to [bend left=20] (a16) to [bend left=10] (a17)  to [bend left=10] (a18) to [bend right=23] (a0);
\draw[fill=gray] (a4) to [bend left=17] (a13) to [bend left=10] (a14) to [bend left=10] (a15) to [bend right=20] (a4);
\draw[fill=gray] (a8) to [bend left=20] (a9) to [bend left=10] (a10) to [bend left=10] (a11) to [bend left=10] (a12) to [bend right=20] (a8);

\draw[fill=gray] (b1) to [bend right=20] (b16) to [bend right=10] (b17) to [bend right=10] (b18) to [bend left=23] (b1);
\draw[fill=gray] (b5) to [bend right=17] (b13) to [bend right=10] (b14) to [bend right=10] (b15) to [bend left=20] (b5);
\draw[fill=gray] (b0) to [bend right=20] (b9) to [bend right=10] (b10) to [bend right=10] (b11)  to [bend right=10] (b12) to [bend left=20] (b0);

\node (a0) at (-1,3) {};
\node (a1) at (0,3) {};
\node (a2) at (1,3) {};
\node (a3) at (2,3) {};

\node (a4) at (3,3) {};
\node (a5) at (4,3) {};
\node (a6) at (5,3) {};
\node (a7) at (6,3) {};
\node (a8) at (7,3) {};
\node (a9) at (10,3) {};
\node (a10) at (11,3) {};
\node (a11) at (12,3) {};
\node (a12) at (13,3) {};
\node (a13) at (14,3) {};
\node (a14) at (15,3) {};
\node (a15) at (16,3) {};
\node (a16) at (17,3) {};
\node (a17) at (18,3) {};
\node (a18) at (19,3) {};
\node (b1) at (0,0) {};
\node (b2) at (1,0) {};
\node (b3) at (2,0) {};
\node (b4) at (3,0) {};
\node (b5) at (4,0) {};
\node (b6) at (5,0) {};
\node (b7) at (6,0) {};
\node (b8) at (7,0) {};
\node (b9) at (10,0) {};
\node (b10) at (11,0) {};
\node (b11) at (12,0) {};
\node (b12) at (13,0) {};
\node (b13) at (14,0) {};
\node (b14) at (15,0) {};
\node (b15) at (16,0) {};
\node (b16) at (17,0) {};
\node (b17) at (18,0) {};
\node (b18) at (19,0) {};
\node (b0) at (8,0) {};

\tikzstyle{every node}=[]
\draw (a0) node [below right]           {\footnotesize$a_0$};
\draw (a1) node [below right]           {\footnotesize$a_1$};
\draw (a2) node [below right]           {\footnotesize$a_2$};
\draw (a3) node [below right]           {\footnotesize$a_3$};
\draw (a4) node [below right ]           {\footnotesize$a_4$};
\draw (a5) node [below right ]           {\footnotesize$a_5$};
\draw (a6) node [below right]           {\footnotesize$a_6$};
\draw (a7) node [below right]           {\footnotesize$a_7$};
\draw (a8) node [below right]           {\footnotesize$a_8$};
\draw (a9) node [below right]           {\footnotesize$a_9$};
\draw (a10) node [below right]           {\footnotesize$a_{10}$};
\draw (a11) node [below right]           {\footnotesize$a_{11}$};
\draw (a12) node [below right]           {\footnotesize$a_{12}$};
\draw (a13) node [below right]           {\footnotesize$a_{13}$};
\draw (a14) node [below right]           {\footnotesize$a_{14}$};
\draw (a15) node [below right]           {\footnotesize$a_{15}$};
\draw (a16) node [below right]           {\footnotesize$a_{16}$};

\draw (a17) node [below right]           {\footnotesize$a_{17}$};
\draw (a18) node [below right]           {\footnotesize$a_{18}$};

\draw (b1) node [above right]           {\footnotesize$b_1$};
\draw (b2) node [above right]           {\footnotesize$b_2$};
\draw (b3) node [above right]           {\footnotesize$b_3$};
\draw (b4) node [above right]           {\footnotesize$b_4$};
\draw (b5) node [above right]           {\footnotesize$b_5$};
\draw (b6) node [above right]           {\footnotesize$b_6$};
\draw (b7) node [above right]           {\footnotesize$b_7$};
\draw (b8) node [above right]           {\footnotesize$b_8$};
\draw (b9) node [above right]           {\footnotesize$b_9$};
\draw (b10) node [above right]           {\footnotesize$b_{10}$};
\draw (b11) node [above right]           {\footnotesize$b_{11}$};
\draw (b12) node [above right]           {\footnotesize$b_{12}$};
\draw (b13) node [above right]           {\footnotesize$b_{13}$};
\draw (b14) node [above right]           {\footnotesize$b_{14}$};
\draw (b15) node [above right]           {\footnotesize$b_{15}$};
\draw (b16) node [above right]           {\footnotesize$b_{16}$};
\draw (b17) node [above right]           {\footnotesize$b_{17}$};
\draw (b18) node [above right]           {\footnotesize$b_{18}$};
\draw (b0) node [above right]           {\footnotesize$b_{0}$};

\begin{scope}[thick, decoration={
    markings,
    mark=at position 0.5 with {\arrow{>}}}
    ]
    \draw[bend left=5,postaction={decorate}] (4.5,4.65) to (a1);
    \draw[bend left=5, postaction={decorate}] (7.5,4.05) to (a5);
    \draw[bend left=20, postaction={decorate}] (9,-1.0) to (b4);
    \draw[bend left=20,postaction={decorate}] (10,-.3) to (b8);

\foreach \from/\to in {a1/a2,a2/a3,a3/a4,a5/a6,a6/a7,a7/a8,b8/b7,b7/b6,b6/b5,b4/b3,b3/b2,b2/b1}
\draw[postaction={decorate}] (\from) to (\to);
\end{scope}
\end{tikzpicture}
\caption{An 18-bar 8-framework.} \label{fig:evenframework}
\end{figure}

The structure we need when $\ell$ is even is shown in figure \ref{fig:evenframework}.
We have vertices $a_0\LL a_k$ ($k=18$ in the figure) and $b_1\LL b_k$, but now there is an extra vertex $b_0$. There are paths $P_i$ between $a_i, b_i$ of
length $\ell/2-1$ for $1\le i\le m$, and length $\ell/2-2$ for $m+1\le i\le k$. ($\ell=8$ and $m=8$ in the figure.) There are upper and 
lower tents as before, but now 
all the tents have apex on the left. There must be an upper tent with apex $a_0$, and one with apex $a_m$, although $m=0$ is permitted.
The upper tents are paired with the lower tents; for each upper tent with base $\{a_i\LL a_j\}$
there is also a lower tent with base $\{b_i\LL b_j\}$, and vice versa. But the apexes shift by one;
if an upper tent has apex $a_i$, the paired lower tent has apex $b_{i+1}$ (or $b_0$ when $i=m$). An important condition, not shown 
in the figure, is: 
\begin{itemize}
	\item for each upper tent-\arb{} $T_i$ say, with apex $a_i$, the paired lower tent-\arb{} $S_{i+1}$ with apex $b_{i+1}$ (or $b_0$, if $i=m$) must be
coarboreal with $T_i$ under the bijection that maps $a_j$ to $b_j$ for each leaf $a_j$ of $T_i$.
\end{itemize}

As before, for each $i\in \{1\LL m\}$, if $a_{i-1}$ is the apex of an upper tent-\arb{} $T_{i-1}$ say, there is a directed edge
from some nonleaf vertex of $T_{i-1}$ (possibly from $a_{i-1}$) to
$a_{i}$; and if $a_{i-1}$ is not the apex of a tent, there is a directed edge from $a_{i-1}$ to $a_i$. So the upper tent-\arbs{} are
connected up to form an \arb{} $T$ with apex $a_0$, and with set of leaves $\{a_{m+1}\LL a_k\}$. Also,
for each $i\in \{1\LL m-1\}$, if $b_{i+1}$ is the apex of a lower tent-\arb{} $S_{i+1}$ say, there is a directed edge
from some nonleaf vertex of $S_{i+1}$ (possibly from $b_{i+1}$) to
$b_{i}$; and if $b_{i+1}$ is not the apex of a tent, there is a directed edge from $b_{i+1}$ to $b_i$. Finally,
there is a directed edge
from some nonleaf vertex of the tent-\arb{} $S_{0}$ with apex $b_0$ (possibly from $b_{0}$ itself) to
$b_{m}$. So the lower tent-\arbs{} are
connected up to form an \arb{} $S$ with apex $b_0$, and with set of leaves $\{b_{m+1}\LL b_k\}$.
We call this graph an {\em $\ell$-framework}. We observe:
\begin{thm}\label{evenadj}
	Let $\ell\ge 5$ be even, and let $F$ be an $\ell$-framework, with notation as above. Let $1\le i<j\le k$. There are three cases:
	\begin{itemize}
		\item If $i,j\le m$, there is a directed path of $T$ between $a_i, a_j$, and there is a directed path of 
			$S$ between $b_i, b_j$.
		\item If $i\le m<j$, either there is a directed path of $T$ between $a_i, a_j$, or there is a directed path of 
			$S$ between $b_i, b_j$, and not both.
		\item If $m< i<j$, there is no directed path of $T$ between $a_i, a_j$, and there is no directed path of 
			$S$ between $b_i, b_j$.
	\end{itemize}
\end{thm}

The {\em transitive closure} $\up{T}$ of an arborescence $T$
is the undirected graph with vertex set $V(T)$ in which vertices $u,v$ are adjacent if and only if some directed path of $T$
contains both of $u,v$.
Let $F$ be an $\ell$-framework (here, $\ell$ may be odd or even). Let $P_1\LL P_k, T,S$ and so on be as in the definition of 
an $\ell$-framework.
Let
$D=\up{T}\cup \up{S}\cup P_1\cupcup P_k$. Thus
$V(D)=V(F)$, and distinct $u,v\in V(D)$ are $D$-adjacent if either they are adjacent in some $P_i$, or there is a directed path of one of $S,T$ between $u,v$.
We say a graph $G$
is a {\em blow-up} of $F$ if
\begin{itemize}
\item $D$ is an induced subgraph of $G$, and for each $t\in V(D)$ there is a clique $W_t$ of $G$, all pairwise 
disjoint and with union $V(G)$; $W_t\cap V(D)=\{t\}$ for each $t\in V(D)$, and $W_t=\{t\}$ for each $t\in V(D)\setminus V(P_1\cupcup P_k)$.
\item For each $t\in V(D)$, there is a linear ordering of $W_t$ with first term $t$, say $(x_1\LL x_n)$ where $x_1=t$. 
It has the property that for all distinct $t,t'\in V(D)$, if $t,t'$ are not $D$-adjacent then $W_t$, $W_{t'}$ are anticomplete, and if $t,t'$ are $D$-adjacent
then $G[W_t, W_{t'}]$ obeys the orderings of $W_t,W_{t'}$, and every vertex of $G[W_t, W_{t'}]$ has positive degree.
(Consequently, if $t,t'$ are $D$-adjacent then $t$ is complete to $W_{t'}$ and vice versa.)
\item If $t,t'\in \{a_1\LL a_k\}$ or $t,t'\in \{b_1\LL b_k\}$, and $t,t'$ are $D$-adjacent, then $W_t$ is complete to $W_{t'}$.
\item For each $t\in V(T)$, if $0\le i\le m$ and $a_i,t$ are $D$-adjacent, then $W_t$ is complete to $W_{a_i}$. For each 
$t\in V(S)$, 
if either
$\ell$ is odd and $i\in \{m+1\LL k\}$, or $\ell$ is even and $i\in \{0\LL m\}$,  and $b_i,t$ are $D$-adjacent, then $W_t$ is complete to 
$W_{b_i}$.
\item For each upper tent-\arb{} $T_j$ with apex $a_j$ say, let $t\in L(T_j)$ and let the path $Q$ of $T$ from $a_0$ to $t$ have
vertices 
$$a_0=y_1\CC y_p\DD a_j\DD z_1\CC z_q=t$$
in order. Then $W_t$ is complete to $\{y_1\LL y_p, a_j\}$; $W_t$ is anticomplete to $\bigcup_{t\in T\setminus V(Q)}W_t$; and
$G[W_t, \{z_1\LL z_{q-1}\}]$ is a 
half-graph that obeys the ordering of $W_t$ and the ordering $z_1\LL z_{q-1}$ of  $\{z_1\LL z_{q-1}\}$. The same holds for lower tent-\arbs{} with $T,a_0$ replaced by $S, b_0$.
\end{itemize}
Our main theorem states:
\begin{thm}\label{mainthm}
Let $G$ be a graph with no clique cutset and no universal vertex, and let $\ell\ge 7$. Then $G$ is $\ell$-holed if and only if either
$G$ is a blow-up of a cycle of length $\ell$, or $G$ is a blow-up of an $\ell$-framework.
\end{thm}
We will often have two graphs $F,G$, and a clique $W_t$ of $G$ for each $t\in V(F)$, pairwise vertex-disjoint. For $C\subseteq F$,
we denote $\bigcup_{t\in V(C)}W_t$ by $W(C)$.

The paper is organized as follows. In the next section we show the easier ``if'' half of \ref{mainthm}, that graphs with the structure specified in the
theorem are $\ell$-holed. The remainder of the paper concerns the ``only if'' half. 
In section \ref{sec:cycleblowup} we show that every $\ell$-holed graph that has no clique cutset or universal vertex, and that contains no 
theta, pyramid or prism (defined later), 
is a blow-up of a cycle. Then we turn to $\ell$-holed graphs with no clique cutset or universal vertex that do contain a theta, pyramid or prism; and they will turn out to be blow-ups of $\ell$-frameworks.

An $\ell$-framework consists of the paths $P_1\LL P_k$, and the transitive closure of two \arbs{},
related somehow according to whether $\ell$ is odd or even. If $G$ is a blow-up of an $\ell$-framework $F$, 
then with notation as before, $G$
can be thought of as consisting of three parts, induced respectively on $W(T)$, $W(P_1\cupcup P_k)$, and $W(S)$. The first and third of these
are disjoint and anticomplete, but the first and second overlap, as do the second and third.
In sections \ref{sec:frames} to \ref{sec:parts}, we show that if a graph $G$ is $\ell$-holed, and has
no clique cutset or universal vertex, and contains a theta, pyramid
or prism, then $G$ decomposes into three parts in this way, where the second is as it should be, and
the intersection of first and second, and between second and third, is as it should be, but
we know little about the interiors of the first and third. Then in section \ref{sec:bordercon} we show that the first and third parts are both blow-ups of transitive closures of \arbs{}.
Finally in sections \ref{sec:oddcase}--\ref{sec:evencase} we ask how these two \arbs{} must be related, and this depends on whether $\ell$ is odd or even.

The work reported in this paper is the product of two groups of researchers, working independently.
Much of it forms part of the 
PhD thesis~\cite{linda} of Linda Cook; and also much of it appears in the PhD theses of Cl\'eoph\'ee Robin~\cite{cleophee} 
and of Jake Horsfield~\cite{jake}.
Several authors of this paper have made their own (equivalent) version of the results available as a manuscript on arXiv~\cite{europe-version}.
Their approach is not the same and may be of interest to the readers of this paper.

\section{The ``if'' half}
In this section we prove the ``if'' half of \ref{mainthm}.
We will need the following lemma:
\begin{thm}\label{treepair}
Let $R$ be a tree, and let $X\subseteq V(R)$, with $|X|=n$ where $n$ is even.  Let us say a partition $\{X_1\LL X_{n/2}\}$ of 
$X$ into sets of size two is {\em feasible} (in $R$) if there are $n/2$ vertex-disjoint paths $P_1\LL P_{n/2}$ of $R$,
such that, for $1\le i\le n/2$, the ends of $P_i$ are the two members of $X_i$. 
There is at most one feasible partition of $X$.
\end{thm}
\Proof
We proceed by induction on $|V(R)|$. If $|V(R)|\le 2$ the result is clear, so we assume that $|V(R)|\ge 3$.
If some leaf $t$ of $R$ does not belong to $X$,
we may delete $t$
without affecting which partitions are feasible, and the result follows from the inductive hypothesis; so we may assume that
every leaf belongs to $X$.
Let $t$ be an end of a longest path of $R$, and let $s$ be its neighbour in $R$.
Thus every neighbour of $s$ is a leaf except possible one.
If $s\notin X$, let $X'=(X\setminus \{t\})\cup \{s\}$; then for every feasible partition $\{X_1\LL X_{n/2}\}$ of $X$, with $t\in X_1$
say, the partition $\{(X_1\setminus \{t\})\cup \{s\}, X_2\LL X_{n/2}\}$ is a feasible partition of $X'$, and the claim follows 
from the inductive hypothesis applied to $R\setminus t$ and $X'$.
Thus we may assume that $s\in X$. If some neighbour of $s$ different from $t$ is a leaf, and hence belongs to $X$, there are no 
feasible
partitions; so we may assume that $s$ has degree exactly two (since $|V(R)|\ge 3$). Consequently $R'=R\setminus \{s,t\}$ is a tree.
If $\{X_1\LL X_{n/2}\}$ is a feasible partition of $X$ in $R$, then $\{s,t\}$ is one of its sets, say $X_1=\{s,t\}$, and
$\{X_2\LL X_{n/2}\}$ is a feasible partition of $X'=X\setminus \{s,t\}$ in 
$R'$; and hence it is unique from the inductive hypothesis. This proves \ref{treepair}.~\bbox

Let $G$ be a graph and let $u,v\in V(G)$ be distinct. We say that $u$ {\em $G$-dominates} $v$ if $u,v$ are $G$-adjacent
and
every $G$-neighbour of $v$ is equal or $G$-adjacent to $u$. We say a vertex $v$ is {\em mixed} on a set $C$ if $v\notin C$ and 
$v$ has a neighbour and a non-neighbour in $C$.
In this section we prove:
\begin{thm}\label{easymainthm}
For $\ell\ge 5$, if $G$ is a blow-up of a cycle of length $\ell$, or $G$ is a blow-up of an $\ell$-framework, then $G$ is $\ell$-holed.
\end{thm}
\Proof The first statement is proved in \cite{kristina}, but since the proof is very short we give a proof anyway.
We assume first that $G$ is a blow-up of an $\ell$-cycle; let $W_1\LL W_{\ell}$ be as in the definition of a 
blow-up of an $\ell$-cycle, and let $C$ be a hole of $G$. 
Suppose that $|V(C)\cap W_i|\ge 2$ for some $i$, say $i = 2$. Since $W_2$
is a clique, there is a four-vertex path $c_1\DD c_2\DD c_3\DD c_4$ of $C$ such that $c_2,c_3\in W_2$ and $c_1,c_4\notin W_2$.
Consequently $c_1,c_4\in W_1\cup W_3$. Since $W_1$ is a clique, not both $c_1,c_4\in W_1$, and similarly they are not both in 
$W_3$; so we may assume that $c_1\in W_1$ and $c_4\in W_3$. But this contradicts that $G[W_2, W_{3}], G[W_2, W_{1}]$ are compatible.
Thus $|V(C)\cap W_i|\le 1$ for $1\le i\le \ell$. Since $C$ is a hole, it follows that $|V(C)\cap W_i|= 1$ for $1\le i\le \ell$,
and hence $C$ has length $\ell$, as required.

Now we assume that $G$ is a blow-up of an $\ell$-framework. Let $F$ be an $\ell$-framework, and let $P_1\LL P_k, T,S$ and so on be as in the definition of
an $\ell$-framework. Thus $F=T\cup S\cup P_1\cupcup P_k$.
Let $D=\up{T}\cup \up{S}\cup P_1\cupcup P_k$, and 
let $A=G[\{a_1\LL a_k\}]$ and $B=G[\{b_1\LL b_k\}]$.
From the definition, we have:
\begin{enumerate}[label={\bf(\roman*)}]
\item $D$ is an induced subgraph of $G$, and for each $t\in V(D)$ there is a clique $W_t$ of $G$, all pairwise
disjoint and with union $V(G)$; $W_t\cap V(D)=\{t\}$ for each $t\in V(D)$, and $W_t=\{t\}$ for each $t\in V(D)\setminus V(P_1\cupcup P_k)$.
\item For each $t\in V(D)$, there is a linear ordering of $W_t$ with first term $t$, say $(x_1\LL x_n)$ where $x_1=t$.
It has the property that for all distinct $t,t'\in V(D)$, if $t,t'$ are not $D$-adjacent then $W_t$, $W_{t'}$ are anticomplete, and if $t,t'$ are $D$-adjacent
then $G[W_t, W_{t'}]$ obeys the orderings of $W_t,W_{t'}$, and every vertex of $G[W_t, W_{t'}]$ has positive degree.
\item If $t,t'\in \{a_1\LL a_k\}$ or $t,t'\in \{b_1\LL b_k\}$, and $t,t'$ are $D$-adjacent, then $W_t$ is complete to $W_{t'}$.
\item For each $t\in V(T)$, if $0\le i\le m$ and $a_i,t$ are $D$-adjacent, then $W_t$ is complete to $W_{a_i}$. For each            
$t\in V(S)$,
if either
$\ell$ is odd and $i\in \{m+1\LL k\}$, or $\ell$ is even and $i\in \{0\LL m\}$,  and $b_i,t$ are $D$-adjacent, then $W_t$ is complete to 
$W_{b_i}$.
\item For each upper tent-\arb{} $T_j$ with apex $a_j$ say, let $t\in L(T_j)$ and let the path $Q$ of $T$ from $a_0$ to $t$ have
vertices 
$$a_0=y_1\CC y_p\DD a_j\DD z_1\CC z_q=t$$
in order. Then $W_t$ is complete to $\{y_1\LL y_p, a_j\}$; $W_t$ is anticomplete to $W(T\setminus V(Q))$; and
$G[W_t, \{z_1\LL z_{q-1}\}]$ is a 
half-graph that obeys the given order of $W_t$ and the order $z_1\LL z_{q-1}$. The same holds for lower tent-\arbs{} with $T,a_0$ replaced by $S, b_0$.
\end{enumerate}

Now suppose that $C$ is a hole of $G$ of length different from $\ell$, and choose $C$ with $V(C)\setminus V(D)$ minimal.
\\
\\
(1) {\em $|V(C)\cap W_t|\le 1$ for each $t\in V(D)$.}
\\
\\
Suppose that $|V(C)\cap W_t|\ge 2$ for some $t\in V(D)$. As before, there is a 
four-vertex path $c_1\DD c_2\DD c_3\DD c_4$ of $C$ such that $c_2,c_3\in W_t$ and $c_1,c_4\notin W_t$.
We may assume that $c_2$ is earlier than $c_3$ in the ordering of $W_t$; but $c_4$ is $G$-adjacent to $c_3$ and not to $c_2$,
contrary to {\bf (ii)} above.
This proves (1).
\\
\\
(2) {\em $C$ is a hole of $D$.}
\\
\\
Suppose not; then there exists $t\in V(D)$ such that $V(C)\cap W_t\not\subseteq \{t\}$; and hence, by (1), $t\notin V(C)$ and 
$V(C)\cap W_t\ne \emptyset$. Choose such a vertex $t$ with $t\notin V(A\cup B)$
if possible. Let $t'\in V(C)\cap W_t$,
 and let $r',s'$ be the neighbours of $t'$ in $C$. Since $|W_t|\ge 2$ it follows that $t\in V(P_i)$ for some $i\in \{1\LL k\}$.
From {\bf (ii)} above, $r',s'$ are $G$-adjacent to $t$. From
the minimality of $V(C)\setminus V(D)$, we cannot replace $t'$ by $t$ to obtain a hole, and so $t$ has a neighbour $q'\in V(C)$,
where $q'\ne r',t',s'$.  
Let $r'\in W_r$, and $s'\in W_s$, and $q'\in W_q$. By (1), $q,r,s,t$ are all different.
Since $r',s',q'$ all have neighbours in $W_t$, 
it follows that $t$ is $D$-adjacent to each of $q,r,s$, and so has degree at least three in $D$; and hence 
$t\in \{a_i,b_i\}\subseteq V(A\cup B)$. From the choice of $t$ it follows that every vertex of $C$ not in $V(D)$ belongs to 
$W(A\cup B)$. 

We assume that $t=a_i$ (the argument
when $t=b_i$ is similar and we omit it).
From the symmetry between $r,s$, we may assume that $s\notin V(P_i)$. 
Thus $s\in V(T)$. Suppose first that $q\notin V(T)$. Thus $q$ is the neighbour of $t=a_i$ in $P_i$. 
Since $q'\in V(C)\setminus W(A\cup B)$, it follows that $q'\in V(D)$ and so $q'=q$; and so $q'$ is complete to $W_t$, a 
contradiction since $q',t'$ are not $G$-adjacent. 
Thus $q\in V(T)$. Since $q'$ is mixed on $W_t$, it follows from {\bf (iv)} above that $i>m$, and so $t=a_i$ is a leaf of $T$.

Let $Q$ be the path of $T$ from $a_0$ to $t$. Since $q,s\in V(T)$, and each is $\up{T}$-adjacent to $t$, 
it follows that $q,s\in V(Q)$, and in particular, $q,s\notin L(T)$. Since $q'\in W_q$ is $G$-adjacent to $t$ and not to $t'$, 
and since $t = a_i$,
{\bf (iii)} above implies that 
$q\notin \{a_0\LL a_m\}$ and so $q\notin A$, and therefore $q'=q$ by {\bf (1)} above. There is a vertex $p'$ of $C$ $G$-adjacent to $q'$ and not to $s'$; let 
$p'\in W_p$
where $p\in V(D)$. (Thus $p\ne q,s,t$, but possibly $p=r$.) Since $q\notin A$, and $p,q$ are $D$-adjacent, it follows that $p\in V(T)$. 

By {\bf(v)} above, since $s'$ (and hence $s$) is $G$-adjacent to $t'$ and $q'=q$ is not $G$-adjacent
to $t'$, it follows that $q$ lies in the interior of the subpath of $Q$ between $s,t$, and so there is a directed subpath of $Q$ 
from $s$ to $q$. There is also a directed path $R$ of $T$ between $p,q$, since $p,q$ are $D$-adjacent. If $R$ is from $p$ to $q$, then
$p\in V(Q)$; and if $R$ is from $q$ to $p$ then there is a directed path of $T$ from $s$ to $p$. In either case $p,s$ are $D$-adjacent.
Since $W_p$ is not complete to $W_s$, it follows from {\bf(iv)} above that $p,s\notin \{a_0\LL a_m\}$. 
Since $s \not \in L(T)$ (as we saw earlier),
and $A \subseteq \{a_0 \LL a_m\} \cup L(T)$, it follows that $s\notin A$; so $|W_s| =1$ by {\bf(i)} and thus $s'=s$. 
Since $p,s$ are $D$-adjacent and $W_p$ is not complete to $W_s$, it follows that $|W_p| > 1$. Since $p \not \in \{a_0 \LL a_m\}$, 
it follows from {\bf(i)} above that $p\in \{a_{m+1}\LL a_k\}$, and so 
$p\in L(T)$. But $p'\in W_p$ is $G$-adjacent to $q'=q$ and not to $s'=s$, contradicting {\bf(v)} above.
This proves (2).

\bigskip

Since $\up{T},\up{S}$ do not contain a 4-vertex induced path or a 4-hole, it follows that $C$ is not a subgraph
of either of them. For $1\le i\le k$, each internal vertex of $P_i$ has degree two in $D$, and so if some such vertex belongs to $C$
then so does the whole of $P_i$. Hence $C$ is the concatenation of paths of $G$, in alternation equal to one of $P_1\LL P_k$ or to a subpath
of one of  $\up{T},\up{S}$ with length one or two (because $\up{T},\up{S}$ contain no induced path of length three). Let us write
$$C=P_{d_1}\cup Q_{1}\cup P_{d_2}\cup Q_2\cup P_{d_3}\cup Q_3\cupcup P_{d_n}\cup Q_n$$
where $d_1,d_2\LL d_n\in \{1\LL k\}$ are distinct, and for $i$ odd $Q_i$ is a path of $\up{T}$ with ends $a_{d_i}, a_{d_{i+1}}$,
and for $i$ even $Q_i$ is a path of $\up{S}$ with ends $b_{d_i}, b_{d_{i+1}}$, where $d_{n+1}=d_1$. Thus $n\ge 2$ and $n$ is even.

Suppose first that $\ell$ is odd. If $n\ge 4$, then there are no edges of $G$ between $P_{d_1}, P_{d_3}$ by \ref{oddadj}, contradicting that
$F$ is an $\ell$-framework. So $n=2$; and exactly one of $Q_1,Q_2$ has length one, from the definition of an $\ell$-framework,
and the other has length two; and since $P_{d_1},P_{d_2}$ both have length $(\ell-3)/2$ from the definition of an $\ell$-framework,
it follows that $C$ has length $\ell$, a contradiction.

Thus $\ell$ is even. Note that in this case there is a symmetry in the definition of an $\ell$-framework that exchanges
$S$ and $T$ (although the bars have to be renumbered). This will be helpful to reduce the number of cases we need to examine.

Suppose that $n=2$. If $P_{d_1}, P_{d_2}$ both have length $\ell/2-1$, then $Q_1,Q_2$ both have length one
(from \ref{evenadj} and the definition of an $\ell$-framework) and so $C$ has length $\ell$; if exactly one of $P_{d_1}, P_{d_2}$ has length 
$\ell/2-1$, then exactly one of $Q_1,Q_2$ has length one (again from \ref{evenadj} and the definition) and so $C$ has length $\ell$; and if
neither of $P_{d_1}, P_{d_2}$ has length $\ell/2-1$, then neither of $Q_1,Q_2$ has length one (again from \ref{evenadj} and the definition) 
and so $C$ has length $\ell$, in all cases a contradiction. So $n\ge 4$. 

Suppose that $P_{d_1}$ has length $\ell/2-1$. There are no edges between $V(P_{d_1})$ and $V(P_{d_3})$, and in particular
$a_{d_1}, a_{d_3}$ are $G$-nonadjacent and $b_{d_1}, b_{d_3}$ are $G$-nonadjacent, contrary to \ref{evenadj}. Thus
all of $P_{d_1}\LL P_{d_n}$ have length $\ell/2-2$, and therefore all of $Q_1\LL Q_n$ have length two. Let $q_i$
be the middle vertex of $Q_i$ for $1\le i\le n$.

Each of $a_{d_1}\LL a_{d_n}$ belongs to the base of one of the tent-\arbs{} $T_0\LL T_m$, and each of $b_{d_1}\LL b_{d_m}$
belongs to the base of one of $S_0\LL S_m$.
\\
\\
(3) {\em There exists $i\in \{0\LL m\}$ such that $a_{d_1}\LL a_{d_n}$ all belong to the base of $T_i$ and $b_{d_1}\LL b_{d_n}$
all belong to the base of $S_{i+1}$.}
\\
\\
Suppose not; then from the symmetry between $S,T$, we may assume that $a_{d_1}\in V(T_g)$ and $a_{d_2}\in V(T_h)$ where $g<h$. 
Hence $b_{d_1}\in S_{g+1}$ and $b_{d_2}\in S_{h+1}$.
Since $q_1$ is $G$-adjacent (and hence $\up{T}$-adjacent) to both $a_{d_1},a_{d_2}$, it follows that $q_1$ belongs to the path 
of $T$
between $a_0, a_{g+1}$. Consequently $q_1$ is $G$-adjacent to every vertex in $V(T_{g+1})\cupcup V(T_m)$; and so none of
$a_{d_3}\LL a_{d_n}, q_3,q_5\LL q_{n-1}$ belong to $V(T_{g+1})\cupcup V(T_m)$. Let $a_{d_3}\in T_{f}$; it follows that
$f\le g$. Moreover, $b_{d_3}\in S_{f+1}$, and since $q_2$ is $G$-adjacent to both $b_{d_2}, b_{d_3}$, it follows that
$q_2$ belongs to the path of $S$ between $b_0$ and $b_h$. But then $q_2$ is $G$-adjacent to $b_{d_1}$, a contradiction. This proves (3).

\bigskip

Choose $i$ as in (3). Since $q_1$ has a neighbour and a non-neighbour in the base of $T_i$, it follows that $q_1\in V(T_i)$,
and similarly $q_3,q_5\LL q_{n-1}\in V(T_i)$ and $q_2,q_4,q_6\LL q_n \in V(S_{i+1})$. 
From the definition of an $\ell$-framework, the arborescences $T_i, S_{i+1}$ are coarboreal under the bijection
that maps $a_j$ to $b_j$ for each leaf $a_j$ of $T_i$. Let $S'$ be the arborescence obtained from $S_{i+1}$ by replacing $b_j$ by $a_j$
for each leaf $a_j$ of $T_i$. Thus $T_i, S'$ are coarboreal; let $R$ be a tree with vertex set $L(T_i)$ in which they both live.
For $j\in \{1\LL n\}$ with $j$ odd, since $q_j$ is $G$-adjacent to $a_{d_j}, a_{d_{j+1}}$ and to no other vertices in $\{a_{d_1}\LL a_{d_n}\}$, and the set of 
$G$-neighbours of $q_j$
in $V(R)$ is the vertex set of a subtree of $R$, there is a path $R_j$ of $R$ between $a_{d_j}, a_{d_{j+1}}$, such that all its
vertices are $\up{T_i}$-adjacent to $q_j$, and 
$$V(R_j)\cap \{a_{d_1}\LL a_{d_n}\}= \{a_{d_j}, a_{d_{j+1}}\}.$$
Similarly, for $j\in \{1\LL n\}$ with $j$ even, there is a path $R_j$ of $R$ between $a_{d_j}, a_{d_{j+1}}$, such that all its
vertices are $\up{S'}$-adjacent to $q_j$, and 
$$V(R_j)\cap \{a_{d_1}\LL a_{d_n}\}= \{a_{d_j}, a_{d_{j+1}}\},$$ 
where $d_{n+1}=d_1$.
The paths $R_1,R_3\LL R_{n-1}$ are pairwise vertex-disjoint, since $q_1,q_3,q_5\LL q_{n-1}$ are pairwise $G$-nonadjacent and therefore
no two of them have a common neighbour in $L(T)$; and similarly $R_2,R_4\LL R_n$ are pairwise vertex-disjoint. 
It follows that 
$$\{\{a_{d_1},a_{d_2}\},\{a_{d_3},a_{d_4}\}\LL \{a_{d_{n-1}},a_{d_n}\}\}$$
is a feasible partition of $\{a_{d_1}\LL a_{d_n}\}$ in $R$, and so is 
$$\{\{a_{d_2},a_{d_3}\},\{a_{d_4},a_{d_5}\}\LL \{a_{d_{n}},a_{d_1}\}\};$$
and they are different since $n\ge 4$, contrary to \ref{treepair}. This proves \ref{easymainthm}.~\bbox
\section{Blow-ups of cycles}\label{sec:cycleblowup}

Our strategy to prove \ref{mainthm} is to choose a maximal ``apexed $\ell$-frame'' (defined in section \ref{sec:gates}), and a maximal 
blow-up of it contained in our graph, and to analyze how the remainder of the
graph can attach to it. But this only works in graphs that contain apexed $\ell$-frames, and here we handle the graphs that do not. 

Here are some types of graph that we will need:
\begin{itemize}
\item A {\em theta} is a graph that is the union of three paths $R_1,R_2,R_3$, each with the same pair of ends,
each of length more than one, and pairwise vertex-disjoint except for their ends.
\item A {\em pyramid} is a graph that is the union of three paths $R_1,R_2,R_3$ and three additional edges $b_1b_2,b_2b_3,b_3b_1$,
where $R_i$ has ends $a\ne b_i$ for $i=1,2,3$, and $R_1,R_2,R_3$ are pairwise vertex-disjoint
except for their common end $a$, and at least two of $R_1,R_2,R_3$
have length at least two.
\item A {\em prism} is a graph that is the union of three paths $R_1,R_2,R_3$ and six additional edges
$a_1a_2,a_2a_3,a_3a_1,b_1b_2,b_2b_3,b_3b_1$, where $R_i$ has ends $a_i\ne b_i$ for $i=1,2,3$, and $R_1,R_2,R_3$ are pairwise
vertex-disjoint.
\item A {\em proper wheel} is a graph consisting of a cycle $C$ of length at least four, and one additional vertex $v$, where $v$ has at least three neighbours and at least one non-neighbour in $V(C)$, and if $v$ has exactly three neighbours in $V(C)$ then they do not induce a path.
\end{itemize}

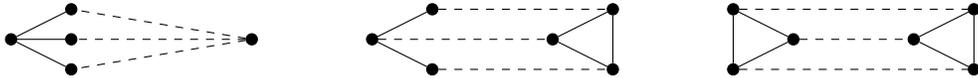
\begin{figure}[H]
\centering

\begin{tikzpicture}[scale=.8,auto=left]
\tikzstyle{every node}=[inner sep=1.5pt, fill=black,circle,draw]
\node (z) at (0,0) {};
\node (a) at (1,.5) {};
\node (b) at (1,0) {};
\node (c) at (1,-.5) {};
\node (d) at (4,0) {};

\foreach \from/\to in {z/a,z/b, z/c}
\draw [-] (\from) -- (\to);
\foreach \from/\to in {d/a,d/b, d/c}
\draw [dashed] (\from) -- (\to);

\node (a) at (6,0) {};
\node (c1) at (7,.5) {};
\node (c3) at (7,-.5) {};
\node (b1) at (10,.5) {};
\node (b2) at (9,0) {};
\node (b3) at (10,-.5) {};

\foreach \from/\to in {a/c1, a/c3,b1/b2,b1/b3,b2/b3}
\draw [-] (\from) -- (\to);
\foreach \from/\to in {c1/b1, a/b2,c3/b3}
\draw [dashed] (\from) -- (\to);

\node (s1) at (12,.5) {};
\node (s2) at (13,0) {};
\node (s3) at (12,-.5) {};
\node (t1) at (16,.5) {};
\node (t2) at (15,0) {};
\node (t3) at (16,-.5) {};
\foreach \from/\to in {s1/s2,s1/s3,s2/s3,t1/t2,t1/t3,t2/t3}
\draw [-] (\from) -- (\to);
\foreach \from/\to in {s1/t1,s2/t2,s3/t3}
\draw [dashed] (\from) -- (\to);

\end{tikzpicture}

\caption{A theta, a pyramid and a prism (dashed lines mean paths of arbitrary positive length)} \label{theta}
\end{figure}

In this paper, we say {\em $G$ contains $H$} to mean $G$ has an induced subgraph isomorphic to $H$.
The following is a corollary of theorem 1.6 of a paper~\cite{kristina} by  V. Boncompagni, I. Penev, and K. Vu\v{s}kovi\'{c}:
\begin{thm}\label{kristina}
If a graph $G$ has no clique cutset or universal vertex, and $G$ contains no theta, pyramid, prism, proper wheel, or 4-hole, then $G$ is a blow-up of a cycle.
\end{thm}

Since every proper wheel has two holes of different lengths, \ref{kristina} has the following immediate consequence:
\begin{thm}\label{cycleblowup}
Let $\ell\ge 5$, and let $G$ be a non-null $\ell$-holed graph, not containing a theta, pyramid or prism, and with no clique cutset
or universal vertex. Then $G$ is a blow-up of a cycle of length $\ell$.
\end{thm}

\section{Frames}\label{sec:frames}

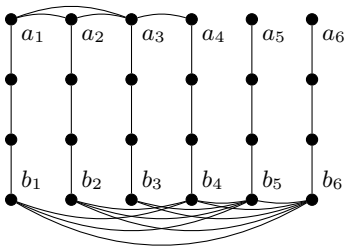
\begin{figure}[H]
\centering
\begin{tikzpicture}[scale=.8,auto=left]
\tikzstyle{every node}=[inner sep=1.5pt, fill=black,circle,draw]

\node (a1) at (0,3) {};
\node (a2) at (1,3) {};
\node (a3) at (2,3) {};
\node (a4) at (3,3) {};
\node (a5) at (4,3) {};
\node (a6) at (5,3) {};

\node (b1) at (0,0) {};
\node (b2) at (1,0) {};
\node (b3) at (2,0) {};
\node (b4) at (3,0) {};
\node (b5) at (4,0) {};
\node (b6) at (5,0) {};

\foreach \from/\to in {a1/a2,a2/a3,a3/a4}
\draw [bend left=15] (\from) to (\to);
\foreach \from/\to in {a1/a3}
\draw [bend left=20] (\from) to (\to);

\foreach \from/\to in {b4/b5, b5/b6}
\draw [bend right=10] (\from) to (\to);
\foreach \from/\to in {b2/b4, b3/b5, b4/b6}
\draw [bend right=15] (\from) to (\to);
\foreach \from/\to in {b1/b4, b2/b5, b3/b6}
\draw [bend right=20] (\from) to (\to);
\foreach \from/\to in {b1/b5, b2/b6}
\draw [bend right=25] (\from) to (\to);
\foreach \from/\to in {b1/b6}
\draw [bend right=30] (\from) to (\to);

\node (m1) at (0,2) {};
\node (n1) at (0,1) {};
\node (m2) at (1,2) {};
\node (n2) at (1,1) {};
\node (m3) at (2,2) {};
\node (n3) at (2,1) {};
\node (m4) at (3,2) {};
\node (n4) at (3,1) {};
\node (m5) at (4,2) {};
\node (n5) at (4,1) {};
\node (m6) at (5,2) {};
\node (n6) at (5,1) {};

\foreach \from/\to in {a1/m1,m1/n1,n1/b1,a2/m2,m2/n2,n2/b2,a3/m3,m3/n3,n3/b3,a4/m4,m4/n4,n4/b4,a5/m5,m5/n5,n5/b5,a6/m6,m6/n6,n6/b6}
\draw (\from) -- (\to);

\tikzstyle{every node}=[]
\draw (a1) node [below right]           {\footnotesize$a_1$};
\draw (a2) node [below right]           {\footnotesize$a_2$};
\draw (a3) node [below right]           {\footnotesize$a_3$};
\draw (a4) node [below right ]           {\footnotesize$a_4$};
\draw (a5) node [below right ]           {\footnotesize$a_5$};
\draw (a6) node [below right ]           {\footnotesize$a_6$};
\draw (b1) node [above right]           {\footnotesize$b_1$};
\draw (b2) node [above right]           {\footnotesize$b_2$};
\draw (b3) node [above right]           {\footnotesize$b_3$};
\draw (b4) node [above right ]           {\footnotesize$b_4$};
\draw (b5) node [above right ]           {\footnotesize$b_5$};
\draw (b6) node [above right ]           {\footnotesize$b_6$};
\end{tikzpicture}
\caption{A $9$-frame.} \label{fig:oddframe}
\end{figure}

The aim of this section is to define an ``$\ell$-frame'', which is part of an $\ell$-framework, but we need several preliminary definitions.
We denote the four-vertex path by $\mathcal{P}_4$ (it is usually denoted by $P_4$, but
we have sets of paths in this paper, and would like to number their members $P_1,P_2\LL$ and so on.). Similarly, $\mathcal{C}_k$
denotes the cycle of length $k$.
A {\em threshold graph} is a graph that does not contain $\mathcal{P}_4$, $\mathcal{C}_4$ or $\overline{\mathcal{C}_4}$. 
(The notation $\overline{G}$ means the complement graph of $G$.) 
Let $\ell\ge 5$ be odd. Let $k\ge 3$ be an integer, and take distinct vertices 
$a_1\LL a_k, b_1\LL b_k.$ 
For $1\le i\le k$ let $P_i$ be a path of length $(\ell-3)/2$ with ends $a_i, b_i$, pairwise vertex-disjoint. 
Let the subgraphs $A$ and $B$
induced on $\{a_1\LL a_k\}$ and on $\{b_1\LL b_k\}$ respectively be threshold graphs, and for $1\le i<j\le k$, let $b_i,b_j$ be adjacent if and only if $a_i, a_j$
are nonadjacent. Moreover, let $A$ either be disconnected or two-connected, and the same for $B$.
(See figure \ref{fig:oddframe}.) Let $F$ be the union of $A,B$ and the paths $P_1\LL P_k$.
For $\ell$ odd, a graph $F$ constructible in this way is called an {\em $\ell$-frame}; all its 
holes have length $\ell$. We call $P_1\LL P_k$ the {\em bars} of the frame, and $A,B$ are its {\em sides}.

Every non-null threshold graph has either a vertex of degree zero, or a vertex adjacent to all other vertices (see~\cite{chvatal}),
and so, since the subgraphs induced on $\{a_1\LL a_k\}$ and $\{b_1\LL b_k\}$ are complementary threshold graphs, one of them
has a vertex of degree 0 (such as $a_6$ in the figure). Thus all $\ell$-frames when $\ell$ is odd have one-vertex clique cutsets.

\begin{figure}[H]
\centering
\begin{tikzpicture}[scale=.8,auto=left]
\tikzstyle{every node}=[inner sep=1.5pt, fill=black,circle,draw]

\node (a1) at (0,3) {};
\node (a2) at (1,3) {};
\node (a3) at (2,3) {};
\node (a4) at (3,3) {};

\node (b1) at (0,0) {};
\node (b2) at (1,0) {};
\node (b3) at (2,0) {};
\node (b4) at (3,0) {};

\node (m1) at (0,2) {};
\node (n1) at (0,1) {};
\node (m2) at (1,2) {};
\node (n2) at (1,1) {};
\node (m3) at (2,2) {};
\node (n3) at (2,1) {};
\node (m4) at (3,2) {};
\node (n4) at (3,1) {};
\foreach \from/\to in {a1/m1,m1/n1,n1/b1,a2/m2,m2/n2,n2/b2,a3/m3,m3/n3,n3/b3,a4/m4,m4/n4,n4/b4}
\draw (\from) -- (\to);

\node (c1) at (5,3) {};
\node (c2) at (6,3) {};
\node (c3) at (7,3) {};
\node (c4) at (8,3) {};

\node (d1) at (5,0) {};
\node (d2) at (6,0) {};
\node (d3) at (7,0) {};
\node (d4) at (8,0) {};

\node (p1) at (5,2.25) {};
\node (q1) at (5,1.5) {};
\node (p2) at (6,2.25) {};
\node (q2) at (6,1.5) {};
\node (p3) at (7,2.25) {};
\node (q3) at (7,1.5) {};
\node (p4) at (8,2.25) {};
\node (q4) at (8,1.5) {};
\node (r1) at (5,.75) {};
\node (r2) at (6,.75) {};
\node (r3) at (7,.75) {};
\node (r4) at (8,.75) {};
\foreach \from/\to in {c1/p1,p1/q1,q1/r1, r1/d1,c2/p2,p2/q2,q2/r2,r2/d2,c3/p3,p3/q3,q3/r3, r3/d3,c4/p4,p4/q4,q4/r4,r4/d4}
\draw (\from) -- (\to);

\foreach \from/\to in {c1/c2,c2/c3,c3/c4}
\draw[bend left=10] (\from) to (\to);
\foreach \from/\to in {c1/c3,c2/c4}
\draw[bend left=15] (\from) to (\to);
\foreach \from/\to in {c1/c4}
\draw[bend left=20] (\from) to (\to);

\foreach \from/\to in {d1/d2,d2/d3,d3/d4}
\draw[bend right=10] (\from) to (\to);
\foreach \from/\to in {d1/d3,d2/d4}
\draw[bend right=15] (\from) to (\to);
\foreach \from/\to in {d1/d4}
\draw[bend right=20] (\from) to (\to);

\foreach \from/\to in {c1/a4}
\draw[bend right=15] (\from) to (\to);
\foreach \from/\to in {c1/a3, c2/a4}
\draw[bend right=20] (\from) to (\to);
\foreach \from/\to in {c1/a2, c2/a3,c3/a4}
\draw[bend right=25] (\from) to (\to);
\foreach \from/\to in {c2/a2,c3/a3,c4/a4}
\draw[bend right=30] (\from) to (\to);

\foreach \from/\to in {d1/b1}
\draw[bend left=15] (\from) to (\to);
\foreach \from/\to in {d2/b1, d3/b2, d4/b3}
\draw[bend left=20] (\from) to (\to);
\foreach \from/\to in {d3/b1, d4/b2}
\draw[bend left=25] (\from) to (\to);
\foreach \from/\to in {d4/b1}
\draw[bend left=30] (\from) to (\to);

\tikzstyle{every node}=[]
\draw (a1) node [below right]           {\footnotesize$a_1$};
\draw (a2) node [below right]           {\footnotesize$a_2$};
\draw (a3) node [below right]           {\footnotesize$a_3$};
\draw (a4) node [below right ]           {\footnotesize$a_4$};
\draw (b1) node [above right]           {\footnotesize$b_1$};
\draw (b2) node [above right]           {\footnotesize$b_2$};
\draw (b3) node [above right]           {\footnotesize$b_3$};
\draw (b4) node [above right ]           {\footnotesize$b_4$};

\draw (c1) node [below right]           {\footnotesize$c_1$};
\draw (c2) node [below right]           {\footnotesize$c_2$};
\draw (c3) node [below right]           {\footnotesize$c_3$};
\draw (c4) node [below right ]           {\footnotesize$c_4$};
\draw (d1) node [above right]           {\footnotesize$d_1$};
\draw (d2) node [above right]           {\footnotesize$d_2$};
\draw (d3) node [above right]           {\footnotesize$d_3$};
\draw (d4) node [above right ]           {\footnotesize$d_4$};
\end{tikzpicture}
\caption{A $10$-frame.} \label{fig:evenframe}
\end{figure}
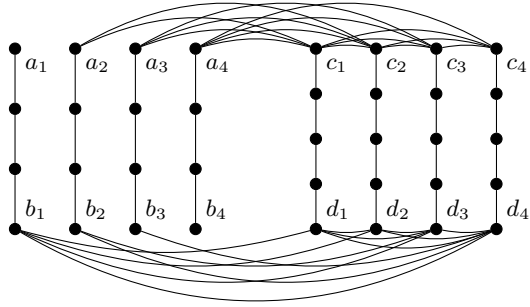

Now the case when $\ell$ is even. Let $m,n\ge 0$ be integers with $n\ge 2$ and $m+n\ge 3$; and let
$$a_1\LL a_n, c_1\LL c_m $$
$$b_1\LL b_n, d_1\LL d_m$$
all be distinct. For $1\le i\le n$ let $P_i$ be a path with ends $a_i, b_i$ of length $\ell/2-2$, and for $1\le i\le m$ let $Q_i$
be a path between $c_i, d_i$ of length $\ell/2-1$, all pairwise vertex-disjoint. 
Let $A,B$ be graphs with vertex sets $\{a_1\LL a_n, c_1\LL c_m\}$ and $\{b_1\LL b_n,d_1\LL d_m\}$ respectively,
with the following properties.
Let $\{c_1\LL c_m\}$ and $\{d_1\LL d_m\}$ be cliques; and let 
the bipartite subgraph 
$A[\{a_1\LL a_n\}, \{c_1\LL c_m\}]$ be a half-graph. For $1\le i\le n$ and $1\le j\le m$, let $b_i, d_j$
be adjacent if and only if $a_i,c_j$ are nonadjacent. Let one of $a_1\LL a_n$ have degree zero in $A$, and let 
one of $b_1\LL b_n$ have degree zero in $B$.
There are no other edges (and thus $\{a_1\LL a_n\}$ and $\{b_1\LL b_n\}$ are stable sets).
It follows that $A,B$ are both disconnected threshold graphs.
Let $F$ be the union of $A,B$ and the paths $P_1\LL P_n$ and $Q_1\LL Q_m$. 
We call such a graph $F$ an {\em $\ell$-frame}. (See figure \ref{fig:evenframe}.)
We call $P_1\LL P_n,Q_1\LL Q_m$ the {\em bars} of the frame, and 
$A,B$ its {\em sides}.
Every hole in an $\ell$-frame has length $\ell$.

\section{Gates}\label{sec:gates}

Let $k\ge 3$ be an integer. A {\em $k$-bar semigate} in a graph $G$ is a $(k+2)$-tuple $(A,B,P_1\LL P_k)$ of induced subgraphs of $G$, with the following properties:
\begin{itemize}
\item $A,B$ are vertex-disjoint connected subgraphs of $G$;
\item for $1\le i\le k$, $P_i$ is a path with ends $a_i, b_i$, and $V(P_i\cap A)=\{a_i\}$ and $V(P_i\cap B)=\{b_i\}$;
\item $P_1\LL P_k$ are all distinct, and $V(P_i\cap P_j)=\{a_i,b_i\}\cap \{a_j,b_j\}$ for $ 1\le i<j\le k$;
\item each vertex of $A$ has at most one  neighbour in $V(B)$, and vice versa;
\item $A\cup B\cup P_1\cupcup P_k$ is an induced subgraph of $G$.
\end{itemize}
We call $A, B$ the {\em sides} of the semigate, and $P_1\LL P_k$ its {\em bars}. (Sadly for the metaphor, we will draw 
semigates with the bars vertical and the sides horizontal, because they fit better on the page that way.) Note that 
a theta makes a 3-bar semigate with $|A|=|B|=1$; and similarly a pyramid makes one with $|A|=1$ and $|B|=3$, and a prism makes one with
$|A|=|B|=3$.

Let $(A,B,P_1\LL P_k)$ and $(A',B', P_1'\LL P_k')$ be $k$-bar semigates in $G$. We say the second is {\em better} than the first
if $A'\subseteq A$ and $B'\subseteq B$ and $P_i\subseteq P_i'$ for $1\le i\le k$, and {\em strictly better} if either $A'\ne A$ 
or $B'\ne B$. (Note that the 
sides are getting smaller and the bars are getting bigger.)
A $k$-bar semigate is a {\em $k$-bar gate} if there is no strictly better $k$-bar semigate. For every $k$-bar semigate, there is a $k$-bar gate 
which is better than it.
Given a gate in this notation, we want to see what we can say about $A$ and the vertices $a_1\LL a_k\in V(A)$, and the same in $B$, when $G$ is $\ell$-holed.

In general, if $u,v$ are vertices of some graph $H$, $d_H(u,v)$ denotes the distance in $H$ between $u,v$.
We begin with the following trivial but very useful observation:
\begin{thm}\label{legfact}
Let  $(A,B,P_1\LL P_k)$ be a $k$-bar gate in $G$. Then for $1\le i\le k$, either $a_i$ has at least two neighbours in $A$, or $a_i=a_j$
for some $j\in \{1\LL k\}\setminus \{i\}$.
\end{thm}
\Proof Suppose that $a_1\ne a_2\LL a_k$ say, and $a_1$ has degree at most one in $A$. Since $A$ is connected and $a_2\in V(A)$, 
it follows that $a_1$ has a neighbour $u\in V(A)$. But then $(A\setminus \{a_1\}, B, P_1',P_2\LL P_k)$ is a strictly better $k$-bar semigate, where $P_1'$
is the path obtained by adding $u$ and the edge $a_1u$ to $P_1$, a contradiction. This proves \ref{legfact}.~\bbox

A {\em cut-vertex} in a graph $A$ is a vertex $v\in V(A)$ such that $A\setminus \{v\}$ is disconnected.
\begin{thm}\label{cutvertex}
Let  $(A,B,P_1\LL P_k)$ be a $k$-bar gate in $G$, and suppose that $A$ has a cut-vertex. There exist distinct
$g,h,i,j\in \{1\LL k\}$ such that 
$$d_A(a_g,a_h)+d_A(a_i,a_j)+2\le  d_A(a_g,a_i)+d_A(a_h,a_j)= d_A(a_g,a_j)+d_A(a_h,a_i).$$
\end{thm}
\Proof
Let $v$ be a cut-vertex of $A$, and let $C,D$ be distinct components of $A\setminus \{v\}$. 
\\
\\
(1) {\em There exist distinct  $i,j\in \{1\LL k\}$ such that $a_i, a_j\in V(C)$ and $d_A(a_i,a_j)<d_A(a_i,v)+d_A(a_j,v)$.}
\\
\\
Let $I$ be the set of $i\in \{1\LL k\}$
with $i\in V(C)$, and for each $i\in I$, let $Q_i$ be an induced path of $A[V(C)\cup \{v\}]$ between $a_i, v$ of length $d_A(a_i,v)$, and let
$X$ be the union of the vertex sets of the paths $Q_i\;(i\in I)$. Since
$$((A\setminus V(C))\cup G[X],B, P_1\LL P_k)$$
is not a strictly better $k$-bar semigate, it follows that $X=V(C)$. In particular, since $V(C)\ne \emptyset$, it follows that $I\ne \emptyset$.
Let $i\in I$ say. If $a_i=a_j$ for some $j\in \{1\LL k\}\setminus \{i\}$, then (1) holds; so we may assume there is no such $j$.
By \ref{legfact}, $a_i$ has degree at least two in $A$. Not 
both these neighbours lie in $V(Q_i)$ since $Q_i$ is induced; let $u\in V(A)\setminus V(Q_i)$ be $A$-adjacent to $a_i$. 
Thus $u\in X$, and so $u\in V(Q_j)$ for some $j\in I\setminus \{i\}$. But then $Q_i\cup Q_j$ is not an induced path, and so 
$d_A(a_i,a_j)<d_A(a_i,v)+d_A(a_j,v)$. This proves (1).

\bigskip

From (1) we may assume that $a_g,a_h\in V(C)$ and $d_A(a_g,a_h)<d_A(a_g,v)+d_A(a_h,v)$.
Similarly we may assume that $a_i,a_j\in V(D)$ and $d_A(a_i,a_j)<d_A(a_i,v)+d_A(a_j,v)$.
Summing, we deduce that
$$d_A(a_g,a_h)+d_A(a_i,a_j)\le d_A(a_g,v)+d_A(a_h,v)+ d_A(a_i,v)+d_A(a_j,v) -2.$$
But $d_A(a_g,a_i)=d_A(a_g,v)+d_A(a_i,v)$ since every path of $A$ between $a_g, a_i$ contains $v$, and there are three similar equations.
This proves \ref{cutvertex}.~\bbox

\begin{thm}\label{distance}
Let $\ell\ge 5$, let $G$ be $\ell$-holed, and let  $(A,B,P_1\LL P_4)$ be a $4$-bar gate in $G$. Define
$d(1,2,3,4) = d_A(a_1,a_2)+d_A(a_3,a_4)$, and define $d(1,3,2,4)$ and $d(1,4,2,3)$
similarly.
Then two of 
$$d(1,2,3,4), d(1,3,2,4), d(1,4,2,3)$$ 
are equal and the third is at most one more.
\end{thm}
\Proof
The graph $A$ is connected. If it has a cut-vertex then the claim follows from \ref{cutvertex}, so we may assume that 
it is two-connected or has at most two vertices. In particular, $V(A)=\{a_1,a_2,a_3,a_4\}$ since any other vertex in $A$ could 
be deleted to obtain a strictly better $4$-bar semigate. 
For each $i\in \{1\LL 4\}$ either $a_i=a_j$ for some $j\in \{1\LL 4\}\setminus \{i\}$, or $a_i$ is adjacent to at least two
other members of $\{a_1\LL a_4\}$, by \ref{legfact}; so one of the following holds:
\begin{itemize}
\item $|V(A)|=1$; then $d(1,2,3,4), d(1,3,2,4), d(1,4,2,3)$ are all equal to zero.
\item $|V(A)|=2$; then we may assume that $a_1=a_2\ne a_3=a_4$, and so $d(1,2,3,4)=0$ and $d(1,3,2,4)= d(1,4,2,3)=2$.
\item $|V(A)|=3$; then we may assume that $a_1=a_2$, and $a_1,a_3,a_4$ are distinct, and all three pairs of these vertices are adjacent
since $A$ is two-connected; but then $d(1,2,3,4)=1$ and $d(1,3,2,4)= d(1,4,2,3)=2$.
\item $|V(A)|=4$; then $A$ has a cycle of length four, since it is two-connected, and since $G$ is $\ell$-holed, 
this cycle is not induced. Consequently $A$ has at least five edges, and we may assume that every two of $a_1,a_2,a_3,a_4$
are adjacent except possibly $a_1,a_4$. Thus $d(1,2,3,4)=d(1,3,2,4)=2$ and $2\le d(1,4,2,3)\le 3$.
\end{itemize}
This proves \ref{distance}.~\bbox

\begin{thm}\label{closeddistance}
Let $\ell\ge 5$, let $G$ be $\ell$-holed, and let  $(A,B,P_1\LL P_4)$ be a $4$-bar gate in $G$. Define
$d(1,2,3,4) = d_A(a_1,a_2)+d_A(a_3,a_4)$, and define $d(1,3,2,4)$ and $d(1,4,2,3)$
similarly.
Then every two of 
$$d(1,2,3,4), d(1,3,2,4), d(1,4,2,3)$$ 
differ by at most one. 
\end{thm}
\Proof Suppose that $d(1,2,3,4)\le d(1,3,2,4)-2$ say. By \ref{distance}, $d(1,3,2,4)=d(1,4,2,3)$.
Let $d'(1,2,3,4)= d_B(b_1,b_2)+d_B(b_3,b_4)$, and define $d'(1,3,2,4)$ and $d'(1,4,2,3)$ similarly.
For $1\le i\le 4$, let $P_i$ have length $\ell_i$. Let $1\le i<j\le 4$. 
Since there is a hole formed by the union of $P_i,P_j$ with a shortest 
path of $A$ between $a_i,a_j$, and a shortest path in $B$ between $b_i,b_j$, it follows that
$$d_A(a_i,a_j)+d_B(b_i,b_j)+\ell_i+\ell_j=\ell.$$
In particular, 
$$d(1,2,3,4)+d'(1,2,3,4)=d_A(a_1,a_2)+d_A(a_3,a_4)+d_B(b_1,b_2)+d_B(b_3,b_4)= 2\ell -\ell_1-\ell_2-\ell_3-\ell_4,$$
and similarly 
$$d(1,3,2,4)+d'(1,3,2,4) =d(1,4,2,3)+d'(1,4,2,3)=2\ell -\ell_1-\ell_2-\ell_3-\ell_4.$$
Since $d(1,3,2,4)=d(1,4,2,3)$, it follows that $d'(1,3,2,4)=d'(1,4,2,3)$, and since
$d(1,2,3,4)\le d(1,3,2,4)-2$ it follows that $d'(1,2,3,4)\ge d'(1,3,2,4)+2$, contrary to \ref{distance} applied in $B$.
This proves that every two of 
$$d(1,2,3,4), d(1,3,2,4), d(1,4,2,3)$$ 
differ by at most one. 
This proves \ref{closeddistance}.~\bbox

\bigskip

Next we extend this to semigates:

\begin{thm}\label{dist}
Let $\ell\ge 5$, let $G$ be $\ell$-holed, and let  $(A,B,P_1\LL P_4)$ be a $4$-bar semigate in $G$. Define
$d(1,2,3,4) = d_A(a_1,a_2)+d_A(a_3,a_4)$, and define $d(1,3,2,4)$ and $d(1,4,2,3)$
similarly.
Then every two of $$d(1,2,3,4), d(1,3,2,4), d(1,4,2,3)$$ 
differ by at most one.
\end{thm}
\Proof
Choose a 4-bar gate $(A',B', P_1', P_2', P_3', P_4')$ that is better than $(A,B,P_1\LL P_4)$, and let $P_i'$ have ends $a_i'\in V(A')$ and 
$b_i'\in V(B')$ for $1\le i\le 4$. 
Define $d'(1,2,3,4) = d_{A'}(a_1',a_2')+d_{A'}(a_3',a_4')$, and define $d'(1,3,2,4)$ and $d'(1,4,2,3)$ similarly. For $1\le i\le 4$,
$P_i$ is a subpath of $P_i'$; let the subpath of $P_i'$ between $a_i, a_i'$ have length $\ell_i$. Let $1\le i<j\le 4$.
Since every induced path of $A$
between $a_i, a_j$ has the same length (because they can both be extended to a hole via $P_i, P_j$ and the same path of $B$ between $b_i, b_j$),
it follows that 
$$d_A(a_i,a_j)=\ell_i+\ell_j+d_{A'}(a_i',a_j').$$
Consequently
$$d(1,2,3,4)=\ell_1+\ell_2+\ell_3+\ell_4+ d'(1,2,3,4),$$
and similar equations hold for $d(1,3,2,4)$ and $d(1,4,2,3)$. Since every two of 
$$d'(1,2,3,4), d'(1,3,2,4), d'(1,4,2,3)$$ 
differ by at most one by \ref{closeddistance}, it follows that every two of 
$$d(1,2,3,4), d(1,3,2,4), d(1,4,2,3)$$ differ by at most one. This proves
\ref{dist}.~\bbox

\begin{thm}\label{border}
Let $\ell\ge 5$, let $G$ be $\ell$-holed, and let  $(A,B,P_1\LL P_k)$ be a $k$-bar gate in $G$. Then:
\begin{itemize}
\item either $|V(A)|=1$, or $|V(A)|\ge 3$ and $A$ is two-connected;
\item $V(A)=\{a_1\LL a_k\}$;
\item $A$ is a threshold graph;
\item at most one vertex of $A$ appears in the list $a_1\LL a_k$ more than once; that is, there do not exist
distinct $g,h,i,j\in \{1\LL k\}$ with $a_g=a_h\ne a_i=a_j$.
\end{itemize}
These statements also hold with $A$ replaced by $B$.
\end{thm}
\Proof
Suppose that $|V(A)|=2$, and $V(A)=\{u,v\}$ say. By \ref{legfact} we may assume that $a_1=a_2=u$ and $a_3=a_4=v$; and 
$(A, B, P_1,P_2,P_3,P_4)$ is a 4-bar semigate, contrary to \ref{dist}. This proves that $|V(A)|\ne 2$. If $|V(A)|=1$ then all four bullets
of the theorem hold; so we may assume that $|V(A)|\ge 3$. Suppose that there exist distinct
$g,h,i,j\in \{1\LL k\}$ such that
$$d_A(a_g,a_h)+d_A(a_i,a_j)+2\le  d_A(a_g,a_i)+d_A(a_h,a_j)= d_A(a_g,a_j)+d_A(a_h,a_i).$$
Then $(A,B,P_g,P_h,P_i,P_j)$ is a 4-bar semigate and \ref{dist} is violated. So there are no such $g,h,i,j$; and so by \ref{cutvertex}, 
this shows that $A$ has no cut-vertex, and hence it is 
two-connected. This proves the first bullet. 

If there exists $v\in V(A)\setminus \{a_1\LL a_k\}$, then since $A\setminus \{v\}$ is connected, it follows that
$(A\setminus \{v\}, B, P_1\LL P_k)$ is a strictly better $k$-bar semigate, a contradiction.  This proves the second bullet.

Suppose that $A$ contains the four-vertex path $\mathcal{P}_4$; say with vertices $a_1,a_2,a_3,a_4$ in order. Thus $d_A(a_1,a_2)+d_A(a_3,a_4)=2$,
and $d_A(a_1,a_3)+d_A(a_2,a_4) =4$, contrary to \ref{dist} applied to the 4-bar semigate $(A,B,P_1,P_2,P_3,P_4)$.
Suppose that $A$ contains $\overline {\mathcal{C}_4}$; then we may assume that $a_1a_2$ and $a_3a_4$ are edges, and 
$\{a_1,a_2\}$ and $\{a_3,a_4\}$ are anticomplete. Then $d_A(a_1,a_2)+d_A(a_3,a_4)=2$ and $d_A(a_1,a_3)+d_A(a_2,a_4) \ge 4$,
again contrary to \ref{dist}. Certainly $A$ does not contain $\mathcal{C}_4$ since $G$ is $\ell$-holed. Thus $A$ is a threshold graph.
This proves the third bullet.

Finally, suppose that $a_1=a_2\ne a_3=a_4$ say. Then $d_A(a_1,a_2)+d_A(a_3,a_4)=0$ and $d_A(a_1,a_3)+d_A(a_2,a_4) \ge 2$,
again contrary to \ref{dist}. This proves the fourth bullet, and so proves \ref{border}.~\bbox

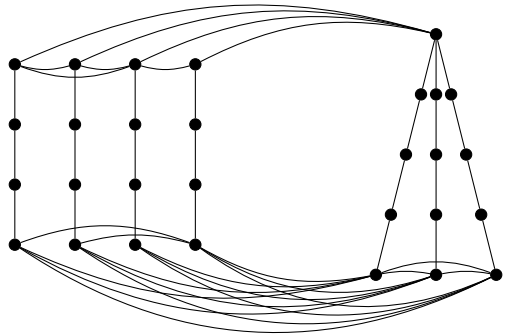
\begin{figure}[H]
\centering
\begin{tikzpicture}[scale=.8,auto=left]
\tikzstyle{every node}=[inner sep=1.5pt, fill=black,circle,draw]

\node (a1) at (0,3) {};
\node (a2) at (1,3) {};
\node (a3) at (2,3) {};
\node (a4) at (3,3) {};

\node (b1) at (0,0) {};
\node (b2) at (1,0) {};
\node (b3) at (2,0) {};
\node (b4) at (3,0) {};

\foreach \from/\to in {a1/a2,a2/a3,a3/a4}
\draw [bend right=15] (\from) to (\to);
\foreach \from/\to in {a1/a3}
\draw [bend right=20] (\from) to (\to);

\foreach \from/\to in {b1/b4}
\draw [bend left=20] (\from) to (\to);
\foreach \from/\to in {b2/b4}
\draw [bend left=15] (\from) to (\to);

\node (m1) at (0,2) {};
\node (n1) at (0,1) {};
\node (m2) at (1,2) {};
\node (n2) at (1,1) {};
\node (m3) at (2,2) {};
\node (n3) at (2,1) {};
\node (m4) at (3,2) {};
\node (n4) at (3,1) {};
\foreach \from/\to in {a1/m1,m1/n1,n1/b1,a2/m2,m2/n2,n2/b2,a3/m3,m3/n3,n3/b3,a4/m4,m4/n4,n4/b4}
\draw (\from) -- (\to);

\node (a0) at (7,3.5) {};
\node (b6) at (6,-.5) {};
\node (b7) at (7,-.5) {};
\node (b8) at (8,-.5) {};

\node (p6) at (6.75,2.5) {};
\node (q6) at (6.5,1.5) {};
\node (r6) at (6.25,0.5) {};
\node (p7) at (7,2.5) {};
\node (q7) at (7,1.5) {};
\node (r7) at (7,0.5) {};
\node (p8) at (7.25,2.5) {};
\node (q8) at (7.5,1.5) {};
\node (r8) at (7.75,0.5) {};
\foreach \from/\to in {a0/p6,p6/q6,q6/r6,r6/b6,a0/p7,p7/q7,q7/r7,r7/b7,a0/p8,p8/q8,q8/r8,r8/b8}
\draw (\from) -- (\to);
\foreach \from/\to in {b6/b7,b7/b8}
\draw [bend left=10] (\from) to (\to);
\foreach \from/\to in {b6/b8}
\draw [bend left=20] (\from) to (\to);

\foreach \to in {a1,a2,a3,a4}
\draw[bend right=20] (a0) to (\to);

\foreach \from in {b1,b2,b3,b4}
\draw[bend right=20] (\from) to (b6);

\foreach \from in {b1,b2,b3,b4}
\draw[bend right=25] (\from) to (b7);

\foreach \from in {b1,b2,b3,b4}
\draw[bend right=30] (\from) to (b8);

\end{tikzpicture}
\caption{A $9$-holed $7$-bar gate.} \label{fig:oddgate}
\end{figure}

We can say more about $k$-bar gates, but we need to treat the $\ell$ even and $\ell$ odd cases separately.
First:
\begin{thm}\label{oddgate}
Let $\ell\ge 5$ be odd, let $G$ be $\ell$-holed, and let  $(A,B,P_1\LL P_k)$ be a $k$-bar gate in $G$. Then, possibly after 
exchanging $A,B$, 
there exists $J\subseteq \{1\LL k\}$ with $|J|\ge 2$ such that all the vertices $(a_i:i\in J)$ are equal, equal to $a_0$ say,
and the following hold:
\begin{itemize}
\item the vertices $a_0$ and $a_i\;(i\in \{1\LL k\}\setminus J)$ are all distinct;
\item $b_1\LL b_k$ are all distinct;
\item the vertices $b_i\;(i\in J)$ are pairwise adjacent, and $P_i$ has length $(\ell-1)/2$ for each $i\in J$;
\item for each $i\in \{1\LL k\}\setminus J$, $a_0$ is adjacent to $a_i$, and $b_j$ is adjacent to $b_i$ for all $j\in J$, and
$P_i$ has length $(\ell-3)/2$;
\item for all distinct $i,j\in \{1\LL k\}\setminus J$, $a_i, a_j$ are adjacent if and only if $b_i, b_j$ are nonadjacent.
\end{itemize}
\end{thm}
(See figure \ref{fig:oddgate}.)

\Proof Let $P_i$ have length $\ell_i$ for each $i$. Since $A$ is a threshold graph (by \ref{border}), and so does not contain $\mathcal{P}_4$, and $A$ is connected,
it follows that every two vertices in $A$ have distance at most two in $A$; and the same in $B$.
For all $i,j\in \{1\LL k\}$, let us define $d(i,j)=d_A(a_i,a_j)+d_B(b_i,b_j)$.
We observe first:
\\
\\
(1) {\em $\ell_i+\ell_j+d(i,j)=\ell$ for $1\le i<j\le k$.}
\\
\\
Take a path of $A$ between $a_i,a_j$ of length $d_A(a_i,a_j)$, and a similar path in $B$; then their union with $P_i$ and $P_j$ 
is a hole of length $\ell$. This proves (1).
\\
\\
(2) {\em There do not exist distinct $i,j\in \{1\LL k\}$ with $a_i=a_j$ and $b_i=b_j$.}
\\
\\
Suppose that $a_1=a_2$ and $b_1=b_2$ say. Thus $\ell_1+\ell_2=\ell$.
From (1), and since $k\ge 3$, it follows that 
$\ell_i+\ell_3+d(i,3)=\ell$ for $i = 1,2$, and since $d(1,3)=d(2,3)$, it follows that $\ell_1=\ell_2$, and so 
$\ell=\ell_1+\ell_2$ is even, a contradiction. This proves (2).
\\
\\
(3) {\em $\ell_i\le (\ell+1)/2$ for $1\le i\le k$.}
\\
\\
Suppose that $\ell_1> (\ell+1)/2$ say, and so $\ell_1\ge (\ell+3)/2$. From (1) and (2), $\ell_1+\ell_i+1\le \ell$ for $i=2,3$; and hence 
$\ell_2+\ell_3\le \ell-5$. But since $d(2,3)\le 4$, (1) implies that $\ell_2+\ell_3\ge \ell-4$, a contradiction.
This proves (3).
\\
\\
(4) {\em $\ell_i\ge (\ell-5)/2$ for $1\le i\le k$.}
\\
\\
Suppose that $\ell_1\le (\ell-7)/2$ say. For $2\le i\le k$, since $\ell_i\le  (\ell+1)/2$ by (3), (1) implies that
$d(1,i)\ge 3$.
In particular, $a_i\ne a_1$. By \ref{legfact},
$a_1$ has at least two neighbours in $A$, say $a_2,a_3$; and since $d(1,i)\le 3$ for $i = 2,3$, it follows from (1)
that
$\ell_i\ge (\ell+1)/2$ for $i = 2,3$. But this
is impossible since $\ell_i+\ell_j\le \ell$ by (1). This proves (4).
\\
\\
(5) {\em $\ell_i\le (\ell-1)/2$ for $1\le i\le k$.}
\\
\\
Suppose that $\ell_1=(\ell+1)/2$ say. 
Suppose first that some vertex in $A$ has distance two from $a_1$; say $a_2$. Since $\ell_2\ge (\ell-5)/2$ by (4), (1) implies that
$d(1,2)\le 2$, and since $d_A(a_1,a_2)=2$ we deduce that $b_1=b_2$, and $\ell_2= (\ell-5)/2$.
Let $a_3$ be a vertex of $A$ adjacent to both $a_1,a_2$. By (1), $\ell_1+\ell_3=\ell-d(1,3)\le \ell-1$, and so $\ell_3\le (\ell-3)/2$.
But also from (1), $\ell_2+\ell_3=\ell-d(2,3)\ge \ell-3$, and so $\ell_3\ge (\ell-1)/2$, a contradiction.

This proves that $a_1$ is adjacent to every other vertex in $A$, and similarly $b_1$ is adjacent to every other vertex in $B$.
Suppose that $\ell_2= (\ell-5)/2$ say. From (1) it follows that $a_1\ne a_2$ and $b_1\ne b_2$.
By \ref{legfact} there exists $i\in \{3\LL k\}$ such that $a_i$ is equal or adjacent to $a_2$, and so $d(2,i)\le 3$. But
(1) implies that $\ell_2+\ell_i= \ell-d(2,i)$;  and from (1) again,
$\ell_1+\ell_i= \ell-d(1,i)$. Consequently $3=\ell_1-\ell_2=d(2,i)-d(1,i)$. But this is impossible since $d(2,i)\le 3$ and $d(1,i)\ge 1$.

Thus $\ell_i\ge  (\ell-3)/2$ for $2\le i\le k$. Since $d(1,i)\ge 1$ by (2), (1) implies that $\ell_1+\ell_i\le \ell-1$, and so
$\ell_i=(\ell-3)/2$ for $2\le i\le k$. From (1), $d(1,i)=1$, and so either $a_1=a_i$ or $b_1=b_i$, 
for each $i\in \{2\LL k\}$. Consequently $d(2,3)\le 2$, and so $\ell_2+\ell_3\ge \ell-2$ by (1), a contradiction. This proves (5).
\\
\\
(6) {\em $\ell_i\ge (\ell-3)/2$ for $1\le i\le k$.}
\\
\\
Suppose that $\ell_1=(\ell-5)/2$ say. If $a_i=a_1$ for some $i\in \{2\LL k\}$, then $d(1,i)\le 2$, and so (1) implies that 
$\ell_i\ge \ell-\ell_1-2\ge (\ell+1)/2$ contrary to (5). Choose a neighbour $a_2$ say of $a_1$. Since $d(1,2)\le 3$ and $\ell_2\le (\ell-1)/2$,
(1) implies that equality holds in both; that is, $b_1,b_2$ are nonadjacent and $\ell_2=(\ell-1)/2$. Choose a neighbour $b_3$ of $b_1$ say.
Then similarly $a_1,a_3$ are nonadjacent and $\ell_3=(\ell-1)/2$. So $a_2\ne a_3$ and $b_2\ne b_3$, and so $d(2,3)\ge 2$, contrary
to (1). This proves (6).

\bigskip
In summary, we have now shown that $\ell_i\in \{(\ell-3)/2,(\ell-1)/2\}$ for each $i\in \{1\LL k\}$. Let $I$ be the set
of all $i\in \{1 \LL k\}$ with $\ell_i=(\ell-3)/2$, and $J=\{1\LL k\}\setminus I$.
\\
\\
(7) {\em There do not exist $i\in I$ and $j\in J$ with $a_i=a_j$.}
\\
\\
Suppose that $1\in I$ and $2\in J$ and $a_1=a_2$ say. By (1), $d_B(b_1,b_2)=2$; choose $b_3\in B$ adjacent to both $b_1,b_2$.
Then $d(1,3)=d(2,3)$, and yet $\ell_1\ne \ell_2$, contrary to (1).
This proves (7).
\\
\\
(8) {\em If $i,j\in I$ are distinct then $P_i, P_j$ are vertex-disjoint, and $a_i,a_j$ are adjacent if and only if $b_i, b_j$
are nonadjacent.}
\\
\\
By (1), $d(i,j)=3$, and both the statements follow. This proves (8).
\\
\\
(9) {\em If $i\in I$ and $j\in J$ then $a_ia_j$ and $b_ib_j$ are edges.}
\\
\\
From (7), $a_i\ne a_j$ and $b_i\ne b_j$; but by (1), $d(i,j)=2$. This proves (9).
\\
\\
(10) {\em $|J|\ge 2$.}
\\
\\
Suppose that $J$ is empty. Since $A,B$ are threshold graphs and hence cographs, and one is isomorphic to the complement of the other by (8),
it follows that one of $A,B$ is disconnected, a contradiction. Thus $J\ne \emptyset$. Now suppose that $J=\{1\}$ say.
Again, the graphs $A\setminus \{a_1\}$, $B\setminus \{b_1\}$ are complementary threshold graphs, with at least two vertices since 
$k\ge 3$, and so one of them is disconnected; and so one of $A,B$ has a cut-vertex, contrary to \ref{border}. This proves (10).
\\
\\
(11) {\em Either all the vertices $a_j\;(j\in J)$ are equal or all the vertices $b_j\;(j\in J)$ are equal.}
\\
\\
In a bipartite graph, either all edges have a common end, or some two of them are disjoint. Thus, if the claim is false, then
there exist distinct $i,j\in J$ with $a_i\ne a_j$ and $b_i\ne b_j$, and so $d(i,j)\ge 2$,
contrary to (1). This proves (11).

\bigskip

From (1), (2), (7), (10) and (11), this proves \ref{oddgate}.~\bbox

These graphs look like they are getting complicated, but there is a better way to think of them. Remove the vertex $a_0$ from $A$, and add to $A$ all the 
neighbours of $a_0$ instead, forming $A'$. All the $A'-B$ paths have the same length and are pairwise vertex-disjoint. 
With $a_0$ deleted, this has become what we called an $\ell$-frame in the previous 
section; the graph $A'$ is not connected, but the extra vertex $a_0$ is adjacent to every vertex in $A'$, and $B$ is two-connected by the first 
bullet of \ref{border}.

We can say this more precisely as follows. Let $F$ be an $\ell$-frame, with $\ell$ odd, and with sides $A,B$ and bars $P_1\LL P_t$.
The graphs $A,B$ are 
complementary threshold graphs, and so exactly one of them is disconnected. If $A$ is disconnected, add a new vertex $a_0$
adjacent to every vertex 
of $A$, and with no other neighbours.
Let us call this new vertex the {\em apex}. Let $A^+$ be the subgraph induced on $V(A)\cup \{a_0\}$, and let $B^+=B$.
If $B$ is disconnected, add a new vertex $b_0$
adjacent to every vertex 
of $B$, and with no other neighbours;
and we call this new vertex the {\em apex}. Let $B^+$ be the subgraph induced on $V(B)\cup \{b_0\}$, and let $A^+=A$.

The graph obtained from $F$ by adding the apex is called an {\em apexed $\ell$-frame}, we call $P_1\LL P_t$ its {\em bars}, and 
$A^+,B^+$ its {\em sides}. 
From \ref{border} and \ref{oddgate}, it follows that:
\begin{thm}\label{oddframe}
Let $\ell\ge 5$ be odd, let $G$ be $\ell$-holed, and let  $(A',B,P_1'\LL P_k')$ be a $k$-bar gate in $G$. 
Let $P_i'$ have ends $a_i', b_i$ for $1\le i\le k$, where $a_i'\in V(A')$. Then, possibly after exchanging $A',B$, there exists $a_0\in V(A')$ with the following property.
For $1\le j\le k$, if $a_j'\ne a_0$ let $a_j=a_j'$ and $P_j=P_j'$, and if $a_j'=a_0$ let $a_j$ be the neighbour of $a_j'$ in $P_j$, and let 
$P_j=P_j'\setminus \{a_j'\}$.
Let $A^+$ be the subgraph of $G$ induced on $\{a_0,a_1\LL a_k\}$, and let $B^+=B$.
Then
$$A'\cup B\cup P_1'\cupcup P_k'=A^+\cup B^+\cup P_1\cupcup P_k$$ 
is an apexed $\ell$-frame, with sides $A^+,B^+$, bars $P_1\LL P_k$ and apex $a_0$.
\end{thm}

\begin{figure}[H]
\centering
\begin{tikzpicture}[scale=.8,auto=left]
\tikzstyle{every node}=[inner sep=1.5pt, fill=black,circle,draw]

\node (a1) at (0,2) {};
\node (a2) at (1,2) {};
\node (a3) at (2,2) {};

\node (b1) at (0,0) {};
\node (b2) at (1,0) {};
\node (b3) at (2,0) {};

\node (n1) at (0,1) {};
\node (n2) at (1,1) {};
\node (n3) at (2,1) {};

\foreach \from/\to in {a1/n1,n1/b1,a2/n2,n2/b2,a3/n3,n3/b3}
\draw (\from) -- (\to);

\node (a4) at (3,3) {};
\node (a5) at (4,3) {};
\node (a6) at (5,3) {};

\node (n4) at (3,1.67) {};
\node (n5) at (4,1.67) {};
\node (n6) at (5,1.67) {};

\node (o4) at (3,0.33) {};
\node (o5) at (4,0.33) {};
\node (o6) at (5,0.33) {};

\node (b4) at (3,-1) {};
\node (b5) at (4,-1) {};
\node (b6) at (5,-1) {};

\foreach \from/\to in {a4/n4,n4/o4,o4/b4,a5/n5,n5/o5,o5/b5,a6/n6,n6/o6,o6/b6}
\draw (\from) -- (\to);

\foreach \from/\to in {a4/a6}
\draw [bend right=20] (\from) to (\to);
\foreach \from/\to in {a4/a5,a5/a6}
\draw [bend right=15] (\from) to (\to);

\foreach \from/\to in {b4/b6}
\draw [bend right=20] (\from) to (\to);
\foreach \from/\to in {b4/b5, b5/b6}
\draw [bend right=15] (\from) to (\to);

\foreach \to in {a1,a2,a3}
\draw (a4) to (\to);

\foreach \from in {b1,b2,b3}
\draw (\from) to (b6);

\foreach \from in {a2,a3}
\draw (\from) to (a5);

\foreach \from in {b1}
\draw (\from) to (b5);

\node (a0) at (-2,2.5) {};

\node (n7) at (-2.17,1.5) {};
\node (n8) at (-1.83,1.5) {};

\node (o7) at (-2.33,0.5) {};
\node (o8) at (-1.63,0.5) {};

\node (b7) at (-2.5,-.5) {};
\node (b8) at (-1.5,-.5) {};

\foreach \from/\to in {a0/n7,n7/o7,o7/b7,a0/n8,n8/o8,o8/b8}
\draw (\from) -- (\to);

\node (b0) at (-4,-.5) {};
\node (o9) at (-3.833,0.5) {};
\node (o10) at (-4.167,0.5) {};

\node (n9) at (-3.67,1.5) {};
\node (n10) at (-4.333,1.5) {};

\node (a9) at (-3.5,2.5) {};
\node (a10) at (-4.5,2.5) {};

\foreach \from/\to in {a9/n9,n9/o9,o9/b0,a10/n10,n10/o10,o10/b0}
\draw (\from) -- (\to);

\foreach \from in {a9,a10}
\draw[bend left=20] (\from) to (a0);

\foreach \to in {b7,b8}
\draw[bend right=20] (b0) to (\to);

\foreach \to in {a1,a2,a3,a4,a5,a6}
\draw[bend left=20] (a0) to (\to);

\foreach \to in {a4,a5,a6}
\draw[bend left=25] (a9) to (\to);

\foreach \to in {a4,a5,a6}
\draw[bend left=30] (a10) to (\to);

\foreach \to in {b1,b2,b3,b4,b5,b6}
\draw[bend right=30] (b0) to (\to);

\foreach \to in {b4,b5,b6}
\draw[bend right=25] (b7) to (\to);

\foreach \to in {b4,b5,b6}
\draw[bend right=30] (b8) to (\to);

\draw[dotted] (-5,0.125) -- (6, 0.125);
\draw[dotted] (-5,1.875) -- (6, 1.875);
\draw[dotted] (-3,4.5) -- (-3, -2);
\draw[dotted] (-1,4.5) -- (-1, -2);
\draw[dotted] (2.5,4.5) -- (2.5, -2);
\tikzstyle{every node}=[-]
\draw node at (-5.5, 2.5)           {$A$};
\draw node at (-5.5, -.5)           {$B$};
\draw node at (-4, 4.5)           {$L$};
\draw node at (-2, 4.5)           {$K$};
\draw node at (1, 4.5)           {$I$};
\draw node at (4, 4.5)           {$J$};

\end{tikzpicture}
\caption{An $8$-holed $10$-bar gate.} \label{fig:evengate}
\end{figure}
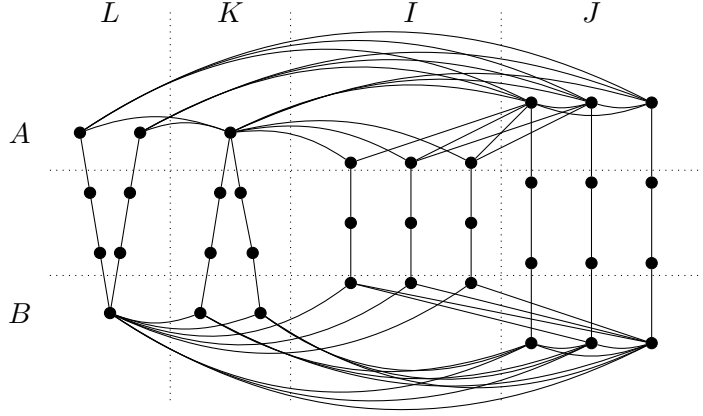

For the $\ell$ even case, we have:
\begin{thm}\label{evengate}
Let $\ell\ge 6$ be even, let $G$ be $\ell$-holed, and let  $(A,B,P_1\LL P_k)$ be a $k$-bar gate in $G$. Then
either $|A|=|B|=1$ and $P_1\LL P_k$ all have length $\ell/2$, or
there is a partition $(I,J,K,L)$ of $\{1\LL k\}$ into four sets with the following properties:
\begin{itemize}
\item $J,K, L\ne \emptyset$ (possibly $I=\emptyset$);
the vertices $a_i\;(i\in K)$ are all equal, with common value $a_0$ say; the vertices $a_i\;(i\in I\cup J\cup L\cup \{0\})$ are all 
distinct; the vertices $b_i\;(i\in L)$ are all equal, with common value $b_0$ say; and the vertices $b_i\;(i\in I\cup J\cup K\cup \{0\})$ are 
all distinct;
\item $\{a_i:i\in I\cup L\}$ and $\{b_i:i\in I\cup K\}$ are stable sets;
\item $\{a_i:i\in J\cup K\}$ and $\{b_i:i\in J\cup L\}$ are cliques;
\item $a_0$ is adjacent to $a_i$ for all $i\in I\cup J\cup L$, and $b_0$ is adjacent to $b_i$ for all $i\in I\cup J\cup K$;
\item the graph
$G[\{a_i:i\in I\cup L\},\{a_j:j\in J\cup K\}]$ is a half-graph; and for each $i\in I\cup L$ and each $j\in J\cup K$, $a_i,a_j$ are adjacent if and only if $b_i,b_j$ are nonadjacent;  and 
\item $P_i$ has  length $\ell/2-2$ for each $i\in I$, and $P_i$ has length $\ell/2-1$ for each $i\in J\cup K\cup L$.
\end{itemize}
\end{thm}
(See figure \ref{fig:evengate}.)

\Proof The proof is very much like that of \ref{oddgate}, and we use the same notation without redefining it. 
In particular, any two vertices in $A$ have distance at most two in $A$, and the same in $B$.
We still have (with the same proof):
\\
\\
(1) {\em $\ell_i+\ell_j+d(i,j)=\ell$ for $1\le i<j\le k$.}
\\
\\
(2) {\em $\ell_i\le \ell/2$ for $1\le i\le k$.}
\\
\\
Suppose that $\ell_1\ge \ell/2+1$ say. For $2\le i\le k$, if $a_i=a_1$ and $b_i=b_1$ then $\ell_i=\ell/2-1$, and so there is at most one
such value of $i$. Since $k\ge 3$ it follows that one of $|A|,|B|>1$, and we assume that $|A|>1$ without loss of generality.
Since $A$ is therefore two-connected by \ref{border}, $a_1$ has at least two neighbours in $A$, and so we may assume that $a_2,a_3$
are both adjacent to $a_1$. By (1), $\ell_1+\ell_i+d(1,i) = \ell$, and so $\ell_i\le \ell/2-1-d(1,i)$ for $i=2,3$. But $\ell_2+\ell_3= \ell-d(2,3)$
by (1), so $d(1,2)+d(1,3)\le  d(2,3)-2$, which is impossible by the triangle inequality. This proves (2).
\\
\\
(3) {\em $\ell_i\ge \ell/2-2$ for $1\le i\le k$.}
\\
\\
Suppose that $\ell_1\le \ell/2-3$ say. For $2\le i\le k$, (1) implies that $\ell_i\ge \ell/2+3-d(1,i)$, and so by (2), $d(1,i)\ge 3$
for $2\le i\le k$. In particular, $a_i\ne a_1$, and so by \ref{border} we may assume that $a_2,a_3$ are distinct neighbours of $a_1$.
Since $d(1,i)\le 3$ for $i = 2,3$ (because $a_1,a_i$ are adjacent and any two vertices in $B$ have distance at most two in $B$), (1) implies that $\ell_i=\ell/2$ for $i = 2,3$; but $\ell_2+\ell_3<\ell$ by (1) since $a_2\ne a_3$,
a contradiction. This proves (3).
\\
\\
(4)  {\em We may assume that $\ell_i\le \ell/2-1$ for $1\le i\le k$.}
\\
\\
Suppose that $\ell_1\ge \ell/2$ say, and so $\ell_1=\ell/2$ by (2). Suppose first that some vertex $a_2\in V(A)$ is different from and nonadjacent to $a_1$. Thus $d(1,2)\ge 2$; but $\ell_2\ge \ell/2-2$ by (3), and so equality holds in both by (1). In particular $b_2=b_1$ and $\ell_2=\ell/2-2$. Choose
$a_3$ adjacent to both $a_1,a_2$. 
Then $d(1,3)=d(2,3)$, and yet $\ell_1\ne \ell_2$, contrary to (1).
Thus $a_1$ is adjacent to every other vertex in $A$, and the same for $b_1$ in $B$.

Suppose that $\ell_2=\ell/2-2$ say. By (1), $d(1,2)=2$, so $a_2\ne a_1$ and $b_2\ne b_1$. Since $A$ is two-connected by \ref{border}, 
there exists $a_3\ne a_1,a_2\in V(A)$
adjacent to $a_2$; and also adjacent to $a_1$ since $a_1$ is adjacent to all other vertices in $A$. Thus $\ell_3= \ell/2-d(1,3)$ by (1),
and $\ell_3=\ell/2+2-d(2,3)$ by (1); and hence $d(2,3)=d(1,3)+2$. But $d(2,3)\le 3$ and $d(1,3)\ge 1$, so equality holds in both.
Since $d(1,3)=1$, it follows that $b_3=b_1$, and so $b_3,b_2$ are adjacent, and so $d(2,3)=2$, a contradiction.
This proves that $\ell_i\ge \ell/2-1$, and hence $d(1,i)\le 1$, for $2\le i\le k$. 

Suppose that $|A|>1$; then we may choose adjacent $a_2,a_3$ both different from $a_1$, by \ref{border}; and both are adjacent to $a_1$
since $a_1$ is adjacent to all the other vertices of $A$. Since $d(1,i)\le 1$ for $i = 2,3$, it follows that $b_2=b_3=b_1$; but then
the union of $P_1,P_2,P_3$ and the three edges $a_1a_2,a_2a_3,a_3a_1$ is a pyramid (because $P_1$ has length $\ell/2\ge 3$ and
$P_2,P_3$ each have length at least $\ell/2-1\ge 2$), and so $G$ has an odd hole, a contradiction.

This proves that $|A|=1$ and similarly $|B|=1$, and so $\ell_i=\ell/2$ for all $i\in \{1\LL k\}$ and the theorem holds. This proves (4).

\bigskip
So $\ell_i\in \{\ell/2-1,\ell/2-2\}$ for $1\le i\le k$. Let $I$ be the set of $i\in \{1\LL k\}$ with $\ell_i=\ell/2-2$. 
From (4), there do not exist distinct $i,j\in \{1\LL k\}$ with $d(i,j)\le 1$; and in particular, either $a_i\ne a_j$
or $b_i\ne b_j$ for all distinct $i,j$.
\\
\\
(5) {\em If $i\in I$ and $j\in \{1\LL k\}\setminus I$ then $a_i\ne a_j$ and $b_i\ne b_j$; and $a_i,a_j$ are adjacent if and only
if $b_i,b_j$ are nonadjacent.}
\\
\\
By (1), $d(i,j)= 3$; and both claims follow. This proves (5).
\\
\\
(6) {\em If $i,j\in I$ are distinct, then $a_i,a_j$ are distinct and nonadjacent, and $b_i, b_j$ are distinct and nonadjacent. If $i,j\in \{1\LL k\}\setminus I$
are distinct, then either $a_ia_j, b_ib_j$ are both edges, or $a_i=a_j$ and $b_i,b_j$ are distinct and nonadjacent, or $b_i=b_j$ and $a_i,a_j$ are
distinct and nonadjacent.}
\\
\\
If $i,j\in I$ are distinct, then $\ell_i+\ell_j=\ell-4$, and so (1) implies that $d(i,j)=4$. If $i,j\in  \{1\LL k\}\setminus I$ are distinct,
then $\ell_i+\ell_j=\ell-2$, and so (1) implies that $d(i,j)=2$. This proves (6).
\\
\\
(7) {\em $|A|,|B|>1$.}
\\
\\
Suppose that $|B|=1$, say. Since $\ell_i\le \ell/2-1$ for each $i$, (1) implies that $d(i,j)\ge 2$ for all distinct 
$i,j\in \{1 \LL k\}$; and since $b_1=\cdots=b_k$, it follows that $a_1\LL a_k$ are all distinct and nonadjacent,
contradicting that $A$ is connected.
This proves (7).

\bigskip
If there exist distinct $i,j\in \{1\LL k\}$ with $a_i=a_j$, let $a_0 = a_i$, and let $K$ be the set of all $h \in \{1\LL k\}$
with $a_h=a_0$. If there are no such $i,j$, 
since $A$ is a connected threshold graph, it has a vertex (indeed, at least two such vertices, since it is two-connected) adjacent 
to all other vertices in $A$. Choose some such vertex $a_i$, and let $K=\{i\}$ and define $a_0=a_i$. Define $b_0, L$ similarly.
In the case when no two of the paths $P_1\LL P_k$ share an end, it follows that $|K|=|L|=1$, and in this case there are two
choices for $b_0$, since at least two vertices in $B$ are adjacent to all the others; so in this case we can additionally
choose $b_0$ such that $K\ne L$. In any case, $|\{a_i:i\in K\}|=1$, and from the final statement of \ref{border}, 
the vertices $a_i; (i\in \{1\LL k\}\setminus K)$ are all
distinct; and similarly $|\{b_i:i\in L\}|=1$, and $b_i\;(i\in \{1\LL k\}\setminus L)$ are all distinct.
\\
\\
(8) {\em $\ell_i=\ell/2-1$ for each $i\in K$; $a_0$ is adjacent to all other vertices in $A$; and the set $\{b_i:i\in I\cup K\}$
is stable. Similarly $b_0$ is adjacent to all other vertices in $B$, and $\{a_i:i\in I\cup L\}$ is stable.}
\\
\\
Suppose first that $|K|\ge 2$; then $\ell_i=\ell/2-1$ for each $i\in K$, by (5) and (6). Suppose that $a_0$
is nonadjacent to some other vertex in $A$; $a_0$ is not adjacent to $a_1$ say, where $a_1\ne a_0$. We may assume
that $a_2$ is adjacent to both $a_0, a_1$. Since $|K|\ge 2$,
we may assume that $3,4\in K$. Since $b_3\ne b_4$, we may assume that $b_3\ne b_1$;
then by (6), it follows that
$1\in I$. By (5), $b_1\ne b_4$, and $b_1$ is adjacent to $b_3, b_4$. 
By (6), $2\notin I$, and $b_2$ is different from $b_3,b_4$;
by (6) again, $b_2$ is adjacent to $b_3, b_4$; and by (5), $b_1,b_2$ are nonadjacent. 
Also $b_3,b_4$ are distinct and nonadjacent by (6), and $\{b_1,b_2,b_3,b_4\}$ induces a 4-hole, a contradiction.
So in this case, $a_0$ is adjacent to all other vertices in $A$; and $\{b_i:i\in K\}$ is stable by (6), and 
$\{b_i:i\in I\}$ is stable by (6), and $\{b_i:i\in K\}$ and $\{b_i:i\in I\}$ are anticomplete by (5). So in this case
(8) holds.

Now we assume that $|K|=1$, $K=\{1\}$ say, and so $a_1\LL a_k$ are all distinct. It follows that $a_1$ is adjacent to all other vertices in $A$ from the definition of $a_0$.
Suppose that $1\in I$. Then $I=\{1\}$ by (6). Choose $b_2$ say in $B$, adjacent to $b_1$ (this is possible by (7) and since $B$ is connected). 
Then 
$a_1a_2$ and $b_1b_2$ are both edges, contrary to (5). Thus $1\notin I$, and so again (8) holds. This proves (8).
\\
\\
(9) {\em $K\cap L=\emptyset$.}
\\
\\
Suppose that $1\in K\cap L$. If $|K|\ge 2$ and $2\in K$ say, then $b_1,b_2$ are distinct and nonadjacent by (6), and yet
$b_1=b_0$ is adjacent to all other vertices in $B$ by (8), a contradiction. So $|K|=|L|=1$; but in this case we were careful 
to choose $b_0,L$
so that $K\ne L$, a contradiction. This proves (9).

\bigskip
Let $J=\{1\LL k\}\setminus (I\cup K\cup L)$.
\\
\\
(10) {\em $J,K,L\ne \emptyset$.}
\\
\\
From their definitions, $K,L\ne \emptyset$. Let $1\in L$ say. By (8), $a_1$ has no neighbour in $\{a_i:i\in I\cup L\}$, and
has only one in $\{a_i:i\in K\}$ since all these vertices are equal; and since $a_1$ has degree at least two in $A$ (from (7) and \ref{border}), it has a neighbour
in $\{a_i:i\in J\}$. Thus $J\ne \emptyset$. This proves (10).
\\
\\
From (6) and since the vertices $a_i\:(i\in J)$ and $b_i\;(i\in J)$ are all distinct, it follows that $\{a_i:i\in J\}$ is a clique, 
and therefore so is $\{a_i:i\in J\cup K\}$; and similarly $\{b_j:i\in J\cup L\}$ is a clique.
Since $A$ is a threshold graph by \ref{border}, 
it follows that $G[\{a_i:i\in I\cup L\},\{a_j: j\in J\cup K\}]$ is a half-graph.
From \ref{border}, (5), (6), (8), (9) and (10), this proves \ref{evengate}.~\bbox

Again, there is a better way to think of these graphs. Let $A'$ be obtained from $A$ by deleting $a_0$ and adding the neighbours 
of $a_0$ in each path $P_i\; (i\in K)$, and define $B'$ similarly. Then all the $A'-B'$ paths are pairwise vertex-disjoint, 
and with $a_0, b_0$ deleted, it becomes an $\ell$-frame; and $a_0$ is adjacent to every vertex in $A'$, and 
$b_0$ to every vertex in $B'$.

More precisely, let $\ell\ge 6$ be even, let $F$ be an $\ell$-frame, with sides $A,B$ and bars $P_1\LL P_k$.
Add two new vertices $a_0, b_0$, where $a_0$ is adjacent to 
each vertex of $A$, and $b_0$ adjacent to each vertex of $B$. Let $A^+$ be the subgraph induced on $V(A)\cup \{a_0\}$, 
and define $B^+$ similarly. The enlarged graph we produce is called
an {\em apexed $\ell$-frame}, the two new vertices are called its {\em apexes}, $P_1\LL P_k$ are {\em bars}, and $A^+,B^+$  are its
{\em sides}.
From \ref{border} and \ref{evengate} we have:

\begin{thm}\label{evenframe}
Let $\ell\ge 6$ be even, let $G$ be $\ell$-holed, and let  $(A',B',P_1'\LL P_k')$ be a $k$-bar gate in $G$. 
Let $P_i'$ have ends $a_i', b_i'$ for $1\le i\le k$, where $a_i'\in V(A')$. Then there exist $a_0\in V(A')$
and $b_0\in V(B')$ with the following properties. For $1\le j\le k$, if $a_j'\ne a_0$ let $a_j=a_j'$ and otherwise let 
$a_j$ be the neighbour of $a_j'$ in $P_j$; and if $b_j'\ne b_0$ let $b_j=b_j'$ and otherwise let 
$b_j$ be the neighbour of $b_j'$ in $P_j$; let $P_j$ be the subpath of $P_j'$ between $a_j, b_j$. Let $A^+$ be the subgraph 
induced on $a_0$ and all its neighbours, and define $B^+$ similarly. Then 
$$A'\cup B'\cup P_1'\cupcup P_k'=A^+\cup B^+\cup P_1\LL P_k$$ 
is an apexed $\ell$-frame, with sides $A^+,B^+$, bars $P_1\LL P_k$ and apexes $a_0, b_0$.
\end{thm}
We remark that, whether $\ell$ is odd or even, if $F$ is an apexed $\ell$-frame, then in the usual notation, if $a_0$ exists then $A$ is disconnected,
and if $a_0$ does not exist then $A$ is two-connected; and the same for $B,b_0$.

\section{Attachments to a maximal frame}

First, we remark that in previous sections, $F$ was an $\ell$-frame, but now it will be an apexed $\ell$-frame.
In this section we prove the following:

\begin{thm}\label{framejumps}
Let $\ell\ge 7$, and let $G$ be an $\ell$-holed graph, containing an apexed $\ell$-frame. Choose an apexed $\ell$-frame $F$ in $G$,
with as many bars as possible, with sides $A^+,B^+$ and bars $P_1\LL P_k$. Let $a_0$ be the vertex of $A^+$ not in $P_1\cupcup P_k$, if there is one, and define $b_0$ similarly.
Then either
\begin{itemize}
\item there exist $x\in V(G)\setminus V(F)$ and $u,v\in V(F)$, such that $u,v$ are nonadjacent, both are adjacent to $x$, and
not both are in $V(A^+)$, and not both are in $V(B^+)$; or
\item there exists $x\in V(G)\setminus V(F)$ such that either $a_0$ exists and the set of neighbours of $x$ in $V(F)$ is $V(A^+)$, or
$b_0$ exists and the set of neighbours of $x$ in $V(F)$ is $V(B^+)$;
or
\item for every component $C$ of $G\setminus V(F)$, 
either $N(C)\subseteq V(A^+)$, or $N(C)\subseteq V(B^+)$, or $N(C)$ is a clique, where $N(C)$ denotes the set of vertices in $V(F)$ with a neighbour 
in $V(C)$.
\end{itemize}
\end{thm}
\Proof Let $A=A^+\setminus \{a_0\}$ if $a_0$ exists, and let $A^+=A$ otherwise; and define $B$ similarly.
We observe first that:
\\
\\
(1) {\em If $X\subseteq V(F)$, such that $X$ is not a clique and $X\not\subseteq V(A^+)$ and $X\not\subseteq V(B^+)$, then
there exist nonadjacent $u,v\in X$, not both in $V(A^+)$ and not both in $V(B^+)$.}
\\
\\
Let $u,v\in X$ be nonadjacent.
If $u$ or $v$ belongs to $V(F)\setminus (V(A^+)\cup V(B^+))$, then $u,v$ satisfy (1);
so we assume that $u,v\in V(A^+)\cup V(B^+)$. If 
$u\in V(A^+)$ and $v\in V(B^+)$, again (1) holds; so we may assume that $u,v\in V(A^+)$. Choose $x\in X\setminus V(A^+)$. It
follows that $x$ is nonadjacent to one of $u,v$, say $v$; but then $x,v$ satisfy (1). This proves (1).

\bigskip

We may assume that the first bullet is false, and so, from (1):
\\
\\
(2) {\em For each $x\in V(G)\setminus V(F)$, the set of vertices of $F$ adjacent to $x$ is either a clique, or a subset of $V(A^+)$, or
a subset of $V(B^+)$.}

\bigskip

We may also assume that the third bullet is false, and so by (1),
there is a component $C$ of $G\setminus V(F)$ and nonadjacent $u,v\in N(C)$, such that 
$\{u,v\}\not\subseteq V(A^+)$ and $\{u,v\}\not\subseteq V(B^+)$.
Consequently there is a minimal connected induced subgraph $S$ of $G\setminus V(F)$, such that there are two nonadjacent
vertices in $V(F)$ both with neighbours in $V(S)$, and not both in $V(A^+)$ and not both in $V(B^+)$. From the minimality
of $S$, it follows that $S$ is an induced path $s_1\CC s_n$ say. Let $N(S)$ denote the set of 
vertices of $F$ with a neighbour in $V(S)$.
Thus $N(S)$ is not a clique, and not a subset of $V(A^+)$, and not a subset of $V(B^+)$.
By (2), $n\ge 2$.
\\
\\
(3) {\em One of $s_1,s_n$ has a neighbour in $V(F)\setminus (V(A^+)\cup V(B^+))$.}
\\
\\
Suppose that all neighbours of $s_1$ in $V(F)$ and all neighbours of $s_n$ in $V(F)$ belong to  $V(A^+)\cup V(B^+)$. 
Since neither of $s_1,s_n$ has a neighbour in $V(A^+)$ and one in $V(B^+)$, we may assume that 
all neighbours of $s_1$ in $V(F)$ belong to $V(A^+)$, and all neighbours of $s_n$ in $V(F)$ belong to $V(B^+)$.

Suppose first that some vertex $w\in V(F)$ is adjacent to an internal vertex of $S$. From the minimality of $S$, 
some $P_i$,
say $P_1$, has length two, and $w$ is its middle vertex, and $a_1$ is the only neighbour of $s_1$ in $V(F)$, and $b_1$ is the only
neighbour of $s_n$ in $V(F)$. There is an induced path $R$ between $a_1,b_1$ consisting of a one- or two-edge path of $A^+$,
the path $P_2$, and a one- or two-edge path of $B^+$; and the union of $R$ and the path $a_1\DD s_1\DD s_2\CC s_n\DD b_1$ 
is a hole.
Also the union of $R,P_1$ is a hole; so the paths $a_1\DD s_1\DD s_2\CC s_n\DD b_1$ and $P_1$ have the same length, and hence $n=1$,
a contradiction.

Thus no internal vertex of $S$ has a neighbour in $V(F)$. Let $A'$ be the subgraph induced on $V(A^+)\cup \{s_1\}$,
and define $B'$ similarly; then 
$$(A', B', S, P_1\LL P_k)$$
is a $(k+1)$-bar semigate, and so $G$ contains a $(k+1)$-bar gate, and hence contains an apexed $\ell$-frame with $k+1$ bars by 
\ref{oddframe} and \ref{evenframe}, contrary to the maximality of $k$. 
This proves (3).

\bigskip

From (3) we may assume that some neighbour of $s_1$ is an internal vertex of $P_1$ say. By (2), all neighbours of $s_1$ in $V(F)$
belong to $V(P_1)$, and either there is just one such neighbour, or there are two and they are adjacent. From the minimality of $S$,
every vertex in $V(F)$ with a neighbour in $V(S)\setminus \{s_n\}$ belongs to $V(P_1)$. 
Let $P_1$ have vertices $q_1\CC q_m$ in order, where $q_1=a_1$ and $q_m=b_1$,
and $\ell/2-1\le m\le \ell/2$ (because $P_1$ has length between $\ell/2-2$ and $\ell/2-1$).
\\
\\
(4) {\em Some neighbour of $s_n$ in $V(F)$ does not belong to $V(P_1)$.}
\\
\\
Suppose that all neighbours of $s_n$ in $V(F)$ belong to $V(P_1)$. Thus $N(S)\subseteq V(P_1)$. By (2), 
$s_n$ has at most two neighbours in $V(F)$ and they form a clique.
Choose $g,h\in \{1\LL m\}$ minimum
and maximum respectively such that $q_g,q_h\in N(S)$. Thus $h\ge g+2$, since $N(S)$ is not a clique. 
From the minimality of $S$, no internal vertex of $S$ is adjacent to $q_g$ or to $q_h$, and so we may assume that
$s_1$ is adjacent to $q_g$ and $s_n$
to $q_h$.
There is a hole $C$ containing $P_1\cup P_2$; and by replacing the path $q_g\CC q_h$ of this hole by the 
path $q_g\DD s_1\CC s_{n}\DD q_h$, we obtain another hole, which therefore has the same length. Thus $h-g=n+1$. Since $n\ge 2$, it
follows that $h-g\ge 3$. From the minimality of $S$, and because $h-g\ge 3$, it follows that
none of $s_2\LL s_{n-1}$ have a neighbour in $V(F)$. If $s_1$ is adjacent to $q_{g+1}$ let $g'=g+1$, and otherwise let $g'=g$;
and if $s_n$ is adjacent to $q_{h-1}$ let $h'=h-1$, and otherwise let $h'=h$. Thus $h'> g'$. There is a hole
$$s_1\CC s_n\DD q_{h'}\DD q_{h'-1}\CC q_{g'}\DD s_1,$$ 
and its length is $n+1+h'-g'$; and so $n+1+h'-g'=\ell$. But $h-g=n+1$, and so $h-g+h'-g'=\ell$, and therefore $h-g\ge \ell/2$.
But $h-g$ is at most the length of $P_1$, which is at most  $\ell/2-1$, a contradiction. This proves (4).
\\
\\
(5) {\em No neighbour of $s_n$ is an internal vertex of any of $P_1\LL P_k$.}
\\
\\
By (4) and (2), no neighbour of $s_n$ is an internal vertex of $P_1$.
Suppose that $s_n$ is adjacent to some internal vertex of $P_2$, say. By (2), 
all neighbours of $s_2$ in $V(F)$ belong to $V(P_2)$ and form a clique; and from the minimality of $S$, none
of $s_2\LL s_{n-1}$ has any neighbour in $V(F)$. 
Let $P_2$ have vertices 
$r_1\CC r_s$ in order, where $r_1=a_2$, and $\ell/2-1\le s\le \ell/2$. Choose $g\in \{1\LL m\}$ minimum such that 
$s_1,q_g$ are adjacent, and if $s_1$ is adjacent to $q_{g+1}$ let $g'=g+1$, and otherwise let $g'=g$. Similarly choose 
$h\in \{1\LL s\}$ minimum such that $s_n,r_h$ are
 adjacent, and if $s_n$ is adjacent to $r_{h+1}$ let $h'=h+1$, and otherwise let $h'=h$. 
Since $A^+$ is a connected threshold graph, there is an induced path $T$ of $A^+$ of length at most two 
between $a_1,a_2$. Since $n\ge 2$, the cycle
$$q_1\CC q_g\DD s_1 \CC s_n\DD r_h\DD r_{h-1}\DD r_1\DD T\DD q_1$$
is a hole, so the path 
$$q_1\CC q_g\DD s_1 \CC s_n\DD r_h\DD r_{h-1}\CC r_1$$
has length at least $\ell-2$. Similarly the path
$$q_m\CC q_{g'}\DD s_1 \CC s_n\DD r_{h'}\CC r_s$$
has length at least $\ell-2$.
Consequently one of the paths 
$$q_1\CC q_g\DD s_1 \CC s_n\DD r_{h'}\CC r_s$$
$$q_m\CC q_{g'}\DD s_1 \CC s_n\DD r_{h}\CC r_1$$
has length at least $\ell-2$.
But each of these can be completed to a hole by adding $P_3$ and at least two further edges (note that, since either $A$ is 
two-connected or $a_0$ exists, there is a path of $A^+$ between $a_1,a_3$ that does not pass through $a_2$, and similarly for 
the other three paths we need), and this hole has length
at least $(\ell-2) +2 +|E(P_3)|>\ell$, a contradiction.
This proves (5).

\bigskip
From (5), and the symmetry between $A,B$, we may assume that every neighbour of $s_n$ in $V(F)$ belongs to $V(A^+)$. Let $I$ be the
set of $i\in \{0\LL k\}$ such that $s_n, a_i$ are adjacent (thus $0\in I$ only if $a_0$ exists). By (4), $I\not\subseteq \{1\}$.
Let $J$ be the set of $i\in \{2\LL k\}$ such that $a_1, a_i$ are nonadjacent.
Let $g',g\in \{1\LL m\}$ be minimum and maximum respectively such that $s_1$ 
is adjacent to $q_{g'}, q_g$.
Since $s_1$ has a neighbour in the interior of $P_1$ it follows that $g\ge 2$. 
\\
\\
(6) {\em No vertex of $V(F)$ has a neighbour in the interior of $S$ except possibly $a_1$, and $a_1$ is adjacent to at least one of $s_1\LL s_{n-1}$, and $g=2$.}
\\
\\
Suppose that 
some internal vertex of $S$ has a neighbour $v\in V(F)$.
From the minimality of $S$, $v$ is adjacent or equal to every neighbour
of $s_1$ in $V(F)$, and so $v\in V(P_1)$; and $v$ is adjacent or equal to every neighbour of $s_n$ in $V(F)$, and
so $v=a_1$ since $I\not\subseteq \{1\}$.  This proves the first statement.

Suppose that $a_1$ has no neighbour in $V(S)\setminus \{s_n\}$, and in particular $g'\ge 2$. 
Let $A'$ be the subgraph induced on $V(A^+)\cup \{s_n\}$, let $B'$ be the subgraph 
induced on
$V(B^+)\cup \{s_1,q_{g'},q_g\LL q_{m}\}$, and  let $P_1'$ be the path $q_1\CC q_{g'}$. Then $(A',B', S, P_1', P_2\LL P_k)$ is a $(k+1)$-bar
semigate, a contradiction. Hence $a_1$ has a neighbour in $V(S)\setminus \{s_n\}$; and by the minimality of $S$ it follows that $a_1$
is adjacent to every neighbour of $s_1$ in $V(F)$, and consequently $g=2$. This proves (6).
\\
\\
(7) {\em $n=2$, and $I\cap \{2\LL k\}\subseteq J$. If $0\in I$ then $I\cap \{2\LL k\}=J$.}
\\
\\
First we show that $I\cap \{2\LL k\}\subseteq J$, and if $I\cap \{2\LL k\}\ne \emptyset$ then $n=2$. 
Both are true if $I\cap \{2\LL k\}= \emptyset$, so we may assume that $i\in I\cap \{2\LL k\}$. There is an induced path 
$R$ between $q_2$ and $a_i$ consisting of $q_2\CC q_m$, 
a one- or two-edge path of $B^+$, and the path $P_i$. This can be completed to a hole via $q_2\DD s_1\CC s_n\DD a_i$, and also via
a path of length at most three between $q_2, a_i$ with interior consisting of $a_1$ and possibly one other vertex of $A^+$.
Consequently both these completions have the same length, and so $n=2$ and $i\in J$ as required. 

It remains to show that if $I\cap \{2\LL k\}= \emptyset$ then $n=2$, and if $0\in I$ then $I\cap \{2\LL k\}=J$.
But if  $I\cap \{2\LL k\}= \emptyset$ then $a_0$ exists and $0\in I$, since $I\not\subseteq \{1\}$. Consequently, for both statements
we may assume that $a_0$ exists and $0\in I$, and under that assumption we must show that $I\cap \{2\LL k\}=J$ and $n=2$. 
Since $a_0$ exists, $A$ is not connected and so there exists some $i\in J$. Choose $i\in J\setminus I$ if possible. 
If $i\in I$, then $I\cap \{2\LL k\}=J$ from the choice of $i$, and $n=2$ from what we already proved, 
since $I\cap \{2\LL k\}\ne \emptyset$. Thus we may assume that $i\notin I$.
There is a path $R$ between $q_2, a_0$ consisting of 
$q_2\CC q_m$,
a one- or two-edge path of $B^+$, the path $P_i$, and the edge $a_ia_0$. Consequently the paths
$q_2\DD s_1\CC s_n\DD a_0$ and $q_2\DD q_1\DD a_0$ have the same length, a contradiction since $n\ge 2$. This proves (7).

\bigskip

Suppose that  $0\notin I$. Then $I\cap J\ne \emptyset$ by (7) since $I\not\subseteq \{1\}$; choose $j\in I\cap J$. Moreover, $A^+$ is a 
connected threshold graph, and so there exists $i\in \{0\LL k\}$ such that $a_i$ is adjacent to all other 
vertices of $A^+$. Since $J\ne \emptyset$,
it follows that $i\ne 1$; and $i\notin J$ since $a_1,a_i$ are adjacent. So $i\in \{0,2\LL k\}\setminus J$, 
and therefore $i\notin I$ by (7). By (6), $a_1$ is adjacent to $s_1$; so one of 
$a_1\DD a_i\DD a_j\DD s_2\DD a_1$, $a_1\DD a_i\DD a_j\DD s_2\DD s_1\DD a_1$ is a hole of length four or five, a contradiction.

This proves that $a_0$ exists and $0\in I$. By (7), $I\setminus \{1\}\subseteq J\cup \{0\}\subseteq I$. Suppose that
there exists $i\in \{2\LL k\}\setminus J$. There is a path $R$ between $q_2,a_i$ consisting of $P_1\setminus \{a_1\}$, $P_i$ 
and a one- or two-edge path of $B^+$. Moreover, $a_i$ is nonadjacent to $s_1,s_2$, and there is an induced path
$q_2\DD s_1\DD s_2\DD a_0\DD a_i$ of length four, and another $q_2\DD a_1\DD a_i$ of length two, and 
for both these paths, their union with $R$ is a hole, a contradiction. This proves that $J=\{2\LL k\}$. If $1\notin I$, there is a 4-hole
or 5-hole with vertex set a subset of $\{a_1,q_2,s_1,s_2,a_0\}$, a contradiction; so $I=\{0\LL k\}$, and the second bullet of 
the theorem holds. 
This proves \ref{framejumps}.~\bbox

\section{Blowing up a frame}

To prove \ref{mainthm} we need to look at 
blow-ups of apexed $\ell$-frames. We will 
handle the $\ell$ odd and $\ell$ even cases together, so the next result lists the properties of apexed $\ell$-frames that we will use;
they hold no matter whether $\ell$ is odd or even.

\begin{thm}\label{frameprop}
Let $k,\ell$ be integers with $k\ge 3$ and $\ell\ge 7$. Let $F$ be an apexed $\ell$-frame with sides $A^+, B^+$ and bars $P_1\LL P_k$, 
and apexes $a_0, b_0$ (if they exist). For $1\le i\le k$ let $P_i$ have ends $a_i\in V(A^+)$ and $b_i\in V(B^+)$.
Let $A=A^+\setminus \{a_0\}$ if $a_0$ exists, and $A=A^+$ otherwise, and define $B$ similarly. Then:
\begin{itemize}
\item $V(A^+), V(B^+)$ are anticomplete, and $A^+, B^+$ are connected;
\item $A$, $B$ are threshold graphs with $|A|=|B|=k$,
such that either $A$ is two-connected and $a_0$ does not exist, or $A$ is not connected and $a_0$ exists; and the same for $B$.
\item 
$V(A)=\{a_1\LL a_k\}$, and $V(B)=\{b_1\LL b_k\}$. 
\item 
For $1\le i\le k$, $P_i$ has length at least $\ell/2-2$ and at most $\ell/2-1$, with no internal vertex in $V(A\cup B)$, and such that
$P_1\LL P_k$ are pairwise vertex-disjoint. 
\item $F$ is $\ell$-holed.
\item there do not exist five vertices $a,b,c,d,x$ of $F$ such that $a\DD b\DD c\DD d$ is an induced path and $x$ is adjacent
to all of $a,b,c,d$,
\end{itemize}\end{thm}

Apexed $\ell$-frames have the following convenient property:
\begin{thm}\label{holeinbed}
If $F$ is an apexed $\ell$-frame and $u,v\in V(F)$, some hole of $F$ contains both $u$ and $v$.
\end{thm}
\Proof
We use the notation of \ref{frameprop}. For $1\le i<j\le k$, there is a hole of $F$ containing $P_i\cup P_j$
(using $a_0$ if $A$ is not connected, and $b_0$ if $B$ is not connected); so we may assume that $u=a_0$ say.
Hence $A$ is not connected; let $I\subseteq \{1\LL k\}$ such that $I\ne \emptyset, \{1\LL k\}$, and 
there are no edges between $\{a_i:i\in I\}$ and $\{a_i:i\in \{1\LL k\}\setminus I\}$. If $v\in V(P_i)$ for some $i\in \{1\LL k\}$,
        choose $j\in \{1\LL k\}$ with exactly one of $i,j$ in $I$;
then there is a hole containing $P_i, P_j$ and $a_0$, as required. So we may assume that also $v=b_0$; and so $B$ is not connected.
        Consequently the set $\{b_i:i\in I\}$ is not complete to $\{b_j:j\in \{1\LL  k\}\setminus I\}$; choose $i\in I$ and
        $j\in \{1\LL k\}\setminus I\}$ with $b_i, b_j$ nonadjacent. Since also $a_i, a_j$ are nonadjacent,
there is a hole contain $a_0,b_0, P_i$ and $P_j$.
This proves \ref{holeinbed}.~\bbox

Let us say two holes $C,C'$ in an $\ell$-holed graph $G$ are {\em close} if $|V(C\cap C')|\ge 4$; and {\em equivalent} if there is a sequence
of holes 
$$C=C_1\LL C_n=C'$$
such that $C_i, C_{i+1}$ are close for $1\le i<n$. We need:
\begin{thm}\label{eqnholes}
Let $\ell\ge 7$, and let $F$ be an apexed $\ell$-frame. Then every two holes in $F$ are equivalent.
\end{thm}
\Proof
We use the notation of \ref{frameprop}. Since $A^+$, $B^+$ are threshold graphs, every hole includes one (and hence two)
of $P_1\LL P_k$. Each of $P_1\LL P_k$ has length at least $\ell/2-2\ge 3/2$, and hence at least two.
\\
\\
(1) {\em Let $C_1,C_2$ be holes of $F$ with some $P_i\subseteq C_1\cap C_2$. Then $C_1,C_2$ are equivalent.}
\\
\\
We may assume that $i=3$. Suppose that $C_1,C_2$ are not equivalent.
In particular, they are not close, so $P_3$ has length two, and $C_1\cap C_2=P_3$.
We may assume that 
$P_1\subseteq C_1$ and $P_2\subseteq C_2$. Let $P_i$ have ends $a_i\in V(A)$ and $b_i\in V(B)$ for $i = 1,2,3$.

Suppose that $a_1,a_3$ are nonadjacent.
If $a_2$ is adjacent to both $a_1,a_3$, then there is a hole close to $C_1$ containing
$P_1,P_3$ and $a_2$, and this hole is also close to $C_2$, a contradiction. Thus $a_2$ is nonadjacent to 
at least one of $a_1,a_3$. Since $A^+$ is a connected threshold graph, it has a vertex $a$ adjacent to all the others; and 
$a\ne a_1,a_2,a_3$. There is a hole close to $C_1$ containing $P_1,P_3$ and $a$, and so we may assume that 
$a\in V(C_1)$. If $a_2,a_3$ are nonadjacent, there is 
a hole containing $P_2,P_3, a$, and this hole is close to both $C_1, C_2$, a contradiction. Thus $a_2,a_3$ are adjacent, and so
$a_1,a_2$ are nonadjacent (because $a_2$ is nonadjacent to
at least one of $a_1,a_3$). There is a hole containing $P_1,P_2$ and $a$, and it is close to $C_1$, so it is not close to 
$C_2$, and therefore $P_2$ has length two; and so $C_2$ has length at most $7$, and therefore $\ell=7$. 
We have shown then that if $a_1,a_3$ are nonadjacent, then $a_2,a_3$ are 
adjacent, and $a_1,a_2$
are nonadjacent, and $P_2$ has length two, and $\ell=7$.

Let us continue to assume that $a_1,a_3$ are nonadjacent. Since $P_2,P_3$ both have length two and $a_2,a_3$ are adjacent, and
$\ell= 7$, it follows that $b_2,b_3$ are nonadjacent. By the argument of the paragraph above, with 
$A,B$ exchanged and $P_1,P_2$ exchanged, it follows that 
$b_1,b_2$
are nonadjacent. But then a hole containing $P_1,P_2$ has length at least eight, a contradiction.

This proves that $a_1,a_3$ are adjacent, and similarly $a_2a_3,b_1b_3,b_2b_3$ are edges. Consequently $P_1,P_2$ both have
length at least three, since $C_1,C_2$ have length $\ell\ge 7$; but there is a hole
containing $P_1,P_2$, and it is close to both $C_1,C_2$, since $P_1,P_2$ have length at least three, a contradiction. This proves (1).

\bigskip
For each pair of distinct $i,j\in \{1\LL k\}$ let $C_{i,j}$ be a hole containing $P_i\cup P_j$. By repeated application of (1)
it follows that all these holes are equivalent; but every hole contains some pair $P_i,P_j$ and so is equivalent to $C_{i,j}$.
This proves \ref{eqnholes}.~\bbox

Let $F$ be an apexed $\ell$-frame, and let $A^+, B^+$ and so on be as in \ref{frameprop}.
Let $G$ be a graph, and for each $t\in V(F)$, let $W_t$ be an ordered non-null clique of $G$.
Suppose that for all distinct $s,t\in V(F)$:
\begin{itemize}
\item $W_s\cap W_t=\emptyset$;
\item if $s,t$ are not $F$-adjacent then $W_s, W_t$ are anticomplete in $G$;
\item if $s,t$ are $F$-adjacent then $G[W_s,W_t]$ is a half-graph that obeys the orderings of $W_s,W_t$, and each of its vertices has positive degree; and
\item if $s,t$ both have degree at least three in $F$,  
and $s,t$ are $F$-adjacent, then $W_s$ is complete to $W_t$ in $G$.
\end{itemize}
We observe:
\begin{thm}\label{adjacency}
If $F,A^+,B^+$ and so on are as above, and 
if $s,t$ are $F$-adjacent, then $W_s$ is complete to $W_t$ unless either $st$ is an edge of some $P_i$, or one of $s,t$ equals $a_0$
and the other is a vertex of $A$ with degree zero in $A$, or the same in $B$.
\end{thm}
We say that the subgraph of $G$ induced on the union of the sets $W_t\;(t\in V(F))$ is a {\em blow-up} of $F$, with {\em bags} 
$W_t\;(t\in V(F))$. We will choose a maximal blow-up of $F$ in $G$ and analyze how the rest of $G$ is attached to this blow-up.

There is a convenient notational simplification, as follows. Given $W_t\;(t\in V(F))$ as above, for each $t\in V(F)$ let $w_t$
be the first term of the ordering of $W_t$. Then the map that sends $t$ to $w_t$ for each $t\in V(F)$ is an isomorphism from $F$
onto an induced subgraph $F'$ of $G$; and it is convenient to index the bags not by vertices of $F$ but by the corresponding vertex of $F'$. Thus, we say that the subgraph of $G$ induced on the union of the sets $W_t\;(t\in V(F))$ is a {\em self-centred blow-up} of $F$,
if it is a blow-up of $F$, and in addition:
\begin{itemize}
\item $F$ is an induced subgraph of $G$; and
\item $t\in W_t$, and $t$ is the first term of the ordering of $W_t$, for each $t\in V(F)$.
\end{itemize}
Whenever we have a blow-up of some apexed $\ell$-frame $F$, it is also a self-centred blow-up of some graph $F'$ isomorphic to $F$, and working 
with a self-centred blow-up is often more convenient than working with a general one. We will often use the following:
\begin{thm}\label{centred}
Let $\ell\ge 7$ be an integer. Let $F$ be an apexed $\ell$-frame, let $G$
be an $\ell$-holed graph, and let $H$ be a self-centred blow-up of $F$ contained in $G$. Let $st\in E(F)$. Then $s$ is $H$-adjacent 
to every vertex in $W_t$ (and $t$ is $H$-adjacent to every vertex in $W_s$).
\end{thm}
\Proof Every vertex $v\in W_t$ has an $H$-neighbour in $W_s$, from the third condition in the definition of a blow-up.
Consequently $v$ is $H$-adjacent to $s$, since the half-graph $G[W_s,W_t]$ obeys the orderings of $W_s,W_t$,
and  $s$ is the first vertex of the ordering of $W_s$.
This proves \ref{centred}.~\bbox

We need:

\begin{thm}\label{blownupeqn}
Let $\ell\ge 7$ be an integer. Let $F$ be an apexed $\ell$-frame, let $G$
be an $\ell$-holed graph, and let $H$ be a blow-up of $F$ that is contained in $G$.
Then every two holes of $H$ are equivalent, and every vertex of $H$ is in a hole of $H$.
\end{thm}
\Proof We use the notation of \ref{frameprop}. By replacing $F$ by an isomorphic graph, we may assume that the blow-up is self-centred.
Let $W_t\;(t\in V(F))$ be the bags of $H$. 
\\
\\
(1) {\em If $C$ is a hole of $H$ then $|W_t\cap V(C)|\le 1$ for each $t\in V(F)$.}
\\
\\
Suppose that $|W_t\cap V(C)|\ge 2$ for some $t\in V(F)$. 
Since $W_t$ is a clique, there are exactly two vertices
of $C$ in $W_t$, say $u,v$; let $u'\DD u\DD v\DD v'$ be a path of $C$. Let $u'\in W_s$ and $v'\in W_r$ say;
then the edges $uu', vv'$ form an induced two-edge matching in $G[W_t, W_r\cup W_s]$. But 
$G[W_t, W_r]$ and $G[W_t,W_s]$ are half-graphs that obey the orderings of $W_r,W_s,W_t$, and so $G[W_t, W_r\cup W_s]$
is a half-graph, and therefore has no induced two-edge matching, a contradiction.
This proves (1).
\\
\\
(2) {\em If $C$ is a hole of $H$ then $\{t:W_t\cap V(C)\ne \emptyset\}$ induces a hole of $F$.}
\\
\\
Let $C$ have vertices $v_1\CC v_{\ell}\DD v_1$ in order, and let $v_i\in W_{t_i}$ for $1\le i\le \ell$. By (1) $t_1\LL t_{\ell}$
are all distinct, and $t_1\DD t_2\CC t_{\ell}\DD t_1$ are the vertices in order of a cycle $D$ of $F$. 
Suppose that it is not induced, and $t_1$ is adjacent in $F$ to $t_i$ say, where $3\le i\le \ell-1$. Hence $G[W_{t_1},W_{t_i}]$ is not complete bipartite, since 
$v_1,v_i$ are nonadjacent in $H$; 
but $t_1,t_i$ both have degree at least three in 
$F$ (because they have degree two in $D$), a contradiction. This proves (2).

\bigskip

Let us call $D$ as in (2) the {\em shadow} of $C$.
\\
\\
(3) {\em If $C$ is a hole of $H$ then it is equivalent to its shadow.}
\\
\\
Let $D$ be the shadow of $C$.
To show that $C,D$ are equivalent we use induction on
the number of vertices of $C$ not in $V(D)$. Let $v$ be such a vertex, with $v\in W_t$ say; then $C$ is close to
the hole $C'$ obtained from $C$ by replacing $v$ by $t$, since $\ell\ge 5$; and $C',D$ are equivalent 
from the inductive hypothesis. This proves (3).

\bigskip

From (3) and \ref{eqnholes}, we deduce that all holes in $H$ are equivalent in $H$. 
This proves the first assertion of the theorem.

For the second, let $v\in V(H)$, with $v\in W_t$ say; we will show that $v$ belongs to a hole of $H$. 
From \ref{holeinbed}, $t$ belongs to some hole $D$ of $F$. Let $C$ be obtained from $D$
by replacing $t$ by $v$; then $C$ is a hole of $H$ containing $v$, by \ref{centred}.
This proves \ref{blownupeqn}.~\bbox

\begin{thm}\label{blownupbed}
Let $k\ge 3$ and $\ell\ge 7$ be integers, and let $F$ be an apexed $\ell$-frame. Let $G$
be an $\ell$-holed graph,  and let $H$ be a maximal blow-up of $F$ contained in $G$.
Let the bags of $H$ be 
$W_t\;(t\in V(F))$.  Then, in the usual notation,  for every vertex 
$x\in V(G)\setminus V(H)$, either:
\begin{itemize}
\item the set of neighbours of $x$ in $V(H)$ is a clique; or
\item $x$ is adjacent to every vertex of $H$; or
\item the set of neighbours of $x$ in $V(H)$ is a subset of $W(A^+)$, and if $a_0$ exists, this set is disjoint from $W_{a_i}$ 
for some $i\in \{1\LL k\}$, or
\item the set of neighbours of $x$ in $V(H)$ is
a subset of $W(B^+)$, and if $b_0$ exists, this set is disjoint from $W_{b_i}$ for some $i\in \{1\LL k\}$.
\end{itemize}
\end{thm}
\Proof
By replacing $F$ by an isomorphic graph, we may assume that the blow-up is self-centred. We use the notation of \ref{frameprop}.
Let $X$ be the set of neighbours of $x$ in $V(H)$, and let $M$ be the set of $t\in V(F)$ with $X\cap W_t\ne \emptyset$.
\\
\\
(1) {\em We may assume that for all $s,t\in M$, if $s,t$ are $F$-nonadjacent then they have a common $F$-neighbour in $M$.
}
\\
\\
Suppose that there exist $s,t\in M$ with distance at least three in $F$. By \ref{holeinbed} there is a hole $D$ of $F$ containing $s,t$.
Choose $x_s\in X\cap W_s$ and $x_t\in X\cap W_t$, and let $C$ be the hole of $G$ induced on 
$\{x_s,x_t\}\cup (V(D)\setminus \{s,t\})$. Thus $C$ has length $\ell$, and $x_s,x_t$ have distance at least three in $C$.
Since $G$ is $\ell$-holed it follows that $V(C)\subseteq X$. But if two holes of $H$ are close and $x$ is $G$-adjacent to all 
vertices of one, then it has at least four neighbours in the other and so is $G$-adjacent to all vertices of the other. By \ref{blownupeqn}
it follows that $x$ is $G$-adjacent to every vertex of $H$ that is in a hole. But every vertex of $H$ is in a hole of $H$, by 
\ref{holeinbed}, and so
the second bullet of the theorem holds. This proves that we may assume that every two vertices in $M$ have distance at most two in $F$. 

Suppose that $s,t\in M$ have distance two in $F$; we claim they have a common neighbour in $M$. Let $r\in V(F)$ be $F$-adjacent to 
them both. Thus $r$ is $H$-adjacent to every vertex in $W_s\cup W_t$, by \ref{centred}.  Choose $x_s\in X\cap W_s$ and 
$x_t\in X\cap W_t$. Since $\{x,x_s,x_t,r\}$ does not induce a 4-hole, it follows
that $r\in X$, and so $r\in M$.
This proves (1).
\\
\\
(2) {\em We may assume that $M$ is not a clique.}
\\
\\
Suppose that $M$ is a clique; and we may assume that there exist $H$-nonadjacent $u,v\in X$, since otherwise
the first bullet of the theorem holds.
Let $u\in W_s$ and $v\in W_t$ say; thus $s,t\in M$, and so they are
distinct (because $W_s$ is a clique of $G$) and $F$-adjacent.
Let $C$ be a hole of $F$ containing $s,t$
(this exists by \ref{holeinbed}). None of its vertices different from $s,t$ belong to $M$ since $M$ is a clique.
By \ref{centred}, the neighbour of $s$ in $C$ different from $t$ is $H$-adjacent to every vertex of $W_s$, and in particular, is $H$-adjacent to $u$; and the same for $t$. Consequently there is an induced path of $H$ between $u,v$, with internal vertices
the vertices of $C\setminus \{s,t\}$; and none of its internal vertices are $G$-adjacent to $x$, since the
corresponding vertices of $C$ are not in $M$.
Adding $x$ to this path gives a hole of length $\ell+1$, a contradiction. This proves (2).

\bigskip
For each $t\in V(F)$, let $N_F(t)$ denote the set of vertices in $F$ that are adjacent to $t$ in $F$, and $N_F[t]=N_F(t)\cup \{t\}$.
\\
\\
(3) {\em There exists $s\in M$ such that $M\subseteq N_F[s]$.}
\\
\\
Choose $s\in M$ $F$-adjacent to as many members of $M$ as possible. We may assume that there exists $s_1\in M$ different 
from and $F$-nonadjacent to $s$. By (1) there exists $t_1\in M$ $F$-adjacent to both $s_1,s$, and therefore $H$-adjacent to 
all vertices in $W_{s}\cup W_{s_1}$, by \ref{centred}. Since $t_1$ is not a better choice than $s$,
there exists $t_2\in M$ $F$-adjacent to $s$ and not to $t_1$. Since $s\DD t_1\DD s_1\DD t_2$ is not a 4-hole of $F$ it follows
that $s_1,t_2$ are $F$-nonadjacent. Hence $s_1\DD t_1\DD s\DD t_2$ is an induced path $P$ of $F$.
By (1), there exists $s_2\in M$ $F$-adjacent to both $s_1,t_2$. Hence $s_2\ne s,t_1$
since each of the latter is nonadjacent to one of $s_1,s_2$, and $s_2$ is adjacent to them both. Since $F$ has no hole of 
length four or five, it follows that $s_2$ is $F$-adjacent to every vertex of the induced path $s_1\DD t_1\DD s\DD t_2$. 
But apexed $\ell$-frames have no four-vertex induced path in the neighbourhood of a vertex, a contradiction.
This proves (3).

\bigskip

Choose $s$ as in (3).
\\
\\
(4) {\em If $p\DD q\DD r$ is an induced path of $F$ with $p,r\in M$, then $q\in M$ and $W_q\subseteq X$. Consequently $W_s\subseteq X$.}
\\
\\
Suppose that $W_q\not\subseteq X$, and choose $y\in W_q\setminus X$. 
Both $x,y$ have a neighbour in the clique $W_p$, so there is an induced path between them with length two or three and with 
interior in $W_p$; and similarly there is such a path with interior in $W_r$. But the union of these paths makes a hole
of length four, five or six, a contradiction. This proves the first assertion. For the second, we may choose nonadjacent
$p,r\in M$, by (2); and so they are both $F$-neighbours of $s$ from the choice of $s$, and hence $W_s\subseteq X$ from
the first assertion. This proves (4).

\bigskip
The vertex $s$ is either an apex of $F$, or an end of one of the paths $P_i$, or an internal vertex of one of these paths, 
and in each case we will show that $x$ can be added to $W_s$, contrary to the maximality of $W(H)$. Before we begin on those cases, let us
see what we need to check to show that adding $x$ to $W_s$ gives a larger blow-up. First, we need:
\begin{itemize}
\item $W_s\cup \{x\}$ is a clique of $G$. 
\end{itemize}
(This is already proved in (4).)
Given this, it follows that for each $F$-neighbour $q$ of $s$, the graph $G[W_s\cup\{x\},W_q]$ is a half-graph. Second, we need 
to check that every vertex of $G[W_s\cup\{x\},W_q]$ has positive degree, and so we need:
\begin{itemize}
\item For each $F$-neighbour $q$ of $s$, $x$ has a $G$-neighbour in $W_q$.
\end{itemize}
But we also need to
make $W_s\cup \{x\}$ an ordered clique. To do so, we need:
\begin{itemize}
\item For every two $F$-neighbours $q,q'$ of $s$, the graphs $G[W_s\cup \{x\},W_q]$ and $G[W_s\cup\{x\},W_{q'}]$ are compatible.
\end{itemize}
Given this, it follows that, if $Q$ is the set of all $F$-neighbours of $s$, then $G[W_s\cup \{x\},\bigcup_{q\in Q}W_q]$
is a half-graph, and hence there is an ordering of $W_s\cup \{x\}$ such that every vertex in $\bigcup_{q\in Q}W_q$
is adjacent to an initial segment of this ordering. This will be the new ordering of $W_s\cup \{x\}$.
But there is another thing to check: we need that for each $q\in Q$, $x$ is adjacent to an initial segment of the ordering of $Q_q$,
and this might not be true; we might need to change the ordering of $W_q$. So finally, we need:
\begin{itemize}
\item For each $q\in Q$ and each $F$-neighbour $r$ of $q$ different from $s$, the graphs $G[W_q,W_s\cup \{x\}]$ and $G[W_q,W_r]$
are compatible.
\end{itemize}
If this is true, then we can choose a new ordering of each $W_q$, in the same way that we did for $W_s\cup \{x\}$. Depending on
the location of $s$, there might be additional conditions that need to be satisfied before we can add $x$ to $W_s$, but we handle them
case-by-case. Note that the enlarged blow-up might not be self-centred.

Now let us turn to the cases.
\\
\\
(5) {\em We may assume that $s$ is not an apex of $F$.}
\\
\\
Suppose that $s=a_0$ say, and so $X\subseteq W(A^+)$ by (3). We must show that $X$ is disjoint from one of $W_{a_1}\LL W_{a_k}$.
Suppose not. By (4)
it follows that $W_{a_0}\subseteq N_H(x)$. We claim that adding $x$ to $W_{a_0}$ gives a larger blow-up of $F$,
a contradiction. To show this, as we saw above, we need to check that
\begin{itemize}
	\item $x$ has a neighbour in $W_{a_i}$ for $1\le i\le k$;
        \item for $1\le i<j\le k$, the graphs $G[W_{a_0}\cup \{x\},W_{a_i}]$ and $G[W_{a_0}\cup \{x\},W_{a_j}]$ are compatible;
        \item for $1\le i \le k$, and every $F$-neighbour $q$ of $a_i$ different from $a_0$, the graphs
                $G[W_{a_i}, W_{a_0}\cup \{x\}]$ and $G[W_{a_i}, W_q]$ are compatible;
	\item for $1\le i\le k$, if $x$ is not complete to $W_{a_i}$ then $a_i$ has no neighbour in $A^+\setminus \{a_0\}$.
\end{itemize}
(We need the final condition because of the special rule for $a_0$ in the definition of a blow-up.)
The first is true by hypothesis.
For the second bullet, suppose it is false; then there is an induced path $a\DD b\DD c\DD d$ where $b,c\in W_{a_0}\cup \{x\}$, and $a\in W_{a_i}$,
and $d\in W_{a_j}$ (and so one of $b,c$ equals $x$). Hence there are induced paths $a\DD a_0\DD d$ and $a\DD b\DD c\DD d$
between $a,d$ of different lengths. Let $c_i$ be the neighbour of $a_i$ in $P_i$; then $c_i$ is complete to $W_{a_i}$, and in 
particular, $c_i$ is $G$-adjacent to $a\in W_{a_i}$. Consequently there is an induced path $P_i'$ of $G$ between $a,b_i$, with 
the same length as $P_i$. Define $P_j'$ similarly, between $d,b_j$; and let $Q$ be an induced path of $B^+$ between
$b_i,b_j$. Then $P_i'\cup Q\cup P_j'$ is an induced path between $a,d$. Moreover, its union with either of $a\DD a_0\DD d$ and 
$a\DD b\DD c\DD d$ is a hole, and these two holes have different lengths, a contradiction.
This proves the second bullet.

We prove the last two statements together. Let $1\le i\le k$. If $x$ is complete to $W_{a_i}$ then both statements are true,
so we may assume that $x$ is not complete to $W_{a_i}$. We assume that $i=1$ without loss of generality. 
Suppose that $a_1a_2$ is an edge of $F$. Since $a_0$ exists, $A$ is disconnected; let $(K,L)$ be a partition of
$V(A)$ into two nonempty sets such that there are no edges of $F$ between $K,L$, where $a_1\in K$. Consequently $a_2\in K$; choose $a_3\in L$,
and for $i = 1,2,3$ let $v_i\in W_{a_i}$ be adjacent to $x$. Choose $u\in W_{a_1}$ nonadjacent to $x$.
Let $p$ be the neighbour of $a_1$
in $P_1$. Thus $p$ is $H$-adjacent to $u$ by \ref{centred}, and so there are induced paths $p\DD v_1\DD x\DD v_3$ and $p\DD u\DD v_2\DD x\DD v_3$ of different lengths between $p,v_3$,
yielding a contradiction as usual. This proves the final bullet, that $a_1$ has degree two in $F$. Let $q$ be an 
$F$-neighbour $a_1$ different from $a_0$. Therefore $q$
is the neighbour of $a_1$ in $P_1$. If $G[W_{a_1}, W_q]$ and $G[W_{a_1}, W_{a_0}\cup \{x\}]$ are not compatible, 
there is an induced path $x\DD b\DD c\DD d$ where $b,c\in W_{a_1}$ and $d\in W_q$.
Choose $j\in \{2\LL k\}$, and choose $v_j\in W_{a_j}$ adjacent to $x$.
Then there are induced paths $v_j\DD a_0\DD a_1\DD d$ and $v_j\DD x\DD b\DD c\DD d$  between $v_j,d$ of different lengths,
yielding a contradiction as usual. This proves the third bullet,
and so proves (5).

\bigskip
By (5), we may assume that $s\in V(P_1)$ say.
\\
\\
(6) {\em $s\in \{a_1,b_1\}$.}
\\
\\
Suppose that $s$ is an internal vertex of $P_1$, and let $r,t$ be its neighbours in $P_1$. Thus $X\subseteq W_r\cup W_s\cup W_t$.
We would like to add $x$ to $W_s$ and obtain a larger blow-up, and as before we need to check that:
\begin{itemize}
        \item $x$ has neighbours in $W_{r}, W_t$;
        \item the graphs $G[W_{s}\cup \{x\},W_{r}]$ and $G[W_{s}\cup \{x\},W_{s}]$ are compatible;
        \item for each $F$-neighbour $q$ of $r$ different from $s$, the graphs
                $G[W_{r}, W_{s}\cup \{x\}]$ and $G[W_{r}, W_q]$ are compatible (and the same for neighbours of $t$).
\end{itemize}
The first holds by (2).
If the second statement is false, 
there is, without loss of generality, a four-vertex path $a\DD b\DD x\DD c$ where $a\in W_r, b\in W_s$ and $c\in W_t$. This is impossible
because there is a hole $D$ of $F$ containing $P_1$ and hence containing $r,s,t$; and the subgraph induced on
$$\{a,b,c,x\}\cup (V(D)\setminus \{r,s,t\})$$
is a hole of $G$ of length $\ell+1$. 

For the final condition,
suppose that $G[W_r, W_s\cup \{x\}]$ and 
$G[W_r,W_q]$ are not compatible. 
Then there is an induced path $x\DD a\DD b\DD c$ where $a,b\in W_r$ and $c\in W_q$; but then, let $D$ be a hole
of $F$ containing $q,r,s$, and then the subgraph induced on
$$\{a,b,c,x\}\cup (V(D)\setminus \{q,r,s\})$$
is a hole of $G$ of length $\ell+1$. 

Consequently we can add $x$ to $W_s$ and obtain a blow-up of $F$ larger than $H$, a contradiction. 
This proves (6).

\bigskip
From (6) and the symmetry, we may assume that $s=a_1$. Let $p$ be the neighbour of $a_1$ in $P_1$. 
\\
\\
(7) {\em We may assume that $X\cap W_p\ne \emptyset$.}
\\
\\
Suppose that $X\cap W_p= \emptyset$, and so $X\subseteq W(A^+)$. If $a_0$ does not exist, then the third outcome of the theorem holds.
If $a_0$ exists, then $A$ is not connected,
and so there exists $i\in \{2\LL k\}$ with $a_i$ nonadjacent to $a_1$. Consequently $a_i\notin M$ by the choice of $s$, and so $X\cap W_{a_i}=\emptyset$,
and again the third outcome holds. This proves (7).

\bigskip

Choose $u\in X\cap W_p$.
\\
\\
(8) {\em $a_0$ does not exist.}
\\
\\
Suppose that $a_0$ exists. 
Suppose first that $a_1$ is nonadjacent in $F$ to all of $a_2\LL a_k$. Thus $M=\{p,a_1,a_0\}$, since $M\subseteq N_F[s]$ and $M$ is not
a clique. Let $q$ be the neighbour of $p$ in $P_1$ different from $a_1$. We claim that adding $x$ to $W_{a_1}$  
gives a larger
blow-up of $F$, and to show this, as before, we must check that 
\begin{itemize}
        \item $p,a_0\in M$;
        \item the graphs $G[W_{a_1}\cup \{x\},W_{a_0}]$ and $G[W_{a_1}\cup \{x\},W_{p}]$ are compatible;
	\item the graphs $G[W_p,W_{a_1}\cup \{x\}],G[W_p,W_q]$ are compatible;
        \item for $2\le i\le k$, the graphs $G[W_{a_0}, W_{a_1}\cup \{x\}]$ and $G[W_{a_0}, W_{a_i}]$ are compatible.
\end{itemize}
We have already seen the first statement. The second holds since otherwise there is an
induced path
$a\DD b\DD c\DD d$ with $a\in W_{p}, b,c\in W_{a_1}\cup \{x\}$ and $d\in W_{a_0}$, which could be extended to a hole of length $\ell+1$ in $G$. 
The third holds since otherwise there is an induced path 
$x\DD b\DD c\DD d$ with $b,c\in W_p$ and $d\in W_q$, which could be extended to a hole of length $\ell+1$ of $G$.
The fourth holds since otherwise there is an induced path      
$x\DD b\DD c\DD d$ with $b,c\in W_{a_0}$ and $d\in W_{a_i}$, which could be extended to a hole of length $\ell+1$ of $G$.
This contradicts the maximality of $H$.

Consequently $a_1$ is adjacent to at least one of $a_2\LL a_k$, and so has degree at least three in $F$; and so $W_{a_1}$
is complete to $W_{a_i}$ for each $i\in \{0,2,3\LL k\}$ such that $a_1,a_i$ are adjacent in $F$.

Suppose that $a_0\notin M$; then we may assume that $a_2\in M$, and so $a_1,a_2$ are adjacent. Choose $v_2\in W_{a_2}\cap X$.
Since $a_0$ exists, $A$ is disconnected, so we may assume that $a_3$ is nonadjacent to $a_1,a_2$, and hence
$a_3\notin M$. Then there are induced paths $u\DD a_1\DD a_0\DD a_3$ and $u\DD x\DD v_2\DD a_0\DD a_3$,
and they can both be extended to holes by taking their union with the same path, but their lengths are different, a contradiction.

So $a_0\in M$. Choose $v_0\in W_{a_0}\cap X$. We claim that $W_{a_i}\subseteq X$ for each $i\in \{2\LL k\}$
such that $a_1,a_i$ are adjacent in $F$. Suppose not; then we may assume that $a_1,a_2$ are adjacent, and $v_{2}\in W_{a_2}$
is not adjacent to $x$. There are paths $u\DD a_1\DD v_2$ and $u\DD x\DD v_0\DD v_2$, giving a contradiction as usual.
Thus $W_{a_i}\subseteq X$ for each $i\in \{2\LL k\}$
such that $a_1,a_i$ are adjacent in $F$.

We claim that $W_{a_0}\subseteq X$. Suppose not, and let $w_0\in W_{a_0}\setminus X$.
Since $a_1$ is adjacent to at least one of $a_2\LL a_k$, we may assume that $a_1,a_2$ are adjacent, and so $W_{a_2}\subseteq X$. 
As before, since $a_0$ exists, $A$ is disconnected, so we may assume that $a_3$ is nonadjacent to $a_1,a_2$, and hence
$a_3\notin M$. There are induced paths $u\DD a_1\DD w_0\DD a_3$ and $u\DD x\DD a_2\DD w_0\DD a_3$, giving a contradiction as usual 
(note that $w_0, a_3$ are adjacent since the blow-up $H$ is self-centred).

This proves that $W_t\subseteq X$ for each $t\in N_F[a_1]\setminus \{p\}$. We claim that adding $x$ to $W_{a_1}$ gives a larger
blow-up of $F$, and to show this we must just check that 
$G[W_p,W_{a_1}\cup \{x\}]$ is compatible with $G[W_p,W_q]$, since all the other conditions are clear.
This holds for the same reason as before, namely
that otherwise there is an induced path 
$x\DD b\DD c\DD d$ with $b,c\in W_p$ and $d\in W_q$, which could be extended to a hole of length $\ell+1$ in $G$, a contradiction.
This proves (8).

\bigskip
Define $N=N_F[a_1]\setminus \{a_1,p\}$. 
\\
\\
(9) {\em For all distinct adjacent $t,t'\in N$, either $W_t\cup W_{t'}\subseteq X$, or $(W_t\cup W_{t'})\cap X=\emptyset$.}
\\
\\
Suppose not; then we may assume that $\{t,t'\}=\{a_2,a_3\}$, and $v_2\in W_{a_2}\setminus X$ and $v_3\in W_{a_3}\cap X$.
But then there are two induced paths $u\DD a_1\DD v_2$ and $u\DD x\DD v_3\DD v_2$, a contradiction as before. This proves (9).

\bigskip

Since $A$ is a two-connected threshold graph (because $a_0$ does not exist by (8)), there are at least two vertices of $A$ that are 
adjacent to all other vertices of $A$, and so one of them belongs to $N$. Consequently $F[N]$
is connected. Since $X\not\subseteq W_s\cup W_p$ by (2), it follows that $M\cap N\ne\emptyset$; and by (9), since
$F[N]$ is connected, it follows that $W_t\subseteq X$ for all $t\in N$. But then we can add $x$ to $W_{a_1}$ as before, 
a contradiction to the maximality of $H$. This proves \ref{blownupbed}.~\bbox

\section{Putting pieces together}\label{sec:parts}

Let $J$ be a graph, and let $A$ be a threshold graph contained in $J$ with $|A|\ge 3$.
For each $t\in V(A)$ let $W_t$ be a non-null ordered clique of $J$ with $t\in W_t$, called a {\em bag},
all pairwise disjoint, such that if $st$ is an edge of $A$ then $W_s$ is complete to $W_t$ in $J$
and otherwise they are anticomplete. Suppose that
\begin{itemize}
\item $J$ is connected, and has no hole of length four or five;
\item each vertex in $V(J)\setminus W(A)$ has two nonadjacent neighbours in $V(A)$;
\item for every induced path $P$ of $J$ with length at least three and with both ends in $W(A)$, some internal vertex of $P$
belongs to the same bag as one of the ends of $P$;
\item for each $t\in V(A)$, $t$ is the first term of the ordering of $W_t$; and for each $v\in V(J)\setminus W_t$,
$v$ is adjacent to an initial segment of the ordering of $W_t$.
\end{itemize}
In this case we call $(J, A, (W_t:t\in V(A)))$ a {\em border}. (Later, we will discuss how to construct borders.)

We were previously working with blow-ups of apexed $\ell$-frames, but now we want blow-ups of $\ell$-frames without the apexes.
With $\ell$ odd or even, let $F$ be an $\ell$-frame, in the usual notation.
For each $t\in V(F)$,
let $W_t$ be an ordered clique where $t$ is the first term of the ordering of $W_t$, all pairwise disjoint, and we will define a graph $H$ with vertex set the union of these cliques.
For every edge $uv$ of $A\cup B$
we make $W_u$ complete to $W_v$ in $H$. For every other edge $uv$ of $F$ let $H[W_u,W_v]$ be a half-graph that obeys the orderings
of $W_u,W_v$.
We call such a graph $H$ a {\em blow-up of an $\ell$-frame}. (Note that here we are including the ``self-centred'' condition in the definition of a blow-up.)

Now take a graph $H$ that is a blow-up of an $\ell$-frame.
Choose $J, K$ such that $(J, A, (W_t:t\in V(A)))$ is a border, with $V(J)\cap V(H)=W(A)$, and
$(K,B,(W_t:t\in V(B)))$ is a border,
with $V(J)\cap V(H)=W(B)$, and $V(K)\cap V(J)=\emptyset$ and  $V(J)$ is anticomplete to $V(K)$. We call the graph $H\cup J\cup K$
a {\em bordered blow-up of an $\ell$-frame}, and say that this graph is the {\em composition} of $H,J,K$.

By combining several of the previous theorems, we obtain the following:

\begin{thm}\label{partition}
Let $\ell\ge 7$, and let $G$ be $\ell$-holed, with no clique cutset or universal vertex. 
Then either $G$ is a blow-up of an $\ell$-cycle, or a bordered blow-up of an $\ell$-frame.
\end{thm}
\Proof
By \ref{cycleblowup}, we may assume that $G$ contains a theta, pyramid or prism, and so $G$ contains a $3$-bar gate.
Choose $k$ maximum such that $G$ contains a $k$-bar gate. By \ref{evenframe} and \ref{oddframe}, $G$ contains an apexed $\ell$-frame
$F$ with sides $A^+, B^+$ and bars $P_1\LL P_k$ say. 
Choose a maximal blow-up $H$ of $F$ contained in $G$. By replacing $F$ by an isomorphic graph, we may assume the blow-up is self-centred.
If some vertex of $A^+$ is not in $V(P_1\cupcup P_k)$, it is unique;
call it $a_0$ (and otherwise $a_0$ is undefined).
Define $b_0$ similarly.  As usual, let $A=A^+\setminus \{a_0\}$ if $a_0$ exists, and $A=A^+$ otherwise, and define $B$ similarly.
Let $H$ have bags $W_t\;(t\in V(F))$.
Let $Z$ be the set of all vertices in $V(G)\setminus V(H)$ that are adjacent to every vertex of $H$. 
Since $H$ has two nonadjacent vertices and $G$ has no 4-hole, it follows that $Z$ is a clique. 
\\
\\
(1) {\em If $x\in V(G)\setminus (V(H)\cup Z)$, and $N_H(x)$ denotes the set of vertices in $V(H)$ adjacent to $x$, then either
\begin{itemize}
\item $N_H(x)$ is a clique, or 
\item $N_H(x)$ is a 
subset of $W(A^+)$, and if $a_0$ exists then $N_H(x)$ is disjoint from $W_{a_i}$ for some $i\in \{1\LL k\}$, or 
\item $N_H(x)$ is a
subset of $W(B^+)$, and if $b_0$ exists then  $N_H(x)$ is disjoint from $W_{b_i}$ for some $i\in \{1\LL k\}$.
\end{itemize}}
\noindent This is immediate from \ref{blownupbed}, since $x$ is not adjacent to every vertex of $H$ by definition of $Z$.
\\
\\
(2) {\em If $C$ is a component of $G\setminus (V(H)\cup Z)$, and $N_H(C)$ denotes the set of vertices in $V(H)$ with a neighbour 
in $V(C)$, then $N_H(C)$ is a subset of one of $W(A^+), W(B^+)$.}
\\
\\
Suppose not. We claim that there exist $u,v\in N_H(C)$, nonadjacent, and not both in $W(A^+)$ and not both in 
$W(B^+)$. The proof is like that of step (1) of the proof of \ref{framejumps}. Since $G$ has no clique cutset, the set of vertices
of $G\setminus V(C)$ with a neighbour in $V(C)$ is not a clique; let $u,v$ be nonadjacent vertices in this set. Thus $u,v\in V(H)\cup Z$,
and $u,v\notin Z$ since every vertex of $Z$ is adjacent to all other vertices of $V(H)\cup Z$. Thus $u,v\in N_H(C)$.
We claim that we can choose $u,v\in N_H(C)$, nonadjacent, and not both in $W(A^+)$ and not both in
$W(B^+)$.
We have already seen that we can choose them in $N_H(C)$.
If $u$ or $v$ belongs to $V(H)\setminus (W(A^+)\cup W(B^+))$, or if $u\in W(A^+)$ and $v\in W(B^+)$, then $u,v$ satisfy the requirement;
so we assume that $u,v\in W(A^+)$. 
Choose $w\in N_H(C)\setminus W(A^+)$. It
follows that $w$ is nonadjacent to one of $u,v$, say $v$; but then $w,v$ satisfy the requirement. This proves that we may choose 
$u,v\in N_H(C)$, nonadjacent, and not both in $W(A^+)$ and not both in                        
$W(B^+)$. Let $u\in W_r$ and $v\in W_s$, where $r,s\in V(F)$. It follows that not both $r,s$ are in $V(A^+)$ and not both are in $V(B^+)$.

We claim that we can choose $u,v,r,s$ as just described, with $r,s$ not $F$-adjacent. 
Because suppose that $r,s$ are $F$-adjacent. Since not both $r,s$ are in $V(A^+)$ and not both are in $V(B^+)$, 
it follows that $rs$ is an edge of one of $P_1\LL P_k$, say $P_1$. Since $P_1$ has length at least two, we may assume
that $s$ is an internal vertex of $P_1$. Let $r\DD s\DD t$ be a subpath of $P_1$. If $N_H(C)$ is not a subset of $W_r\cup W_s\cup W_t$,
then we may choose $q\in V(F)$, $F$-nonadjacent to $s$, such that $N_H(C)\cap W_q\ne \emptyset$ and the claim holds. Thus we assume that
$N_H(C)\subseteq W_r\cup W_s\cup W_t$.  If $N_H(C)\cap W_t\ne \emptyset$ then again the claim holds, with $r,t,u$ and a vertex of $N_H(C)\cap W_t$;
so we assume that $N_H(C)\subseteq W_r\cup W_s$. But there is a hole $D$ of $F$ containing $P_1$; let $q\DD r \DD s \DD t$ be a path of $D$. By deleting
$r,s$ from $D$ and adding $u,v$ and the edges $uq,vt$, we obtain an induced path between $u,v$ of length $\ell-1$; and the union of this path
with an induced path between $u,v$ with interior in $C$ makes a hole of length more than $\ell$, a contradiction.

This proves that we can choose $u,v,r,s$ with $r,s$ not $F$-adjacent. Because $r,s$ are not $F$-adjacent, it follows that
(in the usual notation), the subgraph induced on $(\{u,v\}\cup V(F))\setminus \{r,s\}$ is an apexed $\ell$-frame $F'$ with sides $A'^+,B'^+$ 
and apexes $a_0', b_0'$ say (if they exist). Thus $a_0'$ exists if and only if $a_0$ exists, and either $a_0'=a_0$, or one of $r,s=a_0$, and the same for $b_0$; and either $A'^+ = A^+$, or one of $r,s\in A^+$ (say $r$) and $A'^+=(A^+\cup \{u\})\setminus r$, and the same for $B'^+$.
By \ref{framejumps}, applied to the subgraph of $G$ induced on $V(C)\cup V(F')$,
either
\begin{itemize}
\item there exist $x\in V(C)$ and $u',v'\in V(F')$, such that $u',v'$ are nonadjacent, both are adjacent to $x$, and
not both are in $V(A'^+)$, and not both are in $V(B'^+)$; or
\item there exists $x\in V(C)$ such that either $a_0$ exists and the set of neighbours of $x$ in $V(F')$ is $V(A'^+)$, or
$b_0$ exists and the set of neighbours of $x$ in $V(F')$ is $V(B'^+)$; or
\item 
either $N_{F'}(C)\subseteq V(A'^+)$, or $N_{F'}(C)\subseteq V(B'^+)$, or $N_{F'}(C)$ is a clique, where $N_{F'}(C)$ denotes the set of 
vertices in $V(F')$ with a neighbour
in $V(C)$.
\end{itemize}
The first is false, because of (1). The third is false, because $u,v\in N_{F'}(C)$. Thus the second holds; and we may assume that
$a_0'$ (and hence $a_0$) exists, and there exists $x\in V(C)\setminus V(F')$ such that the set of neighbours of $x$ in $V(F')$ is $V(A'^+)$.
But then $x$ has a neighbour in each of $W_{a_1}\LL W_{a_k}$, contrary to (1).
This proves (2).

\bigskip

For every component $C$ of $G\setminus (V(H)\cup Z)$, $N_H(C)\ne \emptyset$, since $G$ has no clique cutset.
Let $X'$ be the union of the vertex sets of all components $C$ of $G\setminus (V(H)\cup Z)$ such that $N_H(C)\subseteq W(A^+)$, and 
define $Y'$ similarly for $W(B^+)$. Let $X=X'\cup W_{a_0}$ if $a_0$ exists, and $X=X'$ otherwise; and define $Y$ similarly.
By (2), $X',Y',Z$ is a partition of $V(G)\setminus V(H)$; and 
$X,Y$ are anticomplete, also by (2).

We will show that $Z=\emptyset$, and every vertex in $X$ has two nonadjacent neighbours in $V(A)$, and the same for $Y$.
The first step is to show that every vertex in $X$ has two nonadjacent neighbours in $W(A)$, so let 
let $V$ be the set of all $x\in X$ that have two nonadjacent neighbours in $W(A)$, and $U=X\setminus V$. 
If $a_0$ exists then $A$ is disconnected, and so every vertex in $W_{a_0}$ has two nonadjacent neighbours in $W(A)$; and so $W_{a_0}\subseteq V$.
\\
\\
(3) {\em $V$ is complete to $Z$.}
\\
\\
Let $x\in V$, and let $z\in Z$. Since $x\in V$, there exist nonadjacent $u,v\in W(A)$ both adjacent to $x$. But they are also both 
adjacent to $z$, and since $G$ has no 4-hole, it follows that $x,z$ are adjacent. This proves (3).
\\
\\
(4) {\em If $Q$ is an induced path of $G$ with both ends in $W(A)$ and with interior in $X$, then $Q$ has length one or two.}
\\
\\
Let $Q$ have ends $u,v$. We may assume that $u,v$ are nonadjacent, since otherwise $Q$ has length one. Let $u\in W_{a_1}$; 
then $v\notin W_{a_1}$ since $W_{a_1}$ is a clique, so we may assume that $v\in W_{a_2}$.
Since $u,v$ are $G$-nonadjacent, and hence $W_{a_1}$ is not complete to $W_{a_2}$, it follows from \ref{adjacency} that 
$a_1,a_2$ are $F$-nonadjacent. 
Now there is a vertex in $A^+$ adjacent to all other vertices of $A^+$, say $a_i$; and $a_i\ne a_1,a_2$. 
There is a path $R$ of $F$ between $a_1,a_2$ consisting of $P_1,P_2$ and a one- or two-edge path of $B^+$. Since
the union of $R$ with the path $a_1\DD a_i\DD a_2$ is a hole of $F$, it follows that $R$ has length $\ell-2$. There is 
a path of $H$ between $u,v$ with interior $V(R)\setminus \{a_1,a_2\}$; and it also has length
$\ell-2$. But its union with $Q$ is a hole, so $Q$ has length two. This proves (4).

\bigskip
If we can prove that $U=\emptyset$, then (3) implies that $Z$ is complete to $X$, and similarly $Z$ is complete to $Y$; but then 
the vertices in $Z$ are all universal vertices, and so $Z=\emptyset$. Thus, we next need to show that $U=\emptyset$.

Suppose that $U\ne \emptyset$, and let $C$ be a component of $G[U]$. Since $G$ has no clique cutset, there exist nonadjacent $u,v\in V(G)\setminus U$, 
both with neighbours in $V(C)$. Consequently $u,v\in W(A)\cup V\cup Z$. By (3), every vertex in $Z$ is adjacent to every other vertex 
in $W(A)\cup V\cup Z$, and so $u,v\notin Z$.
Thus $u,v\in W(A)\cup V$. 

Either $u,v\in W(A)$, or one is in $W(A)$ and the other in $V$, or they are both in $V$.
We handle these cases separately. Choose an induced path $P$ with ends $u,v$ and interior in $V(C)$.
\\
\\
(5) {\em Not both $u,v$ belong to $W(A)$.}
\\
\\
This is immediate from (4), since if $P$ has length two then its middle vertex is in $V$ by definition of $V$, a contradiction. This proves (5).
\\
\\
(6) {\em Neither of $u,v$ is in $W(A)$.}
\\
\\
Suppose that $u\in W(A)$ say; and so $v\in V$, and so $v$ has nonadjacent neighbours $x,y\in W(A)$. Since $u,v$ are nonadjacent
and $G$ has no 4-hole, it follows that $u$ is nonadjacent to one of $x,y$, say $x$. The path obtained by adding
$x$ and the edge $xv$ to $P$ might not be induced, but its vertex set includes that of an induced path between $u,x$,
which therefore has length two by (4); and consequently the neighbour of $u$ in $P$ is adjacent to $x$, a contradiction, since this neighbour is in $U$. This proves (6).
\\
\\
(7) {\em Not both $u,v$ are in $V$.}
\\
\\
Suppose that $u,v\in V$, and so they both have two nonadjacent neighbours in $W(A)$; say $x,y$ for $u$, and
$x',y'$ for $v$. Choose $x,y,x',y'$ not all distinct if possible.
Now $v$ is not adjacent to both $x,y$, since $G$ has no 4-hole; so we may assume that $v$ is nonadjacent to $x$,
and similarly we may assume that $u$ is nonadjacent to $x'$. So either $x,y,x',y'$ are different, or $y=y'$. Suppose first
that $y=y'$. Since $G$ has no 5-hole it follows that $x,x'$ are nonadjacent; and so by (4) there is an induced path of length two
between $x,x'$
with middle vertex in $V(P)$, a contradiction. (This middle vertex cannot be $u$ or $v$, since $u$ is not adjacent to $x'$ and $v$
is not adjacent to $x$; 
and cannot belong to the interior of $P$ since all those vertices are in $U$.)

This proves that $x,y,x', y'$ are all different; and it is not possible to choose them not all different.
Consequently, if $z\in \{x,y\}$ and $z'\in \{x',y'\}$ are nonadjacent, then neither of $u,v$ is adjacent to both $z,z'$. 
If both of $x,y$ are adjacent to both of $x',y'$ then $G$ has a 4-hole; so we may
assume that some $z\in \{x,y\}$ is nonadjacent to some $z'\in \{x',y'\}$; and so neither of $u,v$ is adjacent to both $z,z'$.
But then by (4), there is an induced path between $z,z'$ of length two
with middle vertex in $V(P)$, and since this middle vertex is not in $U$, and not $u$ or $v$, this is a contradiction.
This proves (7).

\bigskip
From (5), (6), (7) it follows that $U=\emptyset$ and so $Z=\emptyset$.
It remains to show that 
$(G[X\cup W(A)],A, (W_t:t\in V(A)))$ is a border (and the same for $Y,B$).
First we show:
\\
\\
(8) {\em For every induced path $P$ of $G[X\cup W(A^+)]$ of length at least three with both ends in $W(A)$, either the first 
two vertices or the last two vertices of $P$ belong to the same bag $W_t$ for some $t\in V(A)$.}
\\
\\
Let $P$ have ends $u\in W_{a_1}$ and $v\in W_{a_2}$ say. By \ref{adjacency}, $a_1,a_2$ are not $A$-adjacent.
Since $H\setminus W(A^+)$ is connected, and $u,v$ both have neighbours in it, there is an induced path $R$
between $u,v$ with all its internal vertices in $V(H)\setminus W(A^+)$. There is a vertex $a_j$ say of $A^+$ adjacent to all other vertices of $A^+$; and so $j\ne 1,2$. Hence there is an induced path
$u\DD a_j\DD v$. The union of this path and $R$ is a hole, so $R$ has length $\ell-2$. 
It follows that $P\cup R$ is not a hole; so some internal vertex of $P$
belongs to one of $W_{a_1}, W_{a_2}$. Since $W_{a_1}, W_{a_2}$ are cliques and $P$ is induced, this proves (8).
\\
\\
(9)  {\em There is no induced path $v_1\DD v_2\DD v_3\DD v_4$ of $G$ with $v_2,v_3\in W_t$ for some $t\in V(A)$.}
\\
\\
Suppose that there is such a path. Since $W_t$ is a clique, it follows that $v_1,v_4\notin W_t$, and since they are both mixed on 
$W_t$, it follows that $v_1,v_4\notin W(A)$.
Suppose that $v_1\notin X$; so we may assume that $t=a_1$ and $v_1\in W_q$
where $a_1\DD q\DD r$ is a subpath of $P_1$. Since 
$v_1,v_4$ are not adjacent, it follows that $v_4\notin W_q$, and so $v_4\in X$. Consequently $v_4$ has two nonadjacent neighbours
$x,y\in W(A)$ (since $U=\emptyset$). Since $\{x,y,v_2,v_4\}$ does not induce a 4-hole, one of $x,y$ is nonadjacent to $v_2$, say $x$; and hence $x$
is also nonadjacent to $v_3$, since $v_2,v_3\in W_{a_1}$ and $x\in W(A)$. Let $x\in W_{a_2}$ say.
There is an induced path $Q$ of $F$ between $r,a_2$ consisting of $P_1\setminus \{a_1,q\},P_2$ and an induced path of $B^+$ with length one or two;
and so there is an induced path $R$ of $G$ between $r,x$ with the same interior as $Q$, by \ref{centred}. This can be completed to a hole
via $r\DD q\DD v_3\DD v_4\DD x$ or via $r\DD v_1\DD v_2\DD v_3\DD v_4\DD x$, a contradiction.

It follows that $v_1,v_4\in X$. Consequently $v_1$ has two nonadjacent neighbours $x_1,y_1$ in $W(A)$ (since $U=\emptyset$), 
and similarly $v_4$ has two nonadjacent 
neighbours $x_4,y_4\in W(A)$. 
Not both $x_1,y_1$ have a neighbour in $\{v_3,v_4\}$, since otherwise $G$ contains a 4-hole or 5-hole; so we assume that 
$x_1$ is nonadjacent to both $v_3,v_4$, and in particular $x_1\ne v_2$. Similarly we may assume that $x_4$ is nonadjacent to
$v_1,v_2$. Thus $x_1,v_1,v_2,v_3,v_4,x_4$ are all distinct, and $x_1v_1,v_1v_2,v_2v_3,v_3v_4,v_4x_4$ are edges, and all other pairs of these six 
vertices are nonadjacent except possible $x_1v_2, x_4v_3$ and $x_1x_4$. If $x_1,x_4$ are adjacent then $G$ contains a 
4-,5- or 6-hole, a contradiction. Consequently there is an induced path with ends $x_1,x_4$ of length three, four or five, contradicting (8).
This proves (9).
\\
\\
(10) {\em We may assume that for each $t\in V(A)$, every vertex in $X$ is adjacent to an initial segment of $W_t$.}
\\
\\
From (9), there is an ordering of $W_t$ (not necessarily with $t$ as its first term) such that every vertex in $V(G)\setminus W_t$
is adjacent to an initial segment of the ordering. Thus, by replacing $F$ with an isomorphic graph, we may assume that the 
claim of (10) holds.  This proves (10).
\\
\\
(11) {\em Every vertex in $X$ has two nonadjacent neighbours in $V(A)$.}
\\
\\
Let $v\in X$; then $v$ has two nonadjacent neighbours $x,y\in W(A)$. Let $x\in W_s$ and $y\in W_t$, where $s,t\in V(F)$.
Since $x,y$ are nonadjacent, it follows that $s\ne t$, and $s,t$ are not adjacent; and by (10), $v$ is adjacent to both $s,t$.
This proves (11).

\bigskip

Let $F'$ be the $\ell$-frame obtained from $F$ by deleting any apexes, and 
let $H'$ be obtained from $H$ by deleting $W_{a_0}, W_{b_0}$ (if they exist).
Then $H'$ is a blow-up of the $\ell$-frame $F'$. 
From (8), (10), (11), this proves \ref{partition}.~\bbox

\section{Border construction}\label{sec:bordercon}

Now we turn to the second part of the proof of \ref{mainthm}, showing how to make a border.
If $T$ is an \arb{}, then for every vertex $v$ of $T$ different from the apex, there is a unique vertex $u$ of $T$
such that $v$ is adjacent from $u$; we call $u$ the {\em $T$-inneighbour} of $v$.

Let $(J,A,(W_t:t\in V(A)))$ be a border, and let $J'$ be the subgraph of $J$ induced on $V(A)\cup (V(J)\setminus W(A))$. 
We say an \arb{} $T$ is a {\em basis} for 
$(J,A,(W_t:t\in V(A)))$ if:
\begin{enumerate}
\item $\up{T}=J'$.
\item $L(T)\subseteq V(A)$.
\item For each $t\in L(T)$, let $S$ be the path of $T$ from $r(T)$ to the $T$-inneighbour of $t$; 
	then for each $v\in W_t$ there exists a subpath $S_v$ of $S$ with first vertex $r(T)$, 
containing all vertices of $V(S)\cap V(A)$, 
such that 
the set of $J$-neighbours of $v$ in $V(J)\setminus W_t$ is $V(S_v)\cup W(S_v\cap A)$.
\item For each $t\in V(A)\setminus L(T)$, all vertices in $W_t$ have the same $J$-neighbours in $V(J)\setminus W_t$ as $t$.
\item For each $t\in V(A)\setminus \{r(T)\}$, the $T$-inneighbour of $t$ either belongs to $V(A)$ or 
has $T$-outdegree at least two. 
\item $r(T)\in V(A)$ if and only if $A$ is connected.
\end{enumerate}
If so, it follows from the definition of a border that,
with $t,S$, and $S_v\;(v\in W_t)$ as in the third bullet above, if $\{x_1\LL x_n\}$ is the ordering of the ordered 
clique $W_t$, then the length of $S_{x_j}$ monotonically (non-strictly) decreases as $j$ increases.

The first main result of this section is that every border has a basis, but to prove that we need several steps.
Thoughout the remainder of this section,
we will assume that 
$(J,A,(W_t:t\in V(A)))$ is a border, and 
$J'$ is the subgraph of $J$ induced on $V(A)\cup (V(J)\setminus W(A))$.
\begin{thm}\label{fourpaths}
$J'$ does not contain $\mathcal{P}_4$.
\end{thm}
\Proof
From the definition of a border, we have 
\\
\\
(1) {\em No induced path of $J'$ with length at least three has both ends in $V(A)$.}
\\
\\
We deduce
\\
\\
(2) {\em There is no four-vertex induced path $v_1\DD v_2\DD v_3\DD v_4$ of $J'$ with $v_1\in V(A)$.}
\\
\\
Suppose there is such a path. By (1), $v_4\in V(J)\setminus V(A)$ and so has two nonadjacent neighbours $x,y\in V(A)$.

Thus $v_1\ne x,y$, since $v_1,v_4$ are nonadjacent. If both $x,y$ have a neighbour in $\{v_1,v_2\}$, then $J'$ 
contains a 4-hole or 5-hole, a contradiction, since $J$ has no 4-hole or 5-hole from the definition of a border; so we may
assume that neither of $v_1,v_2$ is adjacent to $x$, and so $x\ne v_3$. 
By (1), $v_1\DD v_2\DD v_3\DD v_4\DD x$
is not an induced path, and so $v_3$ is adjacent to $x$. This contradicts (1) applied to $v_1\DD v_2\DD v_3\DD x$, and so proves (2).

\bigskip
Now suppose that $v_1\DD v_2\DD v_3\DD v_4$ is a copy of $\mathcal{P}_4$ in $J'$. By (2), 
$v_1,v_4\notin V(A)$.
Hence $v_4$ has two nonadjacent neighbours $x,y$ in $V(A)$.
As before we may assume that
$v_1,v_2$ are nonadjacent to $x$, and hence $x\ne v_3$.
But then one of $v_1\DD v_2\DD v_3\DD x$, $v_2\DD v_3\DD v_4\DD x$ violates (2). This proves \ref{fourpaths}.~\bbox

We use the following theorem of Wolk~\cite{wolk}:
\begin{thm}\label{laminar}
Let $G$ be a graph. Then $G$ is the transitive closure of some arborescence if and only if $G$ is non-null, connected, and 
does not contain $\mathcal{P}_4$ or $\mathcal{C}_4$.
\end{thm}

The graph $J'$ does not contain $\mathcal{P}_4$ or $\mathcal{C}_4$, by \ref{fourpaths} and the definition of a border,
so by \ref{laminar}, there is 
an arborescence $T$ such that
$J'=\up{T}$. 
Thus, $V(T)=V(J')$, and for every two vertices in $V(T)$, they are 
$J'$-adjacent if and only if they are joined by a directed path of $T$.

Here is a useful observation (the proof is clear, since every vertex has indegree at most one):
\begin{thm}\label{twopaths}
If $P,Q$ are directed paths of an arborescence $T$, and the last vertex of $P$ belongs to $V(Q)$, then $P\cup Q$ is a directed path of $T$.
\end{thm}

\begin{thm}\label{leaves}
Let $T$ be an arborescence with $J'=\up{T}$. Then
$L(T)\subseteq V(A)$.
\end{thm}
\Proof
Suppose that $v\in L(T)\setminus V(A)$.
From the definition of a border, there are $x,y\in V(A)$, not $J'$-adjacent to each other, and $J'$-adjacent to $v$.
Consequently $x,y\in V(T)$, and each of them is joined to $v$ by a directed path of $T$ (since $v$ has zero outdegree in $T$),
say $P$ and $Q$ respectively. By \ref{twopaths},
$P\cup Q$ is a directed path, and so $x,y$ are adjacent in $\up{T}$, a contradiction. This proves \ref{leaves}.~\bbox

\begin{thm}\label{terminal}
Let $T$ be an arborescence with $J'=\up{T}$, and let $t\in V(A)$.
If there exists $v\in V(J)\setminus W_t$ mixed on $W_t$, then $t\in L(T)$.
\end{thm}
\Proof
Suppose that $v\in V(J)\setminus W_t$ is mixed on $W_t$. Hence $v$ is $J$-adjacent to $t$, and $J$-nonadjacent to some $u\in W_t$.
It follows that $v\notin W(A)$, from the definition of a border.
Hence there exist nonadjacent $x,y\in V(A)$ adjacent to $v$. Not both $x,y$ are adjacent to $u$, since otherwise $\{x,y, u,v\}$ induces
a 4-hole. We assume that $x,u$ are nonadjacent.
Consequently $x,t$ are distinct and nonadjacent, since $x$ is not mixed on $W_t$ (because $x\in V(A)$).

Suppose that $t\notin L(T)$, and
let $s\in L(T)$ such that there is a directed path of $T$ from $t$ to $s$.
Since $J' = \overrightarrow{T}$, it follows that $t,s$ are $J'$-adjacent and therefore $J$-adjacent.
Hence by the definition of border, $W_t$ is complete to $W_s$ since $t,s \in V(A)$.
In particular, $u,s$ are adjacent.
Moreover,
$t$ $J'$-dominates $s$; and
so $x,s$ are nonadjacent since $x,s,t\in V(J')$ and $x,t$ are nonadjacent.
Since $x,v,t,s\in V(J')$, and $J'$ contains no four-vertex induced path (by \ref{fourpaths}), it follows that
$v,s$ are $A$-adjacent. But then $x\DD v\DD s\DD u$ is an induced path that violates the
definition of a border.
This proves \ref{terminal}.~\bbox

Let $T$ be an arborescence with $J'=\up{T}$.
Let us say the {\em big cost} of $T$ is the sum, over all $t\in V(A)$, of the $T$-distance between
$r(T)$ and $t$; and the {\em little cost} of $T$ is the sum, over all $v\in W(A)\setminus V(A)$ and all $u\in V(T)$ $J$-adjacent
to $v$, of the $T$-distance between $r(T), u$.
There may be several choices of the arborescence $T$ that have the same transitive closure $\up{T}$. Let us choose $T$
with apex in $V(A)$ if possible; subject to that, with minimum big cost; and subject to that, with minimum little cost.
Such a choice of $T$ is said to be {\em optimal}.

\begin{thm}\label{supercubic}
Let $T$ be an optimal arborescence with $J'=\up{T}$.
Then:
\begin{itemize}
\item  $r(T)\in V(A)$ if and only if $A$ is connected. 
\item If $r(T)\notin V(A)$ and $r(T)$ has a unique $T$-outneighbour, this outneighbour does not belong to $V(A)$.
\item For each $t\in V(A)\setminus \{r(T)\}$, let $v$ be the $T$-inneighbour of $t$; then either $v\in V(A)$ or
the $T$-outdegree of $v$ is at least two.
\item For each $t\in L(T)$, let $S$ be the path of $T$ between $r(T)$ and the $T$-inneighbour of $t$; then for each $v\in W_t$,
there exists a subpath $S_v$ of $S$ with first vertex $r(T)$,
containing all vertices of $V(S)\cap V(A)$,
such that
$V(S_v)$ is the set of $J$-neighbours of $v$ in $V(T)\setminus \{t\}$.
\end{itemize}
\end{thm}
\Proof
Since $r(T)$ is $\up{T}$-adjacent to all other vertices of $T$, it follows that if $r(T)\in V(A)$ then $A$ is connected. Now suppose
that $r(T)\notin V(A)$, and $A$ is connected. Consequently there exists $a\in V(A)$ $A$-adjacent to every other vertex of $A$, 
since $A$ is a threshold graph. Since every vertex $z\in V(T)\setminus V(A)$
has two nonadjacent neighbours in $V(A)$, and $J'$ contains no 4-hole, it follows that $a$ is $J'$-complete to $V(T)\setminus \{a\}$; 
and so $r(T),a$
are adjacent twins of $J'$. There is an automorphism of $J'$ exchanging $r(T),a$ and fixing all other vertices of $J'$. 
Let $T'$ be the image of $T$
under this automorphism. Then $T'$ is an arborescence with $\up{T'}=\up{T}$, and with apex in $V(A)$, contrary to the optimality of $T$. 
This proves the first assertion of the theorem.

For the second, suppose that $r(T)\notin V(A)$ and that $r(T)$ has a unique $T$-outneighbour $a_1$ say, and $a_1\in V(A)$.
Thus 
$a_1$ is $\up{T}$-adjacent to all other vertices in $A$, and so $A$ is connected, contrary to the first assertion.

For the third assertion, let $t\in V(A)\setminus \{r(T)\}$, and let $v$ be the $T$-inneighbour of $t$. 
We may assume that 
$v\notin V(A)$. 
Suppose that $t$ is the unique $T$-outneighbour of $v$. Consequently $t, v$ are 
adjacent $J'$-twins; and replacing $T$ by the image of $T$ under the automorphism of $J'$
that exchanges $v, t$ and fixes all other vertices, contradicts the optimality of $T$. This proves the third bullet.

Finally, let $t\in L(T)$ and $v\in W_t\setminus \{t\}$, and let $S$ be the path of $T$ between $r(T)$ and the $T$-inneighbour of $t$. Since $t$ $J$-dominates $v$
and all $J$-neighbours of $t$ in $V(T)$ belong to $V(S)$, it follows that all $J$-neighbours of $v$ in $V(T)$ belong to $V(S)$.
We claim that this set of neighbours, excepting $t$, is the vertex set of a subpath of $S$ starting at $r(T)$. Suppose not; then 
there exist $u,w\in V(S)\setminus \{t\}$, such that $w$ is a $T$-outneighbour of $u$, and $v$ is $J$-adjacent to $w$ and not to $u$.

Suppose that $u$ has $T$-outdegree more than one, and let $w'\in V(T)$ be a $T$-outneighbour of $u$ different from $w$. Let $t'\in L(T)$
such that there is a directed path of $T$ from $w'$ to $t'$. It follows that $w,t'$ and $t,t'$ are not $J$-adjacent, since there is no directed path
of $T$ between them; and hence $v,t'$ are not $J$-adjacent, since $t'$ is not mixed on $W_t$. But $t'\in V(A)$ by \ref{leaves}; and so $v\DD w\DD u\DD t'$ is an induced
path of $J$, contradicting the definition of a border.
Thus $w$ is the unique $T$-outneighbour of $u$.

Now suppose that there exists $v'\in W(A)\setminus V(A)$ that is $J$-adjacent to $u$ and not to $w$. If $v,v'$ are $J$-adjacent,
then $\{v,v',u,w\}$ induces a 4-hole, and if not then $v\DD w\DD u\DD v'$ violates the definition of a border. Thus there is no such
$v'$. Let $T'$ be the image of $T$ under the automorphism of $J'$
that exchanges $u,w$ and fixes all other vertices. The big cost of $T'$ equals that of $T$, but its little cost is strictly 
smaller, contrary to the optimality of $T$.  

This proves that there is a subpath $S_v$ of $S$ with first vertex $r(T)$, such that $V(S_v)$ is the set of $J$-neighbours of 
$v$ in $V(T)\setminus \{t\}$. Since each vertex in $V(S)\cap V(A)$ is adjacent to $t$ and is not mixed on $W_t$, it follows that
each such vertex belongs to $S_v$. 
This proves \ref{supercubic}.~\bbox

Now we can prove the first main result of this section:
\begin{thm}\label{newforcecon}
Let $(J,A,(W_t:t\in V(A)))$ be a border. Then there is a basis $T$ for $(J,A,(W_t:t\in V(A)))$. 
\end{thm}
\Proof
With $J'$ as before, let $T$ be an optimal arborescence with $J'=\up{T}$. (This exists by \ref{fourpaths} and \ref{laminar}.)
We must check that:
\begin{enumerate}
\item $\up{T}=J'$.
\item $L(T)\subseteq V(A)$.
\item For each $t\in L(T)$, let $S$ be the path of $T$ from $r(T)$ to the $T$-inneighbour of $t$;
        then for each $v\in W_t$ there exists a subpath $S_v$ of $S$ with first vertex $r(T)$,
containing all vertices of $V(S)\cap V(A)$,
such that
the set of $J$-neighbours of $v$ in $V(J)\setminus W_t$ is $V(S_v)\cup W(S_v\cap A)$.
\item For each $t\in V(A)\setminus L(T)$, all vertices in $W_t$ have the same $J$-neighbours in $V(J)\setminus W_t$ as $t$.
\item For each $t\in V(A)\setminus \{r(T)\}$, the $T$-inneighbour of $t$ either belongs to $V(A)$ or
has $T$-outdegree at least two.
\item $r(T)\in V(A)$ if and only if $A$ is connected.
\end{enumerate}
The first evidently holds, and the second by \ref{leaves}. 
The fourth statement follows from \ref{terminal}.
For the third,
let $t\in L(T)$, let $S$ be the path of $T$ between $r(T)$ and the $T$-inneighbour of $t$, and let $v\in W_t$.
By \ref{supercubic}, 
there exists a subpath $S_v$ of $S$ with first vertex $r(T)$,
containing all vertices of $V(S)\cap V(A)$,
such that
$V(S_v)$ is the set of $J$-neighbours of $v$ in $V(T)\setminus \{t\}$. From the fourth statement above,
the set of $J$-neighbours of $v$ in $V(J)\setminus W_t$ equals $V(S_v)\cup W(S_v\cap A)$.
The fifth and sixth statements follow from \ref{supercubic}. This proves \ref{newforcecon}.~\bbox

A {\em split} of a graph $G$ is a partition $(X,Y)$ of $V(G)$ into a clique $X$ and a stable set $Y$. 
With notation as before, if $(X,Y)$ is a split of $A$, and $T$ is an arborescence with $\up{T}=J'$, 
an {\em $(X,Y)$-spine of $T$} is a directed path $R$ of $T$, with first vertex $r(T)$
and last vertex in $X\cup \{r(T)\}$, such that
$X=V(A)\cap V(R)$ and $Y=L(T)\setminus V(R)$.

If $(X,Y)$ is a split of $A$, there may be bases that have no $(X,Y)$-spine, but the next result shows that there is a 
basis that does have such a spine.
\begin{thm}\label{spine}
Let $(J,A, (W_t:t\in V(A)))$ be a border, and let 
$(X,Y)$ be a split of $A$.
Then there is a basis for $(J,A, (W_t:t\in V(A)))$ that has an $(X,Y)$-spine.
\end{thm}
\Proof
As before, let $J'$ be the subgraph of $J$ induced on $V(A)\cup (V(J)\setminus W(A))$.
By \ref{newforcecon} there is a basis $T$ for $(J,A, (W_t:t\in V(A)))$, and therefore $J' = \overrightarrow{T}$.
Choose $T$ with $|L(T) \cap Y|$ maximum. We claim that $T$ has an $(X,Y)$-spine.

Choose a directed path $R$ of $T$ with first vertex $r(T)$ and last vertex in $X\cup \{r(T)\}$, as long as possible, and let
$x$ be the last vertex of $R$. We claim that $X\subseteq V(R)$; for suppose that
there exists $x'\in X\setminus V(R)$. Since $x,x'$ are $\up{T}$-adjacent, there is a directed path $Q$ of $T$ between $x,x'$,
and by \ref{twopaths} $Q\cup R$ is a directed path, contrary to the maximality of the length of $R$. This proves that $X\subseteq V(R)$.

Now we must show that $Y=L(T)\setminus V(R)$.
Certainly $Y\supseteq L(T)\setminus V(R)$, since if $t\in L(T)\setminus V(R)$, then $t\in V(A)$ by \ref{leaves}, and so $t\in Y$, since $t\notin V(R)\supseteq X$. 

Suppose that there exists $t\in Y\setminus (L(T)\setminus V(R))$.
Now $t\notin L(T)\cap V(R)$, since $L(T)\cap V(R)\subseteq L(T)\cap \{x\}\subseteq X$, and $t\notin X$. Consequently
$t\notin L(T)$. Let $R'$ be a directed path of $T$ from $t$ to some vertex $s\in L(T)$.
\\
\\
(1) {\em $R'$ is a subpath of $R$; and all its vertices have $T$-outdegree at most one; and all its vertices belong to $V(A)$.
Moreover, $t\in V(R)$, and $x=s\in L(T)$.}
\\
\\
By the definition of basis, $s \in V(A)$. 
Since $t\ne s$, and $t,s$ are $\up{T}$-adjacent and hence $A$-adjacent (because $s,t\in V(A)$), and $t\in Y$, and $Y$ is $A$-stable, 
it follows that $s\notin Y$. So $s\in X\subseteq V(R)$, and consequently $R'$ is a subpath of $R$, and in particular $t\in V(R)$, and
$x=s\in L(T)$. The uniqueness of $R'$ implies that each of its vertices has $T$-outdegree at most one; and from the fifth bullet
in the definition of a basis, $V(R')\subseteq V(A)$. This proves (1).

\bigskip

By (1), there is an automorphism of $J'$ that exchanges $t,x$ and fixes all other elements of $V(J')$; let $T'$ be the image of $T$ under this automorphism.
We claim that $T'$ is a basis for $(J,A, (W_t:t\in V(A)))$.
To show this, we need to check the six conditions in the definition of a basis, namely:
\begin{enumerate}
\item $\up{T'}=J'$.
\item $L(T')\subseteq V(A)$.
\item For each $t'\in L(T')$, let $S$ be the path of $T'$ from $r(T')$ to the $T'$-inneighbour of $t'$;
        then for each $v\in W_{t'}$ there exists a subpath $S_v$ of $S$ with first vertex $r(T')$,
containing all vertices of $V(S)\cap V(A)$,
such that
the set of $J$-neighbours of $v$ in $V(J)\setminus W_{t'}$ is $V(S_v)\cup W(S_v\cap A)$.
\item For each $t'\in V(A)\setminus L(T')$, all vertices in $W_{t'}$ have the same $J$-neighbours in $V(J)\setminus W_{t'}$ as $t'$.
\item For each $t'\in V(A)\setminus \{r(T')\}$, the $T'$-inneighbour of $t'$ either belongs to $V(A)$ or
has $T'$-outdegree at least two.
\item $r(T')\in V(A)$ if and only if $A$ is connected.
\end{enumerate}
Since $\up{T'} = \up{T}$, the first condition holds; and since $y\in V(A)$, and $L(T)\subseteq V(A)$, the second condition holds. 
To check the third condition, we may assume that $t'=t$; but 
from \ref{terminal} applied to $T$, there is no vertex in $V(J)\setminus W_t$ mixed on $W_t$; and so the third condition holds, taking
$S_v=S$ for each $v\in W_t$.

For the fourth condition, we may assume that $t'=x$; but by \ref{terminal} applied to $T'$, there is no vertex in $V(J)\setminus W_x$ mixed on $W_x$, and so the fourth condition holds.

To check the fifth condition, we may assume that the $T'$-inneighbour of $t'$ is different from its $T$-inneighbour, and hence 
$t'$ is one of $t,x$ or the $T$-outneighbour $t''$ (say) of $t$. For $t$ and $t''$, its $T'$-inneighbour belongs to $V(A)$. For
$x$, its $T'$-inneighbour is the $T$-inneighbour of $t$, and so either belongs to $V(A)$ or has $T$-outdegree (and hence $T'$-outdegree) at least two.

Finally, for the last condition, we may assume that $r(T')\ne r(T)$ and so $t=r(T)$; but then $r(T') = x\in V(A)$, and $A$ is a complete graph and hence connected.
This proves that $T'$ is a basis for $(J,A, (W_t:t\in V(A)))$; but this contradicts the choice of $T$, and so proves \ref{spine}.~\bbox

\bigskip

We would like it to be true that for each $t\in X$, no vertex of $J$ is mixed on $W_t$. From the definition of a basis, this is true 
if $t\notin L(T)$, but we must be careful if the last vertex of $R$ is a leaf of $T$.
We will have two applications of \ref{spine}, and in both we will have to check separately that 
if the last vertex of $R$ is a leaf of $T$, then no vertex of $J$ is mixed on $W_t$.

\section{The case when $\ell$ is odd}\label{sec:oddcase}

In this section we complete the proof of \ref{mainthm} when $\ell$ is odd. 
If $t$ is a vertex of a graph $A$, we say that $t$ is {\em $A$-non-terminal} (or just ``non-terminal'')
if there exists $s\in V(A)\setminus \{t\}$
such that $t$ $A$-dominates $s$, and {\em $A$-terminal} (or just ``terminal'') otherwise.
We need:

\begin{thm}\label{termcomp}
Let $A$ be a threshold graph with $|A|>1$. Then $t\in V(A)$ is $A$-terminal if and only if $t$ is not
$\overline{A}$-terminal.
Consequently every two $A$-terminal vertices are nonadjacent, and every two $A$-non-terminal vertices in $A$ are adjacent.
\end{thm}
\Proof
Let $N$ be the set of neighbours of $t$ in $A$, and $M$ its
set of non-neighbours. Thus $\{t\}, N,M$ partition $V(A)$. 
Suppose that there
exist $n_1,n_2\in N$ and $m_1,m_2\in M$ such that $m_1n_1$ and $m_2n_2$ are edges of $A$, and $m_1$ is not $A$-adjacent to $n_2$,
and $m_2$ is not $A$-adjacent to $n_1$. Since $m_1\DD n_1\DD t\DD n_2$ is not an induced path it follows that $n_1,n_2$
are $A$-adjacent; but then $A[\{m_1,n_1,m_2,n_2\}]$ is a copy of $\mathcal{P}_4$ or $\mathcal{C}_4$, a contradiction. Thus there are no such
$m_1,n_1,m_2,n_2$, and so 
the bipartite graph $A[N,M]$ is a half-graph. 

Consequently either there exists $m\in M$ adjacent to every vertex in $N$, or there exists $n\in N$ with no neighbour in $M$.
In the first case, $t$ is $A$-terminal, since $m$ is adjacent to all neighbours of $t$ and not adjacent to $t$;
and $t$ is $\overline{A}$-non-terminal, since it $\overline{A}$-dominates $m$. In the second case, the same argument
applied in $\overline{A}$ shows that $t$ is $\overline{A}$-terminal and $t$ is $A$-non-terminal. This proves the first claim.

Suppose that $t_1,t_2\in V(A)$ are $A$-terminal and $A$-adjacent. Since neither $A$-dominates the other, there exists
$s_1,s_2\in V(A)\setminus \{t_1,t_2\}$, such that $s_1$ is adjacent to $t_1$ and not to $t_2$, and $s_2$
is adjacent to $t_2$ and not to $t_1$. In particular $s_1\ne s_2$; and $\{s_1,s_2,t_1,t_2\}$ induces $\mathcal{P}_4$ or $\mathcal{C}_4$,
a contradiction. This proves that $A$-terminal vertices are nonadjacent. The final statement follows by taking complements.
This proves \ref{termcomp}.~\bbox

By \ref{easymainthm}, to complete the proof of \ref{mainthm} when $\ell$ is odd, it suffices to prove the following:
\begin{thm}\label{mainodd}
Let $\ell\ge 7$ be odd, and let $G$ be an $\ell$-holed graph with no clique cutset or universal vertex.
Then $G$
is either a blow-up of a cycle of length $\ell$, or a blow-up of an $\ell$-framework.
\end{thm}
\Proof
By \ref{partition}, we may assume that $G$ is the bordered blow-up of an $\ell$-frame, and $G$ is the composition of $H,J,K$
where:
\begin{itemize}
\item $H$ is a blow-up of an $\ell$-frame $F$, where $F$ has sides $A,B$;
\item $(J, A, (W_t:t\in V(A)))$ and
$(K,B,(W_t:t\in V(B)))$ are borders; and
\item $V(H\cap J) = W(A)$, and $V(H\cap K)=W(B)$, and $V(J),V(K)$ are disjoint and anticomplete.
\end{itemize}
Let $F$ have $k$ bars $P_1\LL P_k$, where $P_i$ has ends $a_i\in V(A)$ and $b_i\in V(B)$, and so
$|A|=|B|=k\ge 3$. 
Since $A,B$ are complementary threshold graphs, exactly one is disconnected, say $A$. 
We will show that there is a choice of $T,S$ such that $S\cup T\cup P_1\cupcup P_k$ is an $\ell$-framework, and
$G$ is a blow-up of this framework.

Let $X$ be the set of all non-terminal vertices of $A$, and $Y=V(A)\setminus X$; thus $(X,Y)$ is a split of
$A$, by \ref{termcomp}. By \ref{spine}, there is a basis $T$ for $(J,A,(W_t:t\in V(A)))$ with an $(X,Y)$-spine $R$.
Let $a_0=r(T)$, and let the vertices of $R$ in $V(A)\cup \{a_0\}$
be $\{a_0,a_1\LL a_m\}$, numbered in order on $R$. Thus $X=\{a_1\LL a_m\}$ since $a_0=r(T)\notin V(A)$ (because $A$ is disconnected).
Let $X'=\{b_1\LL b_{m}\}$, and $Y'=V(B)\setminus X'$. Since $a_i, a_j$ are adjacent if and only if $b_i, b_j$ are nonadjacent
for $1\le i<j\le k$,
\ref{termcomp} implies that $Y'$ is the set of non-terminal vertices of $B$.
By \ref{spine} there is a basis $S$ for $(K,B,(W_t:t\in V(B)))$
with an $(Y',X')$-spine $R'$.
Hence $V(R')\cap V(B)=\{b_{m+1}\LL b_k\}$,
and so we may assume (by renumbering $b_{m+1}\LL b_k$) that $b_k$ is the apex of $S$ and $b_k, b_{k-1}\LL b_{m+1}$
are in order in $R'$.

For $0\le i\le m$ let $I_i$ be the set of vertices $v\in \{a_{m+1}\LL a_k\}$ such that there is a directed path of $T$ from 
$a_i$ to $v$, not using $a_{i+1}$ if $i<m$.
For $m+1\le j\le k$, let $L_j$ be the set of $v\in \{b_1\LL b_m\}$
such that there is a directed path of $S$ from $b_j$ to $v$, not using $b_{j-1}$ if $j>m+1$.
The sets $I_0\LL I_m$ are pairwise disjoint and have union 
$\{a_{m+1}\LL a_k\}$, the set of terminal vertices of $A$; and $L_{m+1}\LL L_k$ are pairwise disjoint and have union 
$\{b_{1}\LL b_m\}$, the set of terminal vertices of $B$. 
\\
\\
(1) {\em For $1\le i\le m$ and $m+1\le p\le k$,  $a_p\in I_0\cup I_{1}\cupcup I_{i-1}$ if and only if $b_i\in L_{m+1}\cupcup L_p$.}
\\
\\
$a_i$ and $a_p$ are adjacent if and only if $a_p\in I_i\cup I_{i+1}\cupcup I_m$;
and $b_p, b_i$ are adjacent if and only if $b_i\in L_{m+1}\cup L_{m+2}\cupcup L_p$. Since $a_i,a_p$ are adjacent if and only 
if $b_i, b_p$ are nonadjacent, this proves (1).
\\
\\
(2) {\em For $0\le i<j\le m$, if $a_p\in I_i$ and $a_q\in I_j$ then $p>q$, and consequently $I_0, I_1\LL I_m$ are 
(possibly null) intervals of $\{a_{m+1}\LL a_k\}$.  Similarly,
for $m+1\le p<q\le k$, if $b_i\in L_p$ and $b_j\in L_q$ then
$i>j$, and so $L_{m+1}\LL L_k$ are (possibly null) intervals of $\{b_1\LL b_m\}$. (See figure \ref{fig:oddframework}.)}
\\
\\
Let $0\le i<j\le m$, and $a_p\in I_i$ and $a_q\in I_j$. Since $a_q\notin I_0\cupcup I_{j-1}$, and $a_p\in I_0\cupcup I_{j-1}$,
(1) implies that $b_j\notin L_{m+1}\cupcup L_q$ and $b_j\in L_{m+1}\cupcup L_p$, and so $p>q$. This proves the first statement.

Now let  $m+1\le p<q\le k$ and $b_i\in L_p$ and $b_j\in L_q$. Thus $i,j\ge 1$. Since 
$b_j\notin L_{m+1}\cupcup L_{p}$, and $b_i\in L_{m+1}\cupcup L_{p}$,
(1) implies that $a_{p}\notin I_0\cup I_{1}\cupcup I_{j-1}$ and $a_p\in I_0\cup I_{1}\cupcup I_{i-1}$, and so $i>j$. This proves the second statement,
and so proves (2).
\\
\\
(3) {\em $I_0\ne \emptyset$. Let $0\le i\le m$ with $I_i\ne \emptyset$, and let $I_i$ be the interval $\{a_p\LL a_q\}$. 
Then $L_{p+1}\LL L_{q}$ are all 
empty, and $L_p \ne \emptyset$ unless $p=m+1$.}
\\
\\
Since $a_0\notin V(A)$, the definition of a basis implies that the subpath of $R$ from $a_0$ to $a_1$ has a vertex 
different from $a_1$ with $T$-outdegree at least two. Consequently $I_0\ne \emptyset$.
Now let $0\le i\le m$ with $I_i\ne \emptyset$, let $I_i$ be the interval $\{a_p\LL a_q\}$, 
and suppose that $p< q'\le q$ and $L_{q'}\ne \emptyset$. Let $b_j\in L_{q'}$. 
If $i< j$, then 
$$a_p\in I_0\cupcup I_i\subseteq I_0\cupcup I_{j-1}$$
and 
$b_j\notin L_{m+1}\cupcup L_{p}$, contrary to (1). If $i\ge j$, then 
$a_{q'}\notin I_0\cupcup I_{j-1}$ and $b_j\in L_{m+1}\cupcup L_{q'}$, contrary to (1). 
This proves that $L_{p+1}\LL L_{q}$ are all
empty.

Now suppose that $p>m+1$. Since $a_{m+1}$ belongs to $I_0\cupcup I_m$ and $a_{m+1}\notin I_i$, it follows from (2) that $i<m$.
Since $a_{p-1}\notin I_0\cupcup I_i$, and $a_p\in I_0\cupcup I_i$, it follows from (1) that 
$b_{i+1}\notin L_{m+1}\cupcup L_{p-1}$, and $b_{i+1}\in L_{m+1}\cupcup L_p$. Consequently $L_p\ne \emptyset$. 
This proves (3).
\\
\\
(4) {\em Let $m+1\le p\le k$ with $L_p\ne \emptyset$, and let $L_p$ be the interval $\{b_i\LL b_j\}$. Then $I_i\LL I_{j-1}$
are all empty, and $I_j\ne \emptyset$ unless $j=m$.}
\\
\\
Let $i\le i'<j$, and suppose that $I_{i'}\ne \emptyset$. Let $a_q\in I_{i'}$. If $p\le q$, then $a_{q}\notin I_0\cupcup I_{i'-1}$,
and $b_{i'}\in L_{m+1}\cupcup L_q$, contrary to (1). If $p>q$, then $a_q\in I_0\cupcup I_{j-1}$, and $b_j\notin L_{m+1}\cupcup L_q$,
contrary to (1). Thus $I_i\LL I_{j-1}$
are all empty. Now suppose that $j\ne m$. Since $b_{j+1}\in L_{m+1}\cupcup L_{p-1}$ and 
$b_j\notin L_{m+1}\cupcup L_{p-1}$, (1) implies that
$a_p\in I_0\cupcup I_j$, and $a_p\notin I_0\cupcup I_{j-1}$; and so $I_{j}\ne \emptyset$. This proves (4).

\bigskip

From (2), (3), (4) it follows that $S\cup T\cup P_1\cupcup P_k$ is an $\ell$-framework. 
\\
\\
(5) {\em For each $t\in X$, no vertex of $J$ is mixed on $W_t$; and for each $t\in Y'$, no vertex of $K$ is mixed on $W_t$.}
\\
\\
From the choice of $R$ it follows that $t\in V(R)$; and the claim is true from the definition of a basis if $t\notin L(T)$;
so we assume that $t\in L(T)$, and hence is the last vertex of $R$,
and therefore $m>0$ and $t=a_m$.
Since $a_m\in X$, it follows that $a_m$ is not terminal in $A$, and so
 $a_m$ $A$-dominates some $a_i$ where $0\le i\le k$ and $i\ne m$. Since $a_m\in L(T)$ and $L(T)$ is stable in $A$, it follows that 
$a_i\notin L(T)$, and so $i<m$.
Let $u$ be the $T$-inneighbour of $a_m$, and suppose that $u\ne a_{i}$.
From the definition of a basis, $u$ has $T$-outdegree at least two, and so there is a leaf $s\ne a_m$ of $T$
such that there is a directed path of $T$ from $u$ to $s$. But then $a_{i}, s$ are $\up{T}$-adjacent and hence $A$-adjacent,
and $a_m,s$ are not $A$-adjacent since they both belong to $L(T)$. This contradicts that $a_m$ $A$-dominates $a_{i}$.

Consequently $a_{i}$ is the $T$-inneighbour of $a_m$, and therefore $i=m-1$. But then, from the definition of a basis, every vertex in $W_t$
is $J$-adjacent to $V(S)\cup W(A\cap S)$, where $S$ is the path of $T$ from $a_0$ to $a_{m-1}$, and has no other $J$-neighbours
in $V(J)\setminus W_t$. This proves the first assertion, and the second follows similarly. This proves (5).

\bigskip

Since $G$ is the composition of $H,J,K$, and $T$ is a basis for $(J,A,(W_t:t\in V(A)))$ (and the same for $K$),
(5) implies that $G$ is a blow-up of the $\ell$-framework $S\cup T\cup P_1\cupcup P_k$.
This proves \ref{mainodd}.~\bbox
\section{A structure theorem for laminar families}

Now we turn to the case of \ref{mainthm} when $\ell$ is even. We will need a 
theorem about bipartite graphs that is proved in this section.

Let $G$ be a bipartite graph, with bipartition $(A,B)$. For $A'\subseteq A$, we say {\em $A'$ is laminar in $G$} if for all distinct $a,a'\in A'$,
either $N(a)\subseteq N(a')$, or $N(a')\subseteq N(a)$, or $N(a)\cap N(a')=\emptyset$, where $N(v)$ denotes the set
of neighbours of a vertex $v\in A$.

The result we will use later is:
\begin{thm}\label{twofam}
Let $G$ be a bipartite graph with bipartition $(A,B)$, and let $A$ be the union of disjoint sets $A_1,A_2$. Suppose that
$A_1,A_2$ are both laminar in $G$. Then the following are equivalent:
\begin{itemize}
\item every hole of $G$ has length four;
\item there is a tree $T$ with $V(T)=B$, such that for each $a\in A$, $N(a)$ is the vertex set of a subtree of $T$.
\end{itemize}
\end{thm}

That is a consequence of the following (proved later in this section):
\begin{thm}\label{subtrees}
Let $G$ be a bipartite graph with bipartition $(A,B)$, where $B\ne \emptyset$ and every vertex in $A$ has positive degree.
Then the following are equivalent:
\begin{itemize}
\item for every hole $C$ of $G$ of length at least six, some vertex in $B$ has at least three neighbours in $V(C)$;
\item there is a tree $T$ with $V(T)=B$, such that for each $a\in A$, $N(a)$ is the vertex set of a subtree of $T$.
\end{itemize}
\end{thm}

We say that $G$ admits a {\em 1-join} $(V_1,V_2)$ if $V_1,V_2$ is a partition of $V(G)$ with $|V_1|, |V_2|\ge 2$, and there are subsets
$X_1\subseteq V_1$ and $X_2\subseteq V_2$, such that $X_1$ is complete to $X_2$, and there are no other edges between $V_1,V_2$.
Let us define $N[a]=N(a)\cup \{a\}$, where $N(a)$ is the set of neighbours of $a$.
To prove \ref{subtrees}, we need the following (thanks to Maria Chudnovsky for this short proof; our original proof was 
much longer.)
\begin{thm}\label{1join}
Let $G$ be a connected bipartite graph, with bipartition $(A,B)$, and with $|A|,|B|\ge 2$, that admits no 1-join.
Then there exists $a\in A$ such that $G\setminus N[a]$ is connected.
\end{thm}
\Proof
If $G$ is complete bipartite, then it admits a 1-join $(A,B)$, since $|A|,|B|\ge 2$. Consequently there exist 
$a\in A$ and $b\in B$ with distance at least three. 
Hence there exist $A_1\subseteq A$ and $B_1\subseteq B$ with $G[A_1\cup B_1]$ connected and with $A_1,B_1$ nonnull, 
such that some vertex in 
$A\setminus A_1$ has no neighbour in $B_1$. Choose $A_1,B_1$ with this property, with $A_1\cup B_1$ maximal.

Let $A_2,B_2$ be respectively the sets of vertices in $A\setminus A_1$ and in $B\setminus B_1$ that have a neighbour 
in $A_1\cup B_1$. Let $A_3 = A\setminus  (A_1\cup A_2)$ and $B_3=B\setminus (B_1\cup B_2)$. By hypothesis, $A_3\ne \emptyset$. 
From the maximality of $A_1\cup B_1$ it follows that $A_2=\emptyset$ (because otherwise we could add $A_2$ to $A_1$).
If some vertex $b_2\in B_2$ is nonadjacent to some vertex in $A_3$, then some vertex in 
$A\setminus A_1$ has no neighbour in $B_1\cup \{b_2\}$; so we can add $b_2$ to $B_1$, contrary to the maximality of $A_1\cup B_1$.
So $B_2$ is complete to $A_3$. The only edges between $A_3\cup B_3$ and $A_1\cup A_2\cup B_1\cup B_2$ are those between 
$A_3$ and $B_2$. Since $|A_1\cup B_1|\ge 2$ and $G$ admits no 1-join, it follows that $|A_3\cup B_3|\le 1$. Since
$A_3\ne \emptyset$, it follows that $|A_3|=1$, $A_3=\{a_3\}$ say, and $B_3=\emptyset$. But then $G\setminus N[a_3]$ is connected.
This proves \ref{1join}.~\bbox

We deduce \ref{subtrees}, which we restate:
\begin{thm}\label{subtrees2}
Let $G$ be a bipartite graph with bipartition $(A,B)$, where $B\ne \emptyset$ and every vertex in $A$ has positive degree.
Then the following are equivalent:
\begin{itemize}
\item for every hole $C$ of $G$ of length at least six, some vertex in $B$ has at least three neighbours in $V(C)$;
\item there is a tree $T$ with $V(T)=B$, such that for each $a\in A$, $N(a)$ is the vertex set of a subtree of $T$.
\end{itemize}
\end{thm}
\Proof
Suppose that the second bullet holds; we will prove the first.
Let $T$ be a tree as in the second bullet.
Let $C$ be a hole of $G$ of length at least six, with vertices 
$$b_1\DD a_1\DD b_2\DD a_2\CC b_k\DD a_k\DD b_1$$
for some $k\ge 3$, where $a_1\LL a_k\in A$ and $b_1\LL b_k\in B$. For $1\le i\le k$ let $T_i$ be the subtree of $T$ with vertex set $N(a_i)$. Thus for $1\le i\le k$, $b_i$ belongs to the trees $T_i, T_{i-1}$ but to no other of $T_1\LL T_k$. (We read the subscripts 
modulo $k$.) We need to show that some vertex of $T$ belongs to three of the trees $T_1\LL T_k$, to prove that the first bullet holds.

Let $P$ be the path of $T$ with ends $b_1,b_2$. Thus $P\subseteq T_1$.
Let $e$ be an edge of $P$. If $e$ is not an edge of any of $T_2\LL T_k$, then each of $T_2\LL T_k$
is a subtree of one of the two components of $T\setminus e$, and so $\{T_2\LL T_k\}$ may be partitioned into two nonempty subsets 
$X,Y$, with $T_k\in X$ (because it contains $b_1$) and $T_2\in Y$ (because it contains $b_2$),
such that each tree in the first subset is disjoint from each tree in the second, a contradiction. Thus each edge of $P$
belongs to one of $T_2\LL T_k$. They do not all belong to the same tree $T_i$ where $2\le i\le k$, since otherwise $b_1,b_2\in V(T_i)$
contradicting that $k\ge 3$. So there are two consecutive edges of $P$ that belong to different trees in the list 
$T_2\LL T_k$. But then the common end of these two edges belongs to three of the trees $T_1\LL T_k$, and so the first bullet holds.

We prove the converse implication by induction on $|V(G)|$. Let $G$ be a bipartite graph with bipartition $(A,B)$, where
$B\ne \emptyset$ and every vertex in $A$ has a neighbour in $B$, such that 
for every hole $C$ of $G$ of length at least six, some vertex in $B$ has at least three neighbours in $V(C)$.
We must show that the second bullet holds.
From the inductive hypothesis we may assume that $G$ is connected. The result is easy if 
$|A|\le 1$ or $|B|=1$, so we assume that $|A|,|B|\ge 2$. 
\\
\\
(1) {\em We may assume that $G$ does not admit a 1-join.}
\\
\\
Suppose that $G$ admits a 1-join $(V_1,V_2)$. Let $X_i\subseteq V_i$ for $i = 1,2$, such that $X_1$ is complete to $X_2$ and there are no other edges between $V_1,V_2$.
Thus $X_1,X_2\ne \emptyset$ since $G$ is connected; and so, since $G$ is bipartite, $X_1$ is a subset of one of $A,B$, and $X_2$ of the other.
We may assume that $X_1\subseteq A$ and $X_2\subseteq B$.

Take a new vertex $x_2$ and add it to the graph $G[V_1]$, making it adjacent to the vertices in $X_1$.
Let this graph just made be $G_1$; it admits a bipartition $(V_1\cap A, (V_1\cap B)\cup \{x_2\})$, and it is connected
since $G$ is connected. Since $G_1$ is isomorphic to an induced
subgraph of $G$, and $(V_1\cap B)\cup \{x_2\}\ne \emptyset$ and every vertex in $V_1\cap A$ has a neighbour 
in $(V_1\cap B)\cup \{x_2\}$ (since $G_1$ is connected), the inductive hypothesis implies that there is a tree $T_1$ with vertex set $(V_1\cap B)\cup \{x_2\}$,
such that for each $a\in V_1\cap A$, $N_{G_1}(a)$ is the vertex set of a subtree of $T_1$.

Take a new vertex $x_1$ and add it to the graph $G[V_2]$, making it adjacent to the vertices in $X_2$.
Let this graph just made be $G_2$; it is connected, and admits a bipartition $((V_2\cap A)\cup \{x_1\}, V_2\cap B)$. 
Since it is isomorphic to an induced
subgraph of $G$, and $V_2\cap B\ne \emptyset$ (because it contains $X_2$) and every vertex in $(V_2\cap A)\cup \{x_1\}$ 
has a neighbour
in $V_2\cap B$, the inductive hypothesis implies that there is a tree $T_2$ with vertex set $V_2\cap B$,
such that for each $a\in (V_2\cap A)\cup \{x_1\}$, $N_{G_2}(a)$ is the vertex set of a subtree of $T_2$. Let $T$
be the tree obtained from the disjoint union of $T_1,T_2$ by identifying $x_2$ ($\in V(T_1)$) and some vertex 
($x_2'$ say) of $X_2$ 
($\subseteq V(T_2)$). Then for each $a\in A$, $N(a)$ is the vertex set of a subtree of $T$. To see this, it is clear if
$a\in V_2\cap A$, or $a\in (V_1\cap A)\setminus X_1$. If $a\in X_1$, its neighbour set in $T_1$ forms a subtree of $T_1$
containing $x_2$, and its neighbour set $X_2$ in $T_2$ is the vertex set of a subtree of $T_2$; and the union of these 
(after identifying $x_2,x_2'$) is a subtree of $T$.  This proves (1).

\bigskip
From \ref{1join} and (1), there exists $b\in B$ such that $G\setminus N[b]$ is connected. We may assume that every vertex in $A$ has
degree at least two, because if one of them has degree one, the result follows easily by deleting it and applying the inductive hypothesis.
\\
\\
(2) {\em For all distinct $a,a'\in N(b)$, there is a vertex in $B\setminus \{b\}$ adjacent to them both.}
\\
\\
Suppose that $a,a'\in N(b)$ have no common neighbour in $B\setminus \{b\}$. But $a,a'$ each have a neighbour 
in $B\setminus \{b\}$; and since the graph $G\setminus N[b]$ is connected, there is a path $P$ with interior in $V(G)\setminus N[b]$
and with ends $a,a'$. 
Take the shortest such path; then it is induced, and has length at least four, since $a,a'$ have no common neighbour
in $B\setminus \{b\}$. Now adding $b$ and the edges $ba,ba'$ to $P$ makes a hole in $G$ of length at least six; and so from
the hypothesis, some vertex $b'\in B$ has at least three neighbours in this hole, and hence in $V(P)$. But $b'\ne b$, since
the interior of $P$ is a subset of $V(G)\setminus N[b]$; and rerouting through $b'$ the subpath of $P$ between the first and last 
neighbours on $P$, gives a path shorter than $P$, a contradiction. This proves (2).

\bigskip
Since $B\setminus \{b\}\ne \emptyset$
and every vertex in $A$ has a neighbour in $B\setminus \{b\}$, the inductive hypothesis implies that
there is a tree $S$ with vertex set $B\setminus \{b\}$, such that for each $a\in A$, its neighbour set in $B\setminus \{b\}$ 
is the vertex set of a subtree of $S$.  Let this subtree be $S_a$.

Since the subtrees $S_a\;(a\in N(b))$ pairwise intersect by (2), the Helly property of subtrees
of a tree implies that 
there exists $b'\in B\setminus \{b\}$ that belongs to all of them. Let $T$ be obtained from $S$ by adding 
the new vertex $b$ and an edge $bb'$; then $T$ satisfies the second bullet. 
This proves \ref{subtrees2}.~\bbox

Now we deduce \ref{twofam}, which we restate:
\begin{thm}\label{twofam2}
Let $G$ be a bipartite graph with bipartition $(A,B)$, and let $A$ be the union of disjoint sets $A_1,A_2$. Suppose that
$A_1,A_2$ are both laminar in $B$. Then the following are equivalent:
\begin{itemize}
\item every hole of $G$ has length four;
\item there is a tree $T$ with $V(T)=B$, such that for each $a\in A$, $N(a)$ is the vertex set of a subtree of $T$.
\end{itemize}
\end{thm}
\Proof
Let $T$ be a tree with $V(T)=B$, such that for each $a\in A$, $N(a)$ is the vertex set of a subtree of $T$.
We claim that the first bullet holds. Suppose that $C$ is a hole in $G$ of length at least six. By \ref{subtrees}, some
vertex $b\in B$ has at least three neighbours in $V(C)$; and so we may assume that at least two of them belong to $A_1$, say.
Let $a,a'\in V(C)\cap A_1$ both be adjacent to $b$. Since $A_1$ is laminar, and $a,a'$ have a common neighbour in $B$,
we may assume that $N(a)\subseteq N(a')$. Consequently both neighbours of $a$ in $C$ are adjacent to $a'$, a contradiction
since $C$ is a hole of length at least six.

For the converse, assume the first bullet of \ref{twofam2} holds. Then trivially the first bullet of \ref{subtrees} holds, so by \ref{subtrees}
the second bullet of \ref{twofam2} holds. This proves \ref{twofam2}.~\bbox
\section{The case when $\ell$ is even}\label{sec:evencase}

In this section we complete the proof of \ref{mainthm} when $\ell$ is even.
In view of \ref{easymainthm} it suffices to prove the following:
\begin{thm}\label{maineven}
Let $\ell\ge 8$ be even, and let $G$ be an $\ell$-holed graph with no clique cutset or universal vertex.
Then $G$
is either a blow-up of a cycle of length $\ell$, or a blow-up of an $\ell$-framework.
\end{thm}
\Proof
By \ref{partition}, we may assume that $G$ is a bordered blow-up of an $\ell$-frame. Let $G= H\cup J\cup K$, where:
\begin{itemize}
\item $H$ is a blow-up of an $\ell$-frame $F$, where $F$ has sides $A,B$;
\item $(J, A, (W_t:t\in V(A)))$ and
$(K,B,(W_t:t\in V(B)))$ are borders; and
\item $V(H\cap J) = W(A)$, and $V(H\cap K)=W(B)$, and $V(J),V(K)$ are disjoint and anticomplete.
\end{itemize}
Let the bars of $F$ be $P_1\LL P_k$ where $P_1\LL P_m$ have length $\ell/2-1$, and $P_{m+1}\LL P_k$ have length $\ell/2-2$, and $m\le k-2$.
Let each bar $P_i$ have ends $a_i\in V(A)$ and $b_i\in V(B)$.
For $1\le i\le k$, we say that $a_i$ is the {\em mate} of $b_i$ and vice versa.
Let $M=\{1\LL m\}$ and $N=\{m+1\LL k\}$, and $A_M$ be the set $\{a_i:i\in M\}$, and define $A_N, B_M, B_N$ similarly. 
Thus $A_N, B_N$ are stable sets of $F$,
$A_M, B_M$ are cliques of $F$, and for $i\in N$ and $j\in M$, $a_i, a_j$ are $F$-adjacent if and only if $b_i, b_j$ are not $F$-adjacent.
From the definition of an $\ell$-frame, it follows that some vertex in $A_N$ has no $F$-neighbour in $V(A)$.

From \ref{spine}, there is a basis $T$ for $(J, A, (W_t:t\in V(A)))$ with an $(A_M, A_N)$-spine $R$, and a basis $S$ for 
$(K,B,(W_t:t\in V(B)))$
with a $(B_M, B_N)$-spine $R'$.
We will show that $S\cup T\cup P_1\cupcup P_k$ is an $\ell$-framework, and
$G$ is a blow-up of it. 

Since neither of $A,B$
is connected, it follows that $r(T)\notin V(A)$, and $r(S)\notin V(B)$. Let 
$a_0, a_1\LL a_m$ be the vertices of $V(A)\cup \{r(T)\}$ in order in $R$, where $a_0=r(T)$.
For $0\le i\le m$ let $T_i$ be the \sub{} of $T$ induced on the set of vertices $v$ of $T$ such that there is a directed path
of $T$ from $a_i$ to $v$, not using $a_{i+1}$ if $i<m$.
Thus the sets 
$$V(A)\cap V(T_0), V(A)\cap V(T_1)\LL V(A)\cap V(T_m)$$ 
partition $V(A)$.
Let $V(A)\setminus V(R)=\{a_{m+1}\LL a_k\}$, numbered so that for $0\le i<i'\le m$, if $a_j\in V(T_i)$ and $a_{j'}\in V(T_{i'})$ then $j>j'$.
Let $V(B)\cap V(R')=\{c_1\LL c_m\}$ where $r(S)=c_0, c_1\LL c_m$ are in order in $R'$. Since $R'$ is a $(B_M, B_N)$-spine, it follows
that $\{c_1\LL c_m\}=\{b_1\LL b_m\}$. For $0\le i\le m$ let $S_i$ be the \sub{} of $S$ induced on the set of vertices $v$ of $S$ 
such that there is a directed path
of $S$ from $c_i$ to $v$, not using $c_{i+1}$ if $i<m$.
\\
\\
(1) {\em $T_0, T_m$ both have at least one leaf, and so the last vertex of $R$ is not a leaf of $T$. Similarly, $S_0, S_m$ both have a leaf, and so the last vertex of $R'$ is not a leaf of $S$.}
\\
\\
From the definition of an $\ell$-frame, there exists $b_i\in B_N$ with no neighbours in $B_M$, and hence $a_i\in A_N$
is complete to $A_M$, and in particular is $J$-adjacent to $a_m$; and consequently $a_i\in T_m$. Thus 
$T_m$ has a leaf. Since $a_0\notin V(A)$, it follows from the definition of a basis that the $T$-inneighbour of $a_1$
has outdegree at least two, and so $T_0$ has a leaf.
The second claim follows similarly.
This proves (1).

\bigskip

Let $0=i_0<i_1<\cdots < i_p=m$ be the values of $i\in \{0\LL m\}$ with $L(T_i)\ne \emptyset$, and let
$0=i_0'<i_1'<\cdots < i_{p'}'=m$ be the values of $i\in \{0\LL m\}$ with $L(S_i)\ne \emptyset$.
\\
\\
(2) {\em $p=p'$; and for $0\le q\le p$, the vertices in $L(T_{i_q})\cap V(A)$ are the mates of the vertices in $L(S_{i_{p-q}'})\cap V(B)$.}
\\
\\
Let us say two vertices in $L(T)=A_N$ are {\em equivalent} if they have the same $J$-neighbours in $\{a_0\LL a_m\}$. The 
equivalence classes are the nonempty sets of the form $L(T_i)\cap V(A)$ for some $i\in \{0\LL m\}$, that is, the sets 
$$L(T_{i_0})\cap V(A), L(T_{i_1})\cap V(A)\LL L(T_{i_p})\cap V(A).$$ 
But two vertices $a_i, a_j\in A_N$ have the same neighbours in $A_M$
if and only if $b_i,b_j$ have the same neighbours in $B_M$; and consequently $p=p'$, and there is a bijection $\phi$
from $\{i_0\LL i_p\}$ onto $\{i_0'\LL i_p'\}$ such that for $0\le q\le p$, $L(S_{\phi(i_q)})\cap V(B)$ is the set of mates of the 
vertices in 
$L(T_{i_q})\cap V(A)$.

Suppose that there exist $q,q'$ with $0\le q<q'\le p$ such that $\phi(i_q)<\phi(i_{q'})$. Let $a_j\in L(T_{i_q})$, and so
$b_j\in L(S_{\phi(i_q)})$. Since $i_q<i_{q'}$, $a_j$ is nonadjacent to $a_{i_{q'}}$. But since $\phi(i_q)<\phi(i_{q'})$,
$b_j$ is nonadjacent to $b_{\phi(i_{q'})}$, a contradiction. Thus there are no such $q,q'$; and so $\phi(i_q)=i'_{p-q}$
for $0\le q\le p$. This proves (2).
\\
\\
(3) {\em $i'_{p-q}=m-i_{q}$ for $0\le q\le p$.}
\\
\\
Let $0\le q<p$. The vertices in $\{a_1\LL a_m\}$ that are complete to $L(T_{i_{q+1}})\cap V(A)$ and anticomplete to 
$L(T_{i_q})\cap V(A)$ are the vertices $a_j$ for $j\in \{i_q+1\LL i_{q+1}\}$; and so by (2), their mates are the vertices in $\{c_1\LL c_m\}$
that are anticomplete to $L(S_{i_{p-q-1}'})\cap V(B)$ and complete to $L(S_{i_{p-q}'})\cap V(A)$, that is, the vertices
$c_{i_{p-q-1}'+1}\LL c_{i_{p-q}'}$.
Consequently $i_{q+1}-i_q= i_{p-q}'-i_{p-q-1}'$, that is, $i_{q+1}+i_{p-q-1}'= i_q+i_{p-q}'$, for $0\le q<p$.
Since $i_0=0$ and $i'_p=m$, we deduce that $ i_q+i_{p-q}'=m$ for $0\le q\le p$. By induction on $q$ it follows that 
$i_{p-q}'=m-i_q$ for $0\le q\le p$.
This proves (3).

\bigskip

Let $1\le q< p$. As in (3), the mates of the vertices $a_{i_{q}+1}\LL a_{i_{q+1}}$ 
are the vertices 
$$c_{i_{p-q-1}'+1}\LL c_{i_{p-q}'}$$ 
that is, $c_{m-i_{q+1}+1}\LL c_{m-i_{q}}$. 
But $c_{m-i_{q+1}+1}\LL c_{m-i_{q}}$ are all adjacent $K$-twins; and so we may choose $S, c_0\LL c_m$
such that $b_j=c_{m-j+1}$ for $i_q<j\le i_{q+1}$. By repeating this for all values of $q$ we may therefore assume that
$b_j=c_{m-j+1}$ for $0<j\le m$.

In summary, 
there exist $0=i_0<i_1<\cdots <i_p=m$ with the following properties:
\begin{itemize}
\item For $0\le i\le m$ let $T_i$ be the \sub{} of $T$ induced on the set of vertices $v$ of $T$ such that there is a directed path
	of $T$ from $a_i$ to $v$, not using $a_{i+1}$ if $i<m$. Then $|T_i|>1$ if and only if $i\in \{i_0,i_1\LL i_m\}$.
\item Let $S_0$ be the \sub{} of $S$ induced on the set of vertices $v$ of $S$ such that there is a directed path
of $S$ from $b_0$ to $v$, not using $b_m$ if $m>0$.
For $1\le i\le m$ let $S_i$ be the \sub{} of $S$ induced on the set of vertices $v$ of $S$ such that there is a directed path
		of $S$ from $b_i=c_{m-i+1}$ to $v$, not using $b_{i-1}=c_{m-i}$ if $i>0$.
Then $|S_i|>1$ if and only if $i\in \{0,i_0+1,i_1+1\LL i_{m-1}+1\}.$
\item For $0\le q< p$, the leaves of $T_{i_q}$ in $V(A)$ are the mates of the leaves of $S_{i_q+1}$ in $V(B)$, 
and the leaves of $T_{i_p}$ in $V(A)$ are the mates of the leaves of $S_0$ in $V(B)$.
\end{itemize}

To complete the proof of \ref{maineven}, it suffices to check the condition about being coarboreal.
Let $(x,y)$ be one of the pairs $(i_0, i_0+1)\LL (i_{p-1}, i_{p-1}+1), (i_p, 0)$. Thus there are \arbs{} $T_x, S_y$
defined, and the set of leaves of $T_x$ in $V(A)$ are the mates of the leaves of $S_y$ in $V(B)$. 
We must show that $T_x, S_y$ are coarboreal under the bijection that sends each vertex in $L(T_x)\cap V(A)$ to its mate.
Let the leaves of $T_x$ be $\{a_r\LL a_s\}$; so the leaves of $S_y$ are $\{b_r\LL b_s\}$. 

Let $X=V(T_x)\setminus L(T_x)$ and $Y=V(T_y)\setminus L(T_y)$.
Make a bipartite graph $D$ with bipartition $(X\cup Y, \{r\LL s\})$,
in which each vertex $v\in X$ is adjacent to 
$i\in \{r\LL  s\}$ if and only if $v$ is $J$-adjacent to $a_i$, and 
each vertex $v\in Y$ is adjacent to   
$i\in \{r\LL  s\}$ if and only if $v$ is $K$-adjacent to $b_i$. 
Then $X,Y$ are both laminar in $D$, so we can apply \ref{twofam}. If the second bullet of \ref{twofam} holds, then
$T_x, S_y$ are coarboreal as required, so we assume that the first bullet of \ref{twofam} does not hold, and so there is a hole $C$ of $D$ with length more than four.

If $u,v\in V(C)\cap (X\cup Y)$, with distance two in $C$, then each of them has a neighbour in $C$ not adjacent to 
the other; and since they have a common neighbour in $C$, and $X, Y$ are both laminar in $D$, it follows that one of $u,v$
is in $X$ and the other in $Y$. Hence the length of $C$ is a multiple of four, and we may label its vertices in order
as 
$$v_1\DD j_1\DD v_2\DD j_2\CC v_{t}\DD j_{t}\DD v_1$$ 
where $t\ge 4$ is even, and $v_1,v_3,v_5,\ldots \in X$, and 
$v_2,v_4,v_6,\ldots\in Y$, and $j_1,j_2\LL j_{t}\in \{r\LL s\}$. The set $\{v_1,v_3,v_5,\ldots\}$
is laminar in $D$, and so no two of its members have a common neighbour in $\{r\LL s\}$; and consequently $v_1,v_3,v_5,\ldots$ 
are pairwise not $J$-adjacent, and similarly $v_2,v_4,v_6,\ldots$ are pairwise not $K$-adjacent. We recall that for $1\le i\le k$, $P_i$ is the bar of the $\ell$-frame $F$
with ends $a_i,b_i$. Then 
$$v_1\DD a_{j_1}\DD P_{j_1}\DD b_{j_1}\DD v_2\DD b_{j_2}\DD P_{j_2}\DD a_{j_2}\DD v_3\CC a_{j_{t}}\DD v_1$$
is a hole of $G$. Since each $P_i$ has length $\ell/2-2$, this hole has length $(\ell/2) t$, a contradiction. This proves that
$T_x, S_y$ are coarboreal under the bijection that sends each vertex in $L(T_x)\cap V(A)$ to its mate. Consequently this
proves \ref{maineven}, and so completes the proof of \ref{mainthm}.~\bbox

\section*{Acknowledgement}
We would like to express our thanks to the referee for a report that was very thorough and very helpful.

\end{document}